\documentclass [11pt]{article}
\usepackage[utf8]{inputenc}

\usepackage{amsmath}
\usepackage{ulem}
\usepackage{amsfonts}
\usepackage{amsthm}
\usepackage{enumerate}
\usepackage[dvipsnames]{xcolor}
\usepackage{amssymb}
\usepackage{graphicx}

\def\ord{{\rm ord\,}}

\def\videbox{\mathbin{\vbox{\hrule\hbox{\vrule height1ex \kern.5em\vrule
height1ex}\hrule}}}
\usepackage[a4paper]{geometry}

\newcommand{\bB}{\mathbb {B}}% relatifs
\def\bC{{\mathbb{C}}}%complexes
\def\bK{{\mathbb{K}}}%reels ou complexes
\def\bE{{\mathbb{E}}}%esperance
\newcommand{\bN}{{\mathbb {N}}}
\newcommand{\bO}{{\mathbb {O}}}% octaves
\newcommand{\bH}{{\mathbb {H}}}% quaternions
\newcommand{\bR}{{\mathbb {R}}}%reels
\def\bS{{\mathbb{S}}}
\newcommand{\bZ}{\mathbb {Z}}% relatifs
\def\cC{{\cal C}}

\def\cG{{\cal G }}
\def\cL{{\cal L }}
\def\cP{{\cal P }}
\def\cF{{\cal  F}}

\def\cH{{\cal  H}}
 \def\cR{{\cal  R}}
 \def\cI{{\cal  I}}
 \def\cN{{\cal  N}}
 \def\CP{{\mathbb{CP}}}
 \def\RP{{\mathbb{RP}}}
 \def\indic{{\bf{1}}}

\def\ncCu  {1}
\def\ncPa  {2}
\def\ncPaTa{3}
\def\ncX   {4}
\def\ncXY  {5}
\def\ncXX  {6}
\def\ncXXY {7}

\def\tr{\textmd{trace}\,}

\def\Id{\textmd{Id}\,}
\def\Re{\textmd{Re}\,}
\def\Im{\textmd{Im}\,}

\def\lag{\langle}
\def\rag{\rangle}

\newcommand{\bEps}{\varepsilon}

\def\dis{\displaystyle }

\def\mult{{\rm mult}}

  \def\LL{{\bf L}}
  \def\Ga{{\bf \Gamma}}
  \def\De{{\bf \Delta}}

\def\bpm{\begin{pmatrix}}
\def\epm{\end{pmatrix}}
\def\pn{\par\noindent}

\def\uu{{\bf u}}
\def\vv{{\bf v}}
\def\ww{{\bf w}}
\def\xx{{\bf x}}
\def\yy{{\bf y}}
\def\zz{{\bf z}}

\newcommand{\footnoteremember}[2]{
\footnote{#2}
\newcounter{#1}
\setcounter{#1}{\value{footnote}}
}
\newcommand{\footnoterecall}[1]{
\footnotemark[\value{#1}]
}

\newtheorem{theoreme}{Theorem}[section]
\newtheorem{lemme}[theoreme]{Lemma}
\newtheorem{definition}[theoreme]{Definition}
\newtheorem{proposition}[theoreme]{Proposition}
\newtheorem{corollaire}[theoreme]{Corollary}

\newenvironment{exemp}{\noindent{\bf Example. --- }}{\par}
\newenvironment{exemps}{\noindent{\bf Examples}\benum}{\eenum\par}
\newtheorem{rmq}[theoreme]{Remark}
\newtheorem{rmqs}[theoreme]{Remarks}

\newenvironment{preuve}{\noindent{\it Proof --- }}
{\hfill\rule{1.3mm}{2mm}\par} 
\newenvironment{application}{\noindent{\bf Application. --- }}{\par}
\newenvironment{applications}{\noindent{\bf Applications. --- 
}\benum}{\eenum\par}

\newcommand{\bex}{\begin{exemp}}
\newcommand{\eex}{\end{exemp}}
\newcommand{\bexs}{\begin{exemps}}
\newcommand{\eexs}{\end{exemps}}

\newcommand{\bprop}{\begin{proposition}}
\newcommand{\eprop}{\end{proposition}}

\newcommand{\bthm}{\begin{theoreme}}
\newcommand{\ethm}{\end{theoreme}}

\newcommand{\bcor}{\begin{corollaire}}
\newcommand{\ecor}{\end{corollaire}}

\newcommand{\blem}{\begin{lemme}}
\newcommand{\elem}{\end{lemme}}

\newcommand{\beq}{\begin{equation}}
\newcommand{\eeq}{\end{equation}}
\newcommand{\beqna}{\begin{eqnarray}}
\newcommand{\eeqna}{\end{eqnarray}}
\newcommand{\beqnas}{\begin{eqnarray*}}
\newcommand{\eeqnas}{\end{eqnarray*}}

\newcommand{\bmat}{\begin{pmatrix}}
\newcommand{\emat}{\end{pmatrix}}

\newcommand{\bpf}{\begin{preuve}}
\newcommand{\epf}{\end{preuve}}

\newcommand{\benum}{\begin{enumerate}}
\newcommand{\eenum}{\end{enumerate}}

\newcommand{\bitem}{\begin{itemize}}
\newcommand{\eitem}{\end{itemize}}

\newcommand{\brmq}{\begin{rmq}}
\newcommand{\ermq}{\end{rmq}}

\newcommand{\brmqs}{\begin{rmqs}}
\newcommand{\ermqs}{\end{rmqs}}

\newcommand{\bapp}{\begin{application}}
\newcommand{\eapp}{\end{application}}

\newcommand{\bapps}{\begin{applications}}
\newcommand{\eapps}{\end{applications}}

\newcommand{\bdefi}{\begin{definition}}
\newcommand{\edefi}{\end{definition}}

\newcommand{\bequa}{\begin{equation}}
\newcommand{\eequa}{\end{equation}}

\theoremstyle{definition}

\def\mP{{\mathbb P}}

\def\zo#1{{\edef\a{#1}\edef\z{0}\if\a\z{$\circ$}\else{$\bullet$}\fi}}
\def\NND#1#2#3#4#5#6{\hskip0pt\left[\text{\edef\a{#1}\edef\b{#2}\edef\c{#3}\edef\d{#4}\edef\e{#5}\edef\f{#6}%
\lower1mm\hbox{%
  \vbox to 4.4mm{%
    \vfill%
    \hbox{%
      \hbox to0pt{\smash{\raise3.4mm\hbox{\zo\a}}}%
      \hbox to0pt{\smash{\raise1.7mm\hbox{\zo\b\zo\c}}}%
      \zo\d\zo\e\zo\f%
    }%
    \vskip-0.8mm%
  }%
}%
}\right]}

\begin{document}

\author{D. Bakry \footnoteremember{IMT}{Institut de Math\'ematiques de Toulouse, Universit\'e Paul Sabatier, 118 
route de Narbonne, 31062 Toulouse, France;  \textsf{dominique.bakry@math.univ-toulouse.fr, stepan.orevkov@math.univ-toulouse.fr}}, S. Orevkov \footnoterecall{IMT}, M. Zani \footnote{Institut Denis Poisson, Universit\'e d'Orl\'eans, Universit\'e de Tours, CNRS, Route de Chartres, B.P. 6759, 45067, Orl\'eans cedex 2, France ; \textsf{marguerite.zani\symbol{64}univ-orleans.fr}}}
%\date{\today}

%\makeatletter \bRenewcommand{\@oddfoot}{\sl \small
 %\hfil \thepage\hfil \today}
%\bRenewcommand{\@oddhead}{\sl \small
% \hfil preprint under construction}
 %\makeatother

\title{Orthogonal polynomials and diffusion operators
}
\maketitle
\begin{abstract}
We study the following problem: describe the triplets $(\Omega,g,\mu)$ where
$g= (g^{ij}(x))$ is the (co)metric associated with the symmetric second order
differential operator $\LL (f) = \frac{1}{\rho}\sum_{ij} \partial_i (g^{ij} \rho
\partial_j f)$ defined on a domain $\Omega$ of $\mathbb R^d$  (that is $\LL$
is a diffusion operator with reversible measure $\mu(dx)= \rho(x) dx$)  and such
that there exists an orthonormal    basis  of $\cL^2(\mu)$ made of polynomials
which are at the same time eigenvectors of $\LL$, where the polynomials are
ranked according to their natural degree.
We reduce this problem to a certain algebraic problem (for any $d$) and we find all
solutions for $d=2$ when $\Omega$ is compact.
Namely, in dimension $d=2$, and up to affine transformations, we find $10$ compact
domains $\Omega$ plus a one-parameter family. The proof that this list is exhaustive
relies on the Pl\"ucker-like formulas for the projective dual curves applied to the
complexification of $\partial\Omega$. We then describe some geometric origins for
these various models. We also give some description of the non-compact cases in this
dimension.
\end{abstract}
\renewcommand{\abstractname}{R\'esum\'e}
\abstract{Nous considérons le problème suivant: décrire les triplets $(\Omega,g,\mu)$ où $g= (g^{ij}(x))$ 
est la (co)métrique associée à l'opérateur différentiel du second ordre symétrique $\LL (f) = \frac{1}{\rho}\sum_{ij} \partial_i (g^{ij} \rho
\partial_j f)$ défini sur un domaine $\Omega$ de $\mathbb R^d$  ( i.e. $\LL$ est un opérateur de diffusion de mesure réversible
$\mu(dx)= \rho(x) dx$)  et tels qu'il existe une base orthonormale de polynômes de $\cL^2(\mu)$ qui sont également vecteurs propres de
$\LL$, les polynômes étant classés par ordre croissant de leur degré naturel. Nous réduisons ce problème à un problème algébrique (pour tout $d$)
et décrivons les solutions pour $d=2$ et $\Omega$ compact. Nous montrons que pour 
$d=2$, et à transformations affines près, il y a $10$ domaines compacts
$\Omega$ et une famille à un paramètre. La preuve de l'exhaustivité de ce classement repose sur des formules de type Pl\"ucker pour les courbes duales projectives appliquées à la complexification de $\partial\Omega$. 
Nous présentons alors une interprétation géométrique de ces différents modèles. Nous donnons également une description des cas non--compacts en dimension $d=2$.
}

{\small
\tableofcontents}

    \section{Introduction\label{intro}}

    \subsection{ Content of the paper}
    In this paper, we investigate the following question: for a given set $\Omega\subset \bR^d$, when does there exist 
 a probability measure  $\mu(dx)$ on $\Omega$, absolutely continuous with respect to the Lebesgue measure, and an elliptic  diffusion operator 
    $$\LL (f) = \sum_{ij} g^{ij}(x) \partial_{ij} f +\sum_i b^i(x) \partial_i f,$$ defined on $\Omega$ such that there exists   an orthonormal basis for $\cL^2(\mu)$,  formed by orthogonal polynomials  ordered according to the total degree
\footnote{This means that the space of polynomials of degree $\le n$
          is $\LL$-invariant for any $n$},
which are  eigenvectors of the operator $\LL$. Moreover, can we describe the sets, the operators and the measures?
    
     In dimension 1, given the measure $\mu$, there is a unique family of associated orthogonal polynomials, up to a choice of sign. Some of them share extra properties, and as such are widely used in many areas. This is in particular the case of  Hermite, Laguerre, and Jacobi polynomials, which correspond respectively to the measures with density $Ce^{-x^2/2}$ on $\bR$, $C_ax^{a-1} e^{-x}$, $a>0$, on $[0, \infty)$ and $C_{a,b}(1-x)^{a-1}(1+x)^{b-1})$, $a,b>0$,  on $[-1, 1]$ (where $C, C_a, C_{a,b}$ are normalizing constants which play no role here). In those three cases,  and only in those ones, the associated polynomials are eigenvectors of some second order differential operator $\LL$:  see \cite{Bochner1929, BakM, Maze}.
      Those families have been extensively studied, since they play a central role in probability, analysis, partial differential equations, geometry, mathematical physics... (see e.g. \cite{Forn,Wil,guo1,guo2,guo3,verk1,verk2, szego1975,andrew.askey.roy.99},  see also \cite{hall, privault} and references therein).

  \par     The differential operator $\LL$ may be replaced by some other generator of  a Markov semigroup (finite difference, or $q$-difference operators)
and the orthogonal polynomial eigenfunctions are Hahn, Krawtchouk, Charlier, Meixner (see \cite{NSU}).
In dimension 1, a classification had been done for such families, see~\cite{hahn1949, Vinet.Zhed.2008}, but there are very few such classification results beyond the dimension 1 case.

The main motivation for this study lies  in probability theory, where such models for diffusion operators are the easiest ones where one may check various quantities relating properties of the generator (curvature, diameter, spectral gap, etc.) to the best possible estimation for the various constants in functional inequalities (logarithmic Sobolev inequalities, Sobolev inequalities, isoperimetric inequalities,  estimates on the heat kernel, e.g.). It turns out that  the  dimension 1 models, where most of the computations may be done explicitly, provide good  models for testing various conjectures. However, there are too few dimension 1 models to really explore all the various questions arising in this area. It seems therefore natural to try to describe more families where such computations may be made. Beyond this, those families provide natural bases into which computations may be made in approximation theory, partial differential equations, etc.
 
  The aim of this paper is then to extend the dimension 1 classification  for differential operators to higher dimensions, and in particular in  dimension 2, to give a precise description of the differential operators, the measures and the domains concerned.

  In $\bR^d$, in order to properly define an orthonormal polynomial basis, we first have to agree on a way of ordering the polynomials, and this is done according to the choice of a degree.  Choosing some positive integers $w_1, \dots, w_d$, a monomial $x_1^{p_1}x_2^{p_2}\cdots x_d^{p_d}$ will  have a degree $
  w_1p_1+ \cdots +w_dp_d$, and the degree of a polynomial is the maximum degree of its monomials (we may of course reduce to the case where those integers $w _i$ have no common factor).  When all the $w_i$ are equal to 1, this is the usual degree. According to this, one defines the finite dimensional  vector space $\cP_n^d$ of polynomials with total degree less than $n$, and a polynomial  orthogonal basis is defined by the choice for each $n$ of an orthonormal basis of the orthogonal complement of $\cP_{n-1}^d$ in $\cP_n^d$. 
  
  Although many of the results of this paper, in particular in section~\ref{cas.general}, could be extended to the general degree case, we stick in this paper to the usual degree.
  
  Given the choice of the degree, for bounded sets $\Omega \subset \bR^d$ with piecewise smooth boundary, one may reduce the problem to some algebraic question  about the boundary.   
   In dimension 2, and for the usual degree, this problem may be  entirely solved (Theorem~\ref{thm.central}): we provide the complete list of $10$ different  bounded sets $\Omega\subset \bR^2$ together  with a one parameter family of domains (coaxial parabolas) which, up to affine transformations,     are the only ones on which this problem have a solution.  We also provide  in Section~\ref{comments.models} a  complete description of the associated measures and operators.  Under stronger requirements  on the sets, we also provide a list of the 7 non compact models which solve the problem in dimension 2. Let us mention that this choice of natural degree is not done for simplicity. There are many other bounded models in dimension $2$ with associated orthogonal polynomials according to other choices for the  degree, but the techniques developed below for classification  may not be easily adapted the the general situation.  In particular, in dimension 2, one may construct orthogonal polynomials from root systems (Heckman-Opdam polynomials,  see \cite{Heckman87, Heck1, HeckOp,Heckman1991,opdam})  or finite subgroups of $O(3)$ ( see \cite{b.meyer54, Dunkl1988,bakry-bressaud} for the construction of such orthogonal polynomial families). Indeed,  we recover in our list  the Hekman-Opdam  polynomials associated with the root systems  $B_2$ (Section~\ref{sec.B2}) and $A_2$ (Section~\ref{sec.A2}), but not the family associated with $G_2$, which corresponds to a degree of $x^py^q$ equal to $2p+q$ (see Section~\ref{sec.A2} for details). Many other models in dimension $2$ arising from finite subgroups of $O(3)$ do not appear  either in our  classification, due again to another   degree in the choice of the degree of the  polynomials.   However, even with the usual degree, the example of Section~\ref{sec:DlePtCub}  shows that root systems and finite subgroups of $O(3)$ do not provide all the possible models.

   Further extension to higher dimensional models are  also given, although a classification seems out of reach with the methods of the 2-dimensional analysis, even with the usual degree.
     \par
 
 \subsection{The general problem\label{sec.gal.pb}}
   
    Orthogonal polynomials are a long standing subject of investigation in mathematics. They yield natural Hilbert bases in $\cL^2(\mu)$ spaces, where $\mu$ is a probability measure on some measurable  set $\Omega$ in $\bR^d$ for which polynomials are dense. As a way to describe functions  $f:\Omega \mapsto \bR$, they are  used in many problems in analysis, for example in  partial differential equations, especially when they present some quadratic nonlinearities: since products are in general easy to compute in such polynomial  bases, approximation schemes which consist in restricting the approximation of functions to a finite number of components in those bases are easy to implement in practice.

    In higher dimension, there are several choices for a basis of orthogonal polynomials, and no canonical choice may be proposed in general.  However, many families have been described in various settings. In particular, multivariate analogues of the classical families,  in particular  those which are eigenvectors   of differential operators,  have been put forward by many authors:  see~\cite{Fernandez.al.2005, Heck1,Heckman1991,Koorn, Koorn1,Koorn2,Koorn3,Koorn4,KoornSchw, KrallS, Rosler1}); see also \cite{alvarezco} for a generalization of the Rodrigues formula.  For a general overview on orthogonal polynomials of several variables, we refer to  Suetin \cite{suetin} and to the  book of Dunkl and Xu \cite{DX}. 
  
    As mentioned above, in dimension $d\geq 2$, one orders  in general polynomials by their total degree:  if $\cP_n^d$  denotes the set of polynomials in $d$ variables of degree not greater than $n$,  we are looking for a Hilbert basis of $\cL^2(\mu)$ such that for each $n$, we get a  finite-dimensional basis of $\cP_n^d$. This basis is not unique in general.  This is what we call a polynomial  orthogonal   basis, and is the object of our study.  As already mentioned,  we stick in this paper with the natural degree, but most of the general considerations developed in Section~\ref{section2} remain valid in the general case.

\par
    
    On the other hand, these polynomial bases are   not always the best choice to expand functions or to obtain good approximation schemes. This is in particular the case in probability theory, when one is concerned with symmetric diffusion processes as they naturally appear as solutions of stochastic differential equations.
    Indeed, a Markov diffusion process  $(X_t)_{t\geq 0}$, with continuous trajectories on an open set of $\bR^d$ or a manifold, has a law entirely characterized  by the family of Markov kernels $(P_t)_{t\geq 0}$:
$$
       P_t(f)(x)= \bE(f(X_t)/X_0=x)\,,\quad x\in\bR^d\,,
$$
    where $f$ is in a suitable class of functions.
    The infinitesimal generator $\LL$ associated with $(P_t)_{t\geq 0}$ is defined by
    $$\LL f=\lim_{t\to 0}\frac{P_tf-f}{t}\,,$$
    whenever this limit exists.   
    
    This operator governs the semigroup in the sense that if
    $F(x,t) = P_t(f)(x)$, then $F$ is the solution of the heat equation
    $$\partial_t F= \LL F, ~F(x,0)= f(x).$$
    It is quite difficult in general to obtain a complete description of $P_t$ in terms of the operator $\LL$, which is in general the only datum that one has at hand from the description of $(X_t)$, for example as the solution of a stochastic differential equation.  This operator $\LL$ is a second order differential operator with no zero order component, moreover semi-elliptic, of the form
     \begin{equation}\label{formgen}\LL(f)=\sum_{ij}g^{ij}(x)\partial _{ij}f+\sum_ib^i(x)\partial_if\,.
 \end{equation}

    Although not easy to compute explicitly,  the operator $P_t$, which describes the law of the random variable $X_t$, has a nice expression at least when $\LL$ is self-adjoint with respect to some measure $\mu$ ($\mu$ is then said to be the reversible measure for $(X_t)$), and when the spectrum is discrete . When $\mu$ has a density $\rho$ which is $\cC^1$ with respect to the Lebesgue measure, and if the coefficients $g^{ij}$ are also assumed to be at least $\cC^1$, then this latter case amounts to look for operators $\LL$ of the form
    \beq\label{sym.diffusion}\LL (f) = \frac{1}{\rho}\sum_{ij} \partial_i (g^{ij} \rho\, \partial_j f).\eeq 
    In this paper, we shall restrict our attention to operators which are elliptic in the interior of the support of $\mu$.  Such an operator described in~\eqref{sym.diffusion} will be called a symmetric diffusion operator.
Notice however that the ellipticity assumption is never used in the paper and all
our results remain true for any non-degenerate (co)metric $(g^{ij})$.
Moreover, in dimension 2, where we give a complete classification,
we see {\it a posteriori\/} that $\LL$ appears to be elliptic
(without this {\it a priori\/} assumption)
each time when $g^{ij}$ is unique up to scalar factor.

      \par 
     In the case  under study,  the spectral decomposition leads to some more or less explicit representation. Namely, if there is an orthonormal basis $(e_n)$ of $\cL^2(\mu)$  composed of eigenvectors of $\LL$,
    $$\LL e_n= -\lambda_n e_n,$$ then one has
    $$P_t(f)(x) = \int f(y) p_t(x,y) d\mu(y),$$ where 
    $$p_t(x,y) = \sum_n e^{-\lambda_n t} e_n(x)e_n(y).$$
   For fixed $x$,  the function $p_t(x,y)$ represents the density with respect to $\mu(dy)$ of the law of $X_t$ when $X_0=x$.
    Of course, this representation is a bit formal, since one has to insure first that this series converges, which requires $P_t$ to be trace class, or Hilbert-Schmidt. However, even if it is quite rare that the eigenvalues $\lambda_n$ and the eigenvectors $e_n$ are explicitly known, it can be of great help to know that such a decomposition exists: it provides a good approximation of $P_t$ when $t$ goes to infinity, and as such allows to control convergence to equilibrium.  But even when one explicitly knows the eigenvectors and eigenvalues, it is not always easy to extract many useful information from the previous description. It is even not immediate to check in general that the previous expansion leads to  nonnegative functions.   
     \par

    Even when $\LL$ is elliptic and symmetric, its knowledge, given on say smooth function compactly supported in $\Omega$, is not enough to describe the associated semigroup $P_t$ or any self-adjoint extension of $\LL$.
    One requires in general some boundary conditions.  This requirement will be useless in our context, since we shall impose the  eigenvectors to be polynomials. As a counterpart, this will impose some boundary condition on the operator itself.

    As mentioned earlier, we are interested in the description of  the situation when the   eigenvector expansion coincides with a family of orthogonal polynomials associated with the reversible measure. Although the situation  is well known and described in dimension $1$,  such description is not known in higher dimension, apart from some generic families. At least when the set $\Omega$ is relatively compact with $\cC^1$ piecewise boundary, and when the reversible measure $\mu$ has a $\cC^1$ density with respect to the Lebesgue measure,  we may turn the complete description of this situation into a problem of algebraic nature: the operators and the measures
can be completely recovered from
%are entirely described by
the boundary of $\Omega$, which is some algebraic surface of degree at most $2d$ in dimension $d$. Then, we completely solve this problem in dimension 2, leading, up to affine transformations,  to the 11 different possible boundaries: 
   the square, the circle, the triangle, the coaxial parabolas, the parabola with one tangent and one secant, the parabola with two tangents, the nodal cubic, the cuspidal cubic with one secant line, the cuspidal cubic with one tangent, the swallow tail and  the deltoid.

\par

   \par     
    Once the boundary is known, the possible measures are completely described.
They depend on some parameters (as many parameters as irreducible components in
the minimal equation of the boundary of $\Omega$). It turns out that in many situations,
for some half integer values of these parameters, the associated operator has a natural
geometric  interpretation in terms of Lie group action on symmetric  spaces.
We then provide explicitly many of these interpretations whenever they are at hand.  
  
\par 
 We also show that when  $\Omega=\bR^2$ (that is when the density $\rho$ of $\mu$ is everywhere positive), the only possible measures are Gaussian. Under some extra hypothesis, we also provide some classification of the non compact models. Further  extensions to higher dimension are also provided.

\par The paper is organized as follows. In Section~\ref{section2}, after some rapid overview of the dimension 1 case, we describe the general setting in any dimension, and, when the set $\Omega$ is relatively compact with piecewise smooth boundary,  we show  how to reduce  the description to the classification of some algebraic surfaces in $\bR^d$.   We also describe the various associated measures from the description of the boundary of
$\Omega$. 

Then, Section~\ref{dim2} is devoted to the classification of the compact dimension 2 models, which leads to 11 different cases up to affine transformations.  Section~\ref{comments.models} provides a more detailed description of the 11 models, with some insight on their geometric content for various values of the parameters.
Section~\ref{no.boundary} describes the case where no boundary is present, and the main result of this section is  that the only possible measures are Gaussian ones.
Section~\ref{non.compact.with.boundary} describes the non compact cases under some extra assumption which extends the natural condition of the compact case. Finally, Section~\ref{2fold.covering} provides some way of constructing 3-dimensional models from 2-dimensional ones.

%%%%%%%%%%%%%%%%%%%%%%%%%%%%%%%%%%%%%%%%%%%%%%%%%%%%%%%%%%%%%%%%%%%%%%%%%%%%
%%%%%%%%%%%%%%%%%%%%%%%%%%%%%%%%%%%%%%%%%%%%%%%%%%%%%%%%%%%%%%%%%%%%%%%%%%%%
%%%%%%%%%%%%%%%%%%%%%%%%%%%%%%%%%%%%%%%%%%%%%%%%%%%%%%%%%%%%%%%%%%%%%%%%%%%%
%%%%%%%%%%%%%%%%%%%%%%%%%%%%%%%%%%%%%%%%%%%%%%%%%%%%%%%%%%%%%%%%%%%%%%%%%%%%
%%%%%%%%%%%%%%%%%%%%%%%%%%%%%%%%%%%%%%%%%%%%%%%%%%%%%%%%%%%%%%%%%%%%%%%%%%%%

 \section{Diffusions associated with orthogonal polynomials\label{section2}}
 \subsection{Dimension 1\label{dim1}}
 As mentioned previously, the one-dimensional case has been completely described
for a long time (see e.g.  \cite{BakM,Bochner1929, Maze}). We recall here briefly
the framework and results.
\par
Let $\mu$ be a finite measure absolutely continuous with respect to the Lebesgue measure on an open  interval $I$ of $\bR$ with $\cC^1$ density $\rho$ (we may assume $\mu$ is a probability measure), for which polynomials are dense in $\cL^2(\mu)$ (this is automatic when $I$ is bounded, but in general it is enough to demand
that $\int \exp( \epsilon |x|) d\mu < \infty$ for some $\epsilon>0$, see~\cite{Berg.Christen.1981,Dvurecen.al.2002}).  Let $(Q_n)_{n\geq 0}$ be the family of orthogonal polynomials obtained from the sequence $(x^n)_{n\geq 0}$ by  orthonormalization, e.g. by the Gram--Schmidt process (the normalization of $Q_n$ plays no role in what follows).  Assume furthermore that some  elliptic diffusion operator $\LL$ of type ~\eqref{sym.diffusion} exists on $I$ (and therefore $\mu(dx)= \rho(x) \,dx$ is its reversible measure, that is   $\LL$ is symmetric in $\cL^2(\mu)$, at least on the set of  smooth  compactly supported functions),  such that for some sequence $(\lambda_n)$ of real numbers, 
\[\LL Q_n = -\lambda_n Q_n\,.\]
\par
Then up to affine transformations, $I$, $\mu$ and $\LL$ may be reduced to one of the three following cases:
\begin{enumerate}
\item The Ornstein--Uhlenbeck operator on 
$I=\bR$ $$H=\frac{d^2}{dx^2}-x\frac{d}{dx}\,,$$ the measure $\mu$ is Gaussian centered:
$\mu(dx)=\frac{e^{-x^2/2}}{\sqrt{2\pi}}\,dx.$
The family $(Q_n)_n$ are the Hermite polynomials,  denoted $H_n(x)$ or $H_n(x/\sqrt{2})$ in the literature, and  $\lambda_n=n$.
\item The Laguerre operator (or squared radial generalized Ornstein--Uhlenbeck operator) on
$I=\bR_+^*$ $$L_a=x\frac{d^2}{dx^2}+ (a-x)\frac{d}{dx}\,,\quad a>0,$$ the measure $\mu_a(dx)=C_ax^{a-1}e^{-x}\,dx.$
The family $(Q_n)_n$ are the Laguerre polynomials, often denoted $L_n^{a-1}(x)$, and  $\lambda_n= n$.
\item The Jacobi operator on $I=(-1,1)$
$$
   J_{a,b}=(1-x^2)\frac{d^2}{dx^2}- \big(a(x+1)+b(x-1)\big)\frac{d}{dx}\,,\quad a,b>0,
$$
the measure 
$\mu_{a,b}(dx)=C_{a,b}(1-x)^{a-1}(1+x)^{b-1}\,dx$,
the family $(Q_n)_n$ are the Jacobi polynomials, often denoted $P_n^{a-1,b-1}(x)$, and $\lambda_n = n(n+a+b-1)$.
\end{enumerate}
  The first two families appear as limits of the Jacobi case. For example, when we chose $a=b$ and let then $a$ go to $\infty$, and scale the space variable $x$ into $x/\sqrt{a}$, the measure $\mu_{a,a}$ converges to the Gauss measure, the Jacobi polynomials converge to the Hermite ones, and $\frac{2}{a} J_{a,a}$ converges to $H$.   
  
  In the same way, the Laguerre setting is obtained from the Jacobi one  fixing $b$, changing $x$ into $\frac{2x}{a}-1$, and letting $a$ go to infinity. Then $\mu_{a,b}$ converges to $\mu_b$, and $\frac{1}{a} J_{a,b}$ converges to   $L_b$. 
  
  Also, when $a$ is a half-integer, the Laguerre operator may be seen as the image of the Ornstein--Uhlenbeck operator in dimension $d$. Indeed, as the product of one-dimensional Ornstein--Uhlenbeck operators, the latter has generator 
  $H_d = \Delta -x.\nabla$. Its  reversible measure  is $e^{-|x|^2/2} dx/(2\pi)^{d/2}$,  its  eigenvectors are the products $Q_{k_1}(x_1)\cdots Q_{k_d}(x_d)$, and its  associated process  $X_t= (X^1_t, \dots, X^d_t)$, is formed of independent one dimensional Ornstein-Uhlenbeck processes, see~\cite{bglbook}. Then, if one
sets $R(x)= |x|^2$, then one may observe that, for any smooth function
$F: \bR_+\mapsto \bR$,
$$
   H_d\big(F(R)\big) = 2L_a(F)(R),
$$
where $a= d/2$.  In the  probabilist interpretation, this amounts to observe that if $X_t$ is a $d$-dimensional Ornstein--Uhlenbeck process, then $|X_{t/2}|^2$ is a Laguerre process with parameter $a= d/2$.  This coincides with the fact that the image measure of the Gaussian measure under this map is the measure $\mu_{d/2}$. 
  
  In the same way, when $a=b= d/2$, $J_{a,a}$ may be seen as the Laplace operator $\De_{S^d}$ on the unit sphere $\bS^d$ in $\bR^{d+1}$ acting on functions depending only on the first coordinate (or equivalently on functions invariant under the rotations leaving $(1, 0, \dots, 0)$ invariant), which may be interpreted as the fact that the first coordinate of a Brownian motion on the unit sphere is a diffusion process with generator $J_{d/2, d/2}$. A similar interpretation is valid for $J_{p/2, q/2}$ for some integers $p$ and $q$. Namely, let us consider functions on $\bS^{p+q-1}$ depending only on
$X= x_1^2+ \cdots x_p^2$.
Then,  setting $Y= 2X-1: \bS^{p+q-1}\mapsto [–1,1]$, for any smooth function $f: [-1,1]\mapsto \bR$,
 $\De_{\bS^{p+q-1}} f(Y)= 4 J_{q/2,p/2}(f) (Y)$.
Once again,   the associated Jacobi  process may be seen as the image of a Brownian motion on the $(p+q-1)$-dimensional sphere through the function $Y= 2X-1$.  This interpretation comes from  Zernike and Brinkman \cite{ZerBrink} and Braaksma and Meulenbeld \cite{BraakMeul} (see also  \cite{Dijksma.Koorn.1971, Koornwinder.Add.1973}).
As in the previous case,
these interpretations are compatible with the fact that the images of the uniform measure on the sphere under these various projections are the corresponding reversible measures of our operators. We shall come back to such interpretations of models as images of other ones in  paragraph~\ref{subsec.generalities}.  
  
  Let us mention that Jacobi polynomials  also play a central role in the analysis on compact Lie groups. Indeed, for $(a,b)$ taking the various values of  $(q/2, q/2)$, $((q-1)/2,1)$, $(q-1,2)$, $(2(q-1), 4)$ and $(4,8)$ the Jacobi operator $J_{a,b}$ appears as the radial part of the Laplace-Beltrami (or Casimir) operator on the compact rank 1 symmetric spaces, that is spheres, real, complex and quaternionic projective spaces, and the special case of the projective Cayley plane (see Sherman \cite{Sherman}).

%%%%%%%%%%%%%%%%%%%%%%%%%%%%%%%%%%%%%%%%%%%%%%%%%%%%%%%
%%%%%%%%%%%%%%%%%%%%%%%%%%%%%%%%%%%%%%%%%%%%%%%%%%%%%%%

 \subsection{General setting\label{cas.general}} 
  We now state our problem in full generality, and describe the framework  we are looking for.   In this section,  we describe the general problem (DOP, Definition~\ref{gal.pb.compact})  as stated above, and  we further consider a more constrained one (SDOP, Definition~\ref{def.SDOP}). It turns out that they are equivalent whenever the domain $\Omega$ is bounded, and that the latter is much  easier to handle.
To start with, we restrict the domains we are considering.
   
   \bdefi  We call a natural domain an open connected set in $\bR^d$ which is the interior of its closure.
% has a piecewise $\cC^1$ boundary and on which Stokes formula applies.
   \edefi

% Notice that 

 \bdefi  Let  $\Omega$  be a natural domain.
A diffusion operator on $\Omega$   with smooth coefficients  is a differential operator
 $\LL$, acting on smooth compactly supported function in $\Omega$, which writes 
\beq\label{defdiffusion}
  \LL(f) = \sum_{ij} g^{ij}(x) \partial_{ij}f + \sum_i b^i(x) \partial_i f,
\eeq
where $g^{ij}$ and $b^i$ are smooth functions  (that is $\cC^\infty$) on $\Omega$, and the matrix $(g^{ij})$ is symmetric, positive definite for any $x\in \Omega$.\edefi 

 The  ellipticity assumption   (i.e. the matrix $(g)$ is positive definite on $\Omega$) could be relaxed  to the weaker one of hypoellipticity. However, it would change a lot of arguments since most of the paper rely in an essential way on it.
So, the non-degeneracy of the quadratic form $(g^{ij})$ is crucial.
In contrast, as we already mentioned in the introduction, its positive definiteness
(i.e. the ellipticity of $\LL$)
is never used in the proofs (except, of course, the negativity of the eigenvalues). However, by miracle (which deserves to be explained),
our classification in dimension two
gives only elliptic solutions
when the metric is determined by $\Omega$ up to rescaling.
 Notice that
 diffusion operators (operators such that the associated semigroups are Markov operators) require at least that $\LL$ is semi-elliptic, that is
the  matrices $(g^{ij})$ are non-negative.
 
In the sequel, we shall make a constant use of the square field operator
(see \cite{bglbook})
\beq\label{defGamma}
  \Ga(f_1,f_2) = \sum_{ij} g^{ij} \partial_i f_1\partial_j f_2
         = \frac{1}{2}\Big( \LL(f_1 f_2)- f_1\LL(f_2)-f_2\LL(f_1)\Big),
\eeq
and observe that for any smooth function $\Phi: \bR^k\mapsto \bR$  and any $k$-tuple of smooth functions ${\bf f} = (f_1, \dots , f_k)$ $f_i: \Omega\mapsto \bR$, one has
\beq
  \label{chgt.de.variables} \LL\big (\Phi(f_1, \dots, f_k) \big)
  = \sum_{i,j=1}^k (\partial_{ij}\Phi)({\bf f})\Ga(f_i,f_j)+ 
    \sum_{i=1}^k (\partial_i \Phi) ({\bf f}) \LL(f_i).
\eeq

We also consider some probability measure $\mu(dx)= \rho(x) dx$ with smooth density $\rho$ on $\Omega$ for which polynomials are dense in $\cL^2(\mu)$. This last assumption is automatic as soon as $\Omega$ is relatively compact (in which case polynomials are even dense in any $\cL^p(\mu)$, $1\leq p<\infty$). It would require some extra-assumption on $\mu$ in the general case. For example, it is enough for this to hold to require that $\mu$ has some exponential moment, that is $\int_\Omega e^{\epsilon\|x\|} d\mu(x) < \infty$ for some $\epsilon >0$, in which case polynomials are also    dense in every $\cL^p(\mu)$, $1\leq p< \infty$ (see~\cite{Dvurecen.al.2002}).

The fundamental question is to study whether there exists a  polynomial orthonormal basis of $\cL^2(\mu)$, say $(P_k)$, for which  the polynomials $P_k$
are eigenvectors for $\LL$, that is that there exist some real numbers $(\lambda(P_k))$ with $ \LL P_k = -\lambda(P_k) P_k$.  Such eigenvalues $(\lambda(P_k))$ turn out to be necessarily non negative (this is a general property of symmetric diffusion operators, as a direct consequence on the non-negativity of $\Ga$).

In dimension  $d$, where $d\geq 2$, one should be precise about the notion of polynomial orthogonal basis, as mentioned in the introduction. 

 \bdefi\label{basis.order}
 Let $\Omega$ be a natural domain and $\mu$ a probability measure on $\Omega$ for which the polynomials are dense in $\cL^2(\mu)$.  Let  $\cP_n^d$ be the finite dimensional space   of polynomials with natural degree less than $n$. A polynomial orthonormal basis for $\cL^2(\mu)$ is  a choice, for each $n$, of an orthonormal basis in the orthogonal complement of  $\cP_{n-1}^d$ in $\cP_n^d$.\edefi

 As mentioned earlier, one  could consider more general situations with weighted degrees.    Although this general situation with integer weights may appear in many situations (see~\cite{bakry-bressaud} for example), our paper depends in a crucial way on the fact that the weights here are chosen to be 1, that is the polynomials are ranked according to their natural degree.

 Denote by $\cH_n^d$ the space of polynomials of total degree $n$,  orthogonal to $\cP_{n-1}^d$
in $\cP_n^d$. Then 
$$
   \mbox{dim }\cP_n^d=\binom{n+d}{d}\,,
   \mbox{ and }\mbox{dim }\cH_n^d=\binom{n+d-1}{d-1}\,.
$$
As mentioned above, the choice  of a polynomial orthonormal basis  in $\cL^2(\mu)$   amounts to the choice of a basis for $\cH_n^d$, for any $n$. We are interested in the case where those polynomials are eigenvectors of the diffusion operator $\LL$ given by
$\LL(P)= -\lambda(P) P$, for any polynomial $P$ in the orthonormal basis, and for some real parameter $\lambda(P)$.

This leads us to state the following problem.

\bdefi[DOP problem]\label{gal.pb.compact}

Let $\Omega$ be a natural domain,  $\mu(dx)= \rho(x) dx$ a probability measure with smooth  positive density on $\Omega$, such that polynomials are dense in $\cL^2(\mu)$,  and let  $\LL $ be a diffusion operator (\ref{defdiffusion})
on $\Omega$.  We say that $(\Omega, \LL, \mu)$ is a solution to the Diffusion Orthogonal Polynomials problem (in short DOP problem) if  there exists an orthonormal  polynomial basis of $\cL^2(\mu)$
(see Definition~\ref{basis.order})
whose elements are at the same time eigenvectors of the operator $\LL$.

\edefi
Let us start with few elementary remarks.
Let $(\Omega,\LL,\mu)$ be a solution of the DOP problem.

The hypothesis on eigenbases in the subspaces $\cP_n^d$ implies
that $\LL$ maps $\cP_n^d$ into $\cP_n^d$ and   $\cH_n^d$ into $\cH_n^d$.
Therefore, when $P\in \cP_n^d$ and $Q\in \cP_m^d$, $\Ga(P,Q)\in \cP_{n+m}^d$. 

The restriction of $\LL$ to $\cP_n^d$ is symmetric for any $n$
(because it is so on an orthogonal eigenbasis), i.e.,
for any pair $(P,Q)$ of polynomials one has
\beq
      \label{IPP1}
      \int_\Omega P \LL(Q) d\mu = \int_\Omega Q\LL(P) d\mu.
\eeq

Using \eqref{IPP1} with $Q=1$ leads to $\int_\Omega \LL(P) \, d\mu=0$ for any polynomial. Applying this to $PQ$ together with the definition of the operator $\Ga$, one gets, for any pair $(P,Q)$ of polynomials
$$
    \int_\Omega \LL(PQ)\,d\mu = \int_\Omega P \LL(Q)\,d\mu
    + \int_\Omega Q\LL(P)\,d\mu  + 2 \int_\Omega \Ga(P,Q)\,d\mu =0,
$$
whence, using \eqref{IPP1}, we obtain
\beq
       \label{IPP2}
       \int_\Omega P\LL (Q)\,d\mu= \int_\Omega Q\LL(P)\,d\mu
        = -\int_\Omega \Ga(P,Q)\, d\mu.
\eeq
Applying when  $P=Q$ is an element of the basis,  since $\Ga(P,P)\geq 0$, one sees that $\lambda(P)\geq 0$.

From equation~\eqref{IPP2}, we see that the restriction of $\LL$ to polynomials is entirely determined by $\Ga$ (hence by  the matrices  $(g^{ij}(x))_{x\in \Omega}$), and the measure $\mu$.

Next, the following important observation, relying on the choice of the natural degree, shows that the DOP problem  is invariant under affine transformations:
\bprop 
If $(\Omega, \LL , \mu)$ is a solution to the DOP problem,  and if $A$ is an affine invertible transformation of $\bR^d$, so is 
$(\Omega_1, \LL_1, \mu_1)$, where $\Omega_1= A(\Omega)$, $\mu_1$ is the image measure through $A$ of $\mu$ and 
$$\LL_1(f)= \LL(f\circ A)\circ(A^{-1}).$$

\eprop

\bpf
Affine transformations map polynomials onto polynomials with the same degree. It suffices then to see that the associated operator $\LL_1(f)= \LL(f\circ A)\circ A^{-1}$ is again a diffusion operator, which has a family of orthogonal polynomials as eigenvectors:  if $P_k$ is an eigenvector for $\LL$, then $P_k\circ A^{-1}$ is an eigenvector of $\LL_1$. Moreover,   orthogonality  for the measure  $\mu$ is carried to orthogonality  for the measure $\mu_1$.

\epf

Moreover, the following Proposition shows that solutions to the DOP problem are stable under products

\bprop If $(\Omega_1, \LL_1, \mu_1)$ and $(\Omega_2, \LL_2, \mu_2)$ are solutions to the DOP problem in $\bR^{d_1}$ and $\bR^{d_2}$ respectively, then $(\Omega_1\times\Omega_2, \LL_1\otimes \Id+ \Id \otimes \LL_2, \mu_1\otimes \mu_2)$ is also a solution.

\eprop
\bpf Here $\LL= \LL_1\otimes \Id+ \Id \otimes \LL_2$ denotes the operator acting separately on $x$ and $y$: $\LL f(x,y) = \LL_x f + \LL_y f$. Similarly,  $\mu_1\otimes\mu_2$ is the product measure. The proof is then immediate: if $(P^{(1)}_k)$ and $(P^{(2)}_q)$ are the associated families of orthogonal polynomials, with eigenvalues $-\lambda_k$ and $-\mu_q$, the polynomials associated with $\LL$ are $P_{k,q}(x,y)= P^{(1)}_k(x)P^{(2)}_q(y)$, with associated eigenvalues $-\lambda_k- \mu_q$.

\epf

Next, we describe the general form of the coefficients of the operator $\LL$.

\bprop\label{prop.gal1.1}
Let $\LL$ be a diffusion operator in a natural domain $\Omega$ and $\mu$ be
a probability measure in $\Omega$ for which the polynomials are dense in $\cL^2(\mu)$.  Then $(\Omega,\LL,\mu)$ is a solution of the DOP problem if and only if 
\benum 
\item In the representation~\eqref{defdiffusion} of $\LL$,
for any $i= 1, \dots, d$, 
$b^i(x)\in \cP_1^d$ and for any $i,j= 1, \dots, d$, one has $g^{ij}(x)\in \cP_2^d$.

\item For any pair $(P,Q)$ of polynomials, equality~\eqref{IPP1} holds.
\eenum
\eprop

\bpf Assume that $(\Omega, \LL, \mu)$ is a solution of the DOP problem.  Since $\LL$ maps $\cP_n^d$ into $\cP_n^d$ for any $n\in \bN$, we have  $b^i(x)= \LL(x_i)\in \cP_1^d$ and $g^{ij}(x) = \Ga(x_i,x_j)\in \cP_2^d$. Moreover, writing any pair of polynomials $(P,Q)$ in the basis of orthogonal polynomials, we immediately obtain
equation~\eqref{IPP1}.

Conversely, if  $b^i(x)\in \cP_1^d$, $i= 1,\dots, d$  and $g^{ij}(x)\in \cP_2^d$,   $i,j= 1, \dots, d$, then $\LL$ maps $\cP_n^d$ into $\cP_n^d$ for any $n\in \bN$. Then, when moreover equation~\eqref{IPP1}
holds, $\LL$ is a symmetric operator on the finite dimensional space $\cP_n^d$, endowed with the scalar product inherited from the $\cL^2(\mu)$ metric. As such, we may find a basis of eigenvectors for it, and so we construct an $\cL^2(\mu)$ orthonormal  basis made of eigenvectors for $\LL$.

\epf

Only for polynomials the integration by parts formula~\eqref{IPP2} is a consequence of $\LL$ being a solution of the DOP problem. It may be interesting (and crucial) to extend it to any smooth compactly supported functions. This leads us to the Strong Diffusion Orthogonal Polynomials problem.

\bdefi[SDOP problem]\label{def.SDOP} The triple $(\Omega, \LL, \mu)$ is a solution to the
Strong Diffusion Orthogonal Polynomial problem  (SDOP in short) if it is a solution to
the DOP problem (Definition~\ref{gal.pb.compact})  and in addition, for any $f_1$ and $f_2$ smooth and compactly supported in $\bR^d$, one has
\beq
   \label{IPP3}
   \int_\Omega f_1\LL(f_2)\, d\mu= \int_\Omega f_2 \LL(f_1) \, d\mu.
\eeq
\edefi

If $(\Omega, \LL, \mu)$ is a solution to the SDOP problem, then,
writing $\mu(dx)= \rho(x) \,dx$,
we may define $\LL(f)$ by formula
\beq
\label
   {rep.L.sym}
   \LL(f) = \frac{1}{\rho} \sum_{ij} \partial_i \Big(g^{ij} \rho\, \partial_j f\Big)
\eeq
(see Proposition \ref{prop.L.rho} below)
and therefore $\LL$ is entirely determined from the (co)metric $g= (g^{ij})$ and the measure density $\rho(x)$. We therefore talk about the triple $(\Omega, g,\rho)$
as a solution of the SDOP problem.

Notice that $\LL$ admits a presentation in the form \eqref{rep.L.sym} under
assumptions weaker than those in Definition~\ref{def.SDOP}.
In Proposition~\ref{prop.L.rho} we do not demand that \eqref{IPP3} holds for
all compactly supported functions but only for those whose support is
contained in $\Omega$.

The equation~\eqref{rep.L.sym} allows to identify
$b^i$ from  $g^{ij}$ and $\rho$ as
\beq
      \label{eq.rho.sym}
      b^i= \sum_j \partial_j g^{ij} + \sum_j g^{ij} \partial_j \log \rho.
\eeq

To justify \eqref{rep.L.sym}, we start with the following two lemmas which
will be used again and again.

%%%%%%%%%%%%%%%%%%%%%%%%%%%%%%%%%%%%%%%%%%

\blem\label{lem.eq.IPP3}
Let $(\Omega,\mu)$ be any domain in $\bR^d$ and a measure on it.
Let $\cF$ be either the algebra of smooth functions compactly supported in $\bR^d$,
or its subalgebra consisting of functions compactly supported in $\Omega$.
Let $\LL$ be of the form \eqref{defdiffusion} with smooth coefficients.
Then the following conditions are equivalent:

\benum
\item The equation \eqref{IPP3} holds for any $f_1,f_2\in\cF$.

\item
The following equation \eqref{eq.IPP3}
holds for any $f_1,f_2\in\cF$:
\beq
    \label{eq.IPP3} 
    \int_\Omega\Big( f_1\LL(f_2)+\Ga(f_1,f_2)\Big)\,d\mu = 0.
\eeq
\eenum
\elem

\bpf
Equation \eqref{eq.IPP3} is derived from \eqref{IPP3} in the same way
as \eqref{IPP2} was derived from \eqref{IPP1}
(to justify the symmetry condition \eqref{eq.IPP3} in the case when $f_1=1$,
we replace $f_1$ by a function from $\cF$ which is equal to $1$ on the support of $f_2$).
The converse implication \eqref{eq.IPP3} $\Rightarrow$ \eqref{IPP3}
follows from the symmetry of $\Ga$.
\epf

\medskip

Let $dx_j^*=*(dx_j)$ be the differential $(d-1)$-form Hodge dual to $dx_j$, i.e.,
$$
     dx_j^* = (-1)^{j-1}\,
     dx_1\wedge\ldots\wedge\widehat{dx_j}\wedge\ldots\wedge dx_d.
$$
We set also
\beq
     \label{def.omega}
     \omega_{f} = \sum_{ij} \rho g^{ij} \partial_i f\, dx_j^*.
\eeq

\blem\label{lem.stokes}
Let $\Omega$ be a relatively compact natural domain
with piecewise smooth boundary,
$\mu=\rho\,dx$ with a smooth $\rho$, and
$\LL$ is given by \eqref{rep.L.sym}.
Let $f_1$ and $f_2$ be smooth functions
such that $f_1\rho$ extends
continuously to $\overline\Omega$.
Then the equation \eqref{eq.IPP3} is equivalent to
\beq
    \label{eq.IPP3.stokes}
    \int_{\partial\Omega} f_1\,\omega_{f_2} = 0.
\eeq
The latter equation can be equivalently rewritten as
$$
    \int_{\partial\Omega} f_1\sum_{ij}g^{ij}(\partial_i f_2)n_j\rho\,d\sigma=0
$$
where $(n_1,\dots,n_d)$ is the normal vector to the boundary
and $\sigma$ the surface measure.
\elem

\bpf
A straightforward computation shows that
\beq
    \label{d.omega}
    d(f_1\,\omega_{f_2}) = \big(f_1\LL(f_2)+\Ga(f_1,f_2)\big)\rho\,dx,
\eeq
thus the equivalence of \eqref{eq.IPP3} and \eqref{eq.IPP3.stokes}
follows from Stokes' formula.
\epf

\medskip

\bprop\label{prop.L.rho}
Let $\LL$ be defined by \eqref{defdiffusion} on a domain $\Omega\in\bR^d$,
and $\mu=\rho\,dx$ be a probability measure on $\Omega$ with a smooth
density $\rho$. Suppose that \eqref{IPP3} holds for any pair of smooth
functions compactly supported in $\Omega$. Then $\LL$ is of the form
\eqref{rep.L.sym}.
\eprop

\bpf
Let us temporarily denote the right hand side of \eqref{rep.L.sym} by $\hat\LL$,
and the corresponding square field operator by $\hat\Ga$.
Then $\hat\LL$ is of the form \eqref{defdiffusion} with the same $(g^{ij})$ but with
the $b^i$'s given by \eqref{eq.rho.sym}. So, we have $\hat\Ga=\Ga$.

Let $f_1$ and $f_2$ be functions compactly supported in $\Omega$.
Let $\Omega_0$ be a bounded domain with piecewise smooth boundary such that
$\overline{\text{supp}(f_1)}\subset\Omega_0$ and 
$\overline{\Omega_0}\subset\Omega$.
Then
\eqref{eq.IPP3.stokes} holds for $\hat L$ and $\Omega_0$,
hence \eqref{eq.IPP3} holds for $\hat L$ and $\Omega$.
On the other hand,
by Lemma \ref{lem.eq.IPP3}, we have \eqref{eq.IPP3} for $L$ as well.
Since $\hat\Ga=\Ga$, we deduce that
$$
   \int_\Omega f_1\LL(f_2)\,d\mu = -\int_\Omega \Ga(f_1,f_2)\,d\mu
   =\int_\Omega f_1\hat\LL(f_2)\,d\mu
$$
for any $f_1,f_2$ compactly supported in $\Omega$ whence $\LL=\hat\LL$.
\epf

\medskip
 The next proposition shows that
 the distinction between DOP and SDOP solution is relevant in the non compact case only.

\bprop \label{lem.IPP3} Whenever $\Omega$ is relatively compact, any solution of the DOP problem is a solution of the SDOP problem.
\eprop

\bpf (See also \cite[p.~155, Cor.~2]{Treves1967}.)
We just have to show that for relatively compact sets $\Omega$, 
equation~\eqref{IPP3} is satisfied for any pair $(f_1,f_2)$
of smooth compactly supported functions.
Since $\Omega$ is relatively compact,  for any $f$ smooth  and  compactly supported in $\bR^d$ (and not necessarily in $\Omega$),  we first choose some compact $K$ which contains both the support of $f$ and $\Omega$, and which is a hyper-rectangle oriented parallel to the coordinate axes. Then, there exists  a polynomial  sequence  $(R_n)$   converging uniformly on $K$ to  $\partial_{11\cdots dd} f$. Then, repeated  integrals of $R_n$ converge uniformly on $K$ to the corresponding derivatives of $f$. Finally, we obtain a sequence
$(P_n)$ of polynomials such that $P_n$ and all it's partial derivatives of
order  1 and 2 converge uniformly on $K$ to the corresponding derivatives of $f$.
Choose such sequences $(P_n)$ and $(Q_n)$ for $f_1$ and $f_2$ respectively.   The functions $g^{ij}$ and $b^i$ being polynomials, are bounded on $K$. Therefore,    $(P_n)$, $(Q_n)$, $\LL(P_n)$, and $\LL(Q_n)$
converge uniformly on $K$  to $f_1$ ,$f_2$, $\LL(f_1)$, and $\LL(f_2)$
respectively. Then, it is  clear that formula~\eqref{IPP1} 
extends immediately to the pair  $(f_1,f_2)$.
\epf

\medskip
\bprop\label{prop.gal1}
  If $(\Omega, g, \rho)$ is a solution to the SDOP  problem, then
  there exist polynomials $L^i\in \cP_1^d $, $i= 1, \dots, d$
  (that is polynomials  of degree at most  $1$)  such that, for any
  $x\in \Omega$  and any $i= 1,\dots, d$,
\beq
\label{eq.rho.1}
               \sum_j g^{ij} \partial_j \log (\rho(x)) = L^i(x).
\eeq
\eprop

\bpf
Combine \eqref{eq.rho.sym} with the fact that $\deg g^{ij}=2$ and $\deg b^i=1$.
\epf

\medskip

As a consequence of Proposition~\ref{prop.gal1}, we get the following
general description of the admissible measures (Proposition~\ref{prop.measure}).
We start with the following fact. % (immediate from a basic calculus course).

\blem\label{lem.calculus}
  Let $C$ and $\Delta=\Delta_1^{m_1}\dots\Delta_s^{m_s}$ be
  polynomials in one complex variable where $m_1,\dots,m_s$ are positive integers.
  Assume that all roots of the product $\Delta_1\dots\Delta_s$ are simple.
  Let $F$ be a holomorphic function on the complement of the roots of $\Delta$
  such that $F'=C/\Delta$. Then $F\prod_q\Delta_q^{m_q-1}$ extends to a
  polynomial of degree at most $\max\big(0,1+\deg C-\sum_q\deg\Delta_q\big)$.
\elem

\bpf
The function $F$ is univalued and $F'$ is rational, hence $F$ is rational.
The multiplicity of poles of $F$ at the zeros of $\Delta_k$ is at most $m_k-1$,
hence $F\prod_q\Delta_q^{m_q-1}$ extends to a polynomial.
For a rational function $f=p/q$ where $p$ and $q$ are polynomials, we
set $\deg_\infty f = \deg p-\deg q$. It remains to observe that
$\deg_\infty f'\le\max(0,\deg f-1)$.
\epf

\medskip
Notice that if a real polynomial is irreducible over $\bR$ but reducible over $\bC$,
then it factors over $\bC$ as $(\cR+i\cI)(\cR-i\cI)$ with $\cR$ and $\cI$ irreducible
over $\bC$, and thus it is equal to $\cR^2+\cI^2$.

\bprop[General form of the measure]
\label{prop.measure}~
Let $(\Omega, g, \rho)$ be a solution of the SDOP problem.
Suppose that the determinant $\Delta$ of $(g^{ij})$ writes 
$\Delta= \Delta_1^{m_1} \cdots \Delta_p^{m_p}$, where $\Delta_k$
are irreducible over the reals. Let $J$ the set of indices $j\in \{1,\dots, p\}$
such that $\Delta_j$ is reducible 
over $\bC$, written  $\Delta_j = \cR_j^2+ \cI_j^2$.
Then there exist real constants
$(\alpha_k)_{k=1,\dots,p}$ and $(\beta_j)_{j\in J}$,
and a polynomial $Q$ such that
\beq
    \label{mes.deg}
    \deg(Q) \leq 2d-\sum_{k=1}^p\deg\Delta_k,\quad
    \deg_{x_i}(Q) \leq 2d-\sum_{k=1}^p\deg_{x_i}\Delta_k,\quad
    i=1,\dots,d,
\eeq
and
\footnote{By $\arctan(\cI_j/\cR_j)$ we mean here a continuous single-valued 
branch of the argument of $\cR_j+i\cI_j$. So, formally speaking, one should
replace $\arctan(\cI_j/\cR_j)$ by $\arctan(\cI_j/\cR_j)+c(x)$ where $c(x)$
is a locally constant function
on $\Omega\setminus\{\cR_j=0\}$ which jumps by $\pm\pi$ when crossing
$\{\cR_j=0\}$.}
\beq
\label{measure.form1}
   \rho =
   \exp\Big( \frac{ Q}{\Delta_1^{m_1-1}\cdots \Delta_p^{m_p-1}}
   +\sum_{j\in J}\beta_j \arctan\frac{\cI_j}{\cR_j} \Big)
   \prod_{k=1}^p |\Delta_k|^{\alpha_k}
\eeq
\eprop

\bpf

With $h= \log \rho$, one has  from equation \eqref{eq.rho.1}

\begin{equation}
                \label{mesure}
                \partial_j h = \sum_i g_{ij}  L^i,
\end{equation}
where $g^{(-1)}=(g_{ij})$ is the inverse matrix of $(g^{ij})$ and $L^i\in \cP_1^d$.
But $g^{(-1)}= \Delta^{-1} \hat g$, where $\hat g$ is the  matrix of co-factors of $g$.
Then each $\hat g_{ij}$ is a polynomial of degree at most $2d-2$,
and therefore $\partial_i h = C_i/\Delta$ where $C_i\in \cP_{2d-1}^d$.

Let us extend the differential form $dh$ to  a closed holomorphic form $\omega$
in the complex domain $\bC^d\setminus\{\Delta=0\}$. 
By Alexander duality (see \cite{Alex,pontry,Lef}), the De Rham cohomology group
$H^1_{DR}(\bC^d\setminus\{\Delta=0\})$ is generated by the 1-forms
$d\log\hat\Delta_q$ where
$\hat \Delta_1,\dots,\hat\Delta_s$, $s=p+|J|$, are the irreducible over $\bC$ factors
of $\Delta$. Hence there exist complex numbers $\gamma_1,\dots,\gamma_s$ such that 
$\omega- \sum_q \gamma_q d\log(\hat \Delta_q)= dF_0$ is exact.
From the definition of $\omega$ we know that
$$
  \partial_i F_0= \frac{C_i}{\Delta}-\sum_q \gamma_q \frac{\partial_i \hat \Delta_q}
  {\hat \Delta_q}= \frac{\hat C_i}{\Delta}
$$
with ${\deg}~\hat C_i \leq 2d-1$.

%By partial fraction decomposition, we know
By Lemma~\ref{lem.calculus},
when fixing generically all variables $x_j$ for
$j\neq i$, then $Q= F_0\prod_q \hat\Delta_q^{m_q-1}$ is a polynomial in $x_i$
of degree at most
$$
  n_i = 2d-\sum_{q=1}^s\deg_{x_i}\hat\Delta_q = 2d-\sum_{k=1}^p\deg_{x_i}\Delta_k.
$$
Therefore $\partial_1^{n_1}\cdots \partial_d^{n_d}Q=0$.
Hence $Q$ is a polynomial, and its $x_i$-degrees are as in \eqref{mes.deg}.
Moreover, since the same remains
true for any coordinate system, the total degree of $Q$ is also as in \eqref{mes.deg}.
\footnote{The anonymous referee pointed out that similar relations between 
exact meromorphic $p$-forms and their primitives are found in \cite{Griffiths1969, Dimca1991} for any $p$.}
Thus we obtain \eqref{mes.deg} and
\beq\label{log.rho}
   \log\rho = h = \frac Q{\hat\Delta_1^{m_1-1}\cdots\hat\Delta_s^{m_s-1}}
     +\sum_q\gamma_q\log\hat\Delta_q\,.
\eeq
We now deal with the real form of $\rho$. Whenever there is an irreducible over
$\bR$ factor $\Delta_k$ of $\Delta$ which is reducible over $\bC$, its irreducible decomposition over $\bC$ writes $\Delta_k = (\cR_k+ i \cI_k)(\cR_k-i\cI_k)$,
and the corresponding summand in $\log\rho$ must be of the form 
\[
  \gamma_q\log(\cR_k+i\cI_k)+\bar\gamma_q\log(\cR_k-i\cI_k)\,,
\] 
which writes in real form
$\alpha_k\log\Delta_k + \beta_k\arctan({\cI_k}/{\cR_k})$.
\epf

\medskip 
\brmq In the case where $\Omega$ is bounded, we did not observe up to now
any model where the admissible measures
has the exponential factor in~\eqref{measure.form1}.
Moreover, only components of the reduced boundary equation (see
Definition~\ref{def.red.eq}) appear in all known examples.
In the unbounded case the exponential term must be present
(otherwise the measure would not integrate all the polynomial functions),
however, the fraction in~\eqref{measure.form1} reduces to a polynomial
after cancellation in all known examples.
\ermq

As we see in the proof of Proposition 2.14, the real
1-form $d\log\rho$ extends to a meromorphic 1-form in $\bC^2$
which may have poles only on the algebraic curve $\det g^{ij}=0$.
By abusing notation, we shall still denote this meromorphic
form by $d\log\rho$.

\bprop\label{rho=infty}
In the setting of Proposition \ref{prop.measure}, assume that
$d\log\rho$ has pole along the hypersurface $\Delta_k=0$ (this means that
either $\alpha_k\ne 0$ for $k\not\in J$, or $(\alpha_j,\beta_j)\ne(0,0)$
for $k\in J$, or $m_k\ge2$ and $\Delta_k^{m_k-1}$ does not divide $Q$).
Then there exist polynomials $S_k^i\in\cP_1^d$, $i=1,\dots,d$,
such that
\beq
    \label{eq.rho=infty}
    \sum_j g^{ij}\partial_j\Delta_k = S_k^i\Delta_k
    \qquad\text{for any $i=1,\dots,d$.}
\eeq
\eprop

\bpf
From the point of view of the geometric intuition, this fact is almost obvious.
Indeed, the condition \eqref{eq.rho.1} means that for any $i$, the derivative of
$\log\rho$ along the vector field $\sum_j g^{ij}\partial_j$ is bounded.
Therefore it is clear that this vector field should be tangent
to the hypersurface $\log\rho=\infty$.

Let us give however a more formal proof. We shall use the notation introduced in
the proof of Proposition \ref{prop.measure}.
Let us differentiate \eqref{log.rho} with respect to $x_j$.
Our assumption about $k$ implies that we obtain
\beq\label{d.log.rho}
   \partial_j\log\rho = \frac{P\partial_j\hat\Delta_k+R_j\hat\Delta_k}
                             {\hat\Delta_k^n\tilde\Delta}
\eeq
where $n>0$, $P$ is a polynomial coprime with $\hat\Delta_k$ (which does not
depend on $j$), $\tilde\Delta$ is a product of some powers of the $\hat\Delta_q$'s
with $q\ne k$.
% and the degree of the numerator is less
% than the degree of the denominator.
After plugging 
\eqref{d.log.rho} into \eqref{eq.rho.1} and multiplying by the denominator,
we obtain
$$
   P\Big(\sum_j g^{ij}\partial_j\hat\Delta_k\Big)
    +\hat\Delta_k\Big(\sum_j g^{ij}R_j\Big)=
                             L^i\hat\Delta_k^n\tilde\Delta
$$
Since $P$ is coprime with $\hat\Delta_k$, we deduce that
$\hat\Delta_k$ divides $\sum_j g^{ij}\partial_j\hat\Delta_k$
and we denote the quotient by $S_k^i$. Since the degree of the left hand side
of \eqref{eq.rho=infty} is at most $1+\deg\hat\Delta_k$, we conclude that
$S_k^i\in\cP_1^d$.
\epf

\medskip
\bcor\label{cor.mes.deg}
In the setting of Proposition~\ref{prop.measure}, the
estimate for $\deg_{x_i}Q$ in~\eqref{mes.deg} can be improved by replacing
$2d$ with $2+\max_j\deg_{x_i}\hat g_{ij}$ where $(\hat g_{ij})$ is the co-matrix
of $g$, i.e., $\hat g_{ij}$ is the complementary minor of the entry $g^{ij}$.
\ecor

\bpf
Let $F_0=Q\prod_k\Delta_k^{1-m_k}$ and let $L_0^i=\sum_j g^{ij}\partial_j F_0$.
By combining equations \eqref{eq.rho.1},
\eqref{log.rho}, and \eqref{eq.rho=infty}, we obtain
$$
  L_0^i
%  \sum_j g^{ij}\partial_j F_0
  =\sum_j g^{ij}\Big(\partial_j\log\rho-\sum_k\gamma_k\partial_j\log\hat\Delta_k\Big)
  = L^i - \sum_k \gamma_k S_k^i\in\cP_1^d.
$$
The rest of the proof is similar to the proof of \eqref{mes.deg}.
Namely, the definition of $L_0^i$ implies that
$\Delta\partial_j F_0 = \sum_i \hat g_{ij} L_0^i$, and the required estimate for
$\deg Q$ follows from Lemma~\ref{lem.calculus}. When applying Lemma~\ref{lem.calculus},
we may get rid of $\max(0,\dots)$ because
$$
   \sum_k \deg_{x_j}\Delta_k \le \deg_{x_j}\Delta
    = \deg_{x_j}\sum_i g^{ij}\hat g_{ij} \le 2+\max_i\deg_{x_j}\hat g_{ij}
$$

%
% $\partial_j F_0$
% is a rational function of the required form.
\epf
\medskip

\bcor\label{cor.measure} Let $(\Omega,g,\rho)$ be a solution to the SDOP problem.
Suppose that $\Omega$ contains a half-cylinder, i.e.,
a domain $\Omega_1\subset\Omega$ given in some affine coordinates by
$x_1>0$ and $x_2^2+\dots+x_d^2<1$.
Then
$$
   \deg_{x_1}\det(g) < 2+\max_j\deg_{x_1}\hat g_{1j}.
$$
\ecor

\medskip
\bpf
Let notation be as in Proposition~\ref{prop.measure}.
Then $\rho$ is given by \eqref{measure.form1}.
Let $F=\Delta_1^{m_1-1}\dots\Delta_p^{m_p-1}$ be the denominator of the
fraction in \eqref{measure.form1}. We may assume that the sign of each $\Delta_k$
is chosen so that $\Delta_k>0$ on $\Omega$.
Write $Q=\sum_{j=0}^nQ_j x_1^j$ and $F=\sum_{j=0}^mF_j x_1^j$ with
$Q_j,F_j\in\bR[x_2,\dots,x_n]$ and $n=\deg_{x_1}Q$, $\,m=\deg_{x_1}\Delta$.
Let $\bB^{d-1}_r=\{(x_2,\dots,x_d)\mid x_2^2+\dots+x_d^2<r^2\}$ with $0<r<1$.
We have $F>0$ in $\Omega$, hence $F_m>0$ in the unit $(d-1)$-ball.
Therefore $F_m>C_1$ in $\bB^{d-1}_r$ for some constant $C_1>0$.
Let $C_2$ be a constant such that $|Q_j(x)|<C_2$, $j=0,\dots,m-1$
and $|F_j(x)|<C_2$, $j=0,\dots,n$, when $x\in\bB_r^{d-1}$.
Then, for some constants $A$ and $C$ depending on $C_1$, $C_2$, we have
$|Q/F|<C x_1^{n-m}$ on $\Omega_2=[A,\infty)\times\bB^{d-1}_r$. Thus, $n>m$ because
otherwise we would have
$\exp(Q/F) > -\exp(C)$ on $\Omega_2$ which contradicts the integrability of polynomials
on $\Omega$.
\epf

\medskip
\bdefi\label{def.red.eq}
Given a natural domain $\Omega$ in $\bR^d$ not coinciding with the whole $\bR^d$,
let $I(\partial\Omega)$ be the ideal of $\bC[x_1,\dots,x_d]$
consisting of polynomials identically vanishing on $\partial\Omega$.
If $I(\partial\Omega)\ne\{0\}$, then the condition that $\Omega$ is the interior
of its closure implies that $I(\partial\Omega)$ is a principal ideal
generated by a single real polynomial $\hat\Delta$ which is, evidently,
reduced (i.e., does not have multiple factors).
In this case we say that 
$\hat\Delta$ is the {\em reduced equation} of $\partial\Omega$.
Each irreducible factor of $\hat\Delta$ vanishes on some
open subset of the set of smooth points of $\partial\Omega$.
\edefi

\medskip
We can now state the main result of this section.

%%%%%%%%%%%%%%%%%%%%%%%% MAIN THEOREM of the section %%%%%%%%%%%%%%%%%%%%%%%%%%

 \bthm \label{thm.gal.bord.Omega}~
 Let $\Omega$ be a natural domain in $\bR^d$, $\rho$ a smooth function
 in $\Omega$ and $g=(g^{ij})$ a positive definite (co)metric in $\Omega$.
 Let $\Delta=\det(g)$. Then $(\Omega,g,\rho)$ is a solution to the
 SDOP problem (recall that it is the same as DOP problem when $\Omega$ is bounded)
 if and only if there exists a reduced (i.e., without multiple factors) real
 polynomial $\hat\Delta$ such that
 $\hat\Delta$ divides $\Delta$ and the following conditions hold:
\benum
\item
     \label{thm.gal.gij}
      For any $(i,j)$, $g^{ij}(x) \in \cP_2^d$;
\item
     \label{thm.gal.de}
     $\partial\Omega$ is contained
     in the algebraic hypersurface $\{\hat\Delta=0\}$.
\item
     \label{thm.gal.bord}
     For any $i=1,\dots,d$, for some $S^i\in\cP_1^d$ one has
     \beq
         \label{eq.thm.gal}
         \sum_j g^{ij} \partial_j\hat\Delta = \hat\Delta S^i
     \eeq
\item
     \label{thm.gal.rho}
     $\rho$ is of the form \eqref{measure.form1} (with the ingredients explained
     in Proposition \ref{prop.measure}),
 and polynomials are dense in $\cL^2(\rho\,dx)$
 (if $\Omega$ is bounded, the last condition is equivalent
 to $\int_\Omega\rho\,dx<\infty$).

\item
     \label{thm.gal.bi}
     $\sum_j g^{ij} \partial_j \log \rho\in \cP_1^d$ 
     for any $i= 1, \dots, d$.
\eenum
\ethm
\def\LastConditionInThmGal{5}

\brmq \label{rmq.thm.gal}
Condition \eqref{thm.gal.bord} can be equivalently reformulated
as follows. Let $\Delta_1\dots\Delta_r$ be a factorization
(not important, over $\bR$ or over $\bC$)
of $\hat\Delta$. Then, for any $k=1,\dots,r$ and for any
$i=1,\dots,d$, there are $S_k^i\in\cP_1^d$ such that
\beq
     \label{eq.thm.gal.1}
     \sum_j g^{ij} \partial_j\Delta_k = \Delta_k S^i_k
\eeq

This is also equivalent to the fact that for any $i$,
the differential $(d-1)$-form
$\sum_j g^{ij}\,dx_j^*$ restricted to (the smooth part of) $\partial\Omega$
identically vanishes.
\ermq

 \brmq Equation~\eqref{eq.thm.gal.1} may be rewritten in a more intrinsic way as 
 $$\Ga(\log \Delta_k, x_i)= S^i_k,$$ and similarly for equation~\eqref{eq.thm.gal}.
 \ermq

\medskip
\bpf

{\bf Necessity.}
Suppose that $(\Omega,g,\rho)$ is a solution to the SDOP problem and
let us prove conditions~(1) -- (\LastConditionInThmGal).
The last two of them and the first one are just a rephrasing of
Propositions~\ref{prop.gal1.1}, \ref{prop.gal1}. and~\ref{prop.measure}.

Let us prove that $\Delta$ vanishes on $\partial\Omega$. Let
$x_0$ be a smooth point of $\partial\Omega$.
If $\log\rho(x_0)=\pm\infty$, then $\Delta(x_0)=0$ 
by Proposition~\ref{rho=infty}. Indeed, in this case $\Delta_k(x_0)=0$ for some
$\Delta_k$ satisfying the hypothesis of that lemma. Hence
$(\partial_j\Delta_k(x_0))_{j=1}^d$ is a non-zero solution of a system
of linear equations with the coefficient matrix $(g^{ij}(x_0))$ whence
$\Delta(x_0)=\det g(x_0)=0$.

Suppose now that $0<\rho(x_0)<\infty$.
Assume first that $\partial\Omega$ is piecewise smooth.
Let us choose a neighborhood $B(x_0,r)$
of $x_0$ on which $\rho<\infty$. Let
$\omega^i = \omega_{x_i}$ (in the notation of \eqref{def.omega}).
Then, for any  function $f$ smooth and compactly supported in $B(x_0,r)$,
for any $i$, one has by Lemmas~\ref{lem.eq.IPP3} and~\ref{lem.stokes} that
$\int_{\partial\Omega} f\,\omega^i = 0$. The last equality can by
rewritten in the form
$$
   \int_{\partial\Omega} f \sum_j g^{ij} n_j \rho\,d\sigma = 0
$$
(in the notation of Lemma~\ref{lem.stokes}). This equality holds for any
$f$ supported in $B(x_0,r)$. Hence, for any $i$ we have
\beq
    \label{eq.thm.gal.2}
    \sum_j g^{ij}(x_0)n_j(x_0)=0
\eeq
and once again we obtain a non-zero solution to a system
of linear equations with coefficients $g^{ij}(x_0)$ whence $\Delta(x_0)=0$.
So, we proved that $\partial\Omega\subset\{\Delta=0\}$.
For the case when $\partial\Omega$ is not a priori assumed to be piecewise smooth,
% the vanishing of $\Delta$ on $\partial\Omega$ is proven in
see Lemma~\ref{lem.bad.bndry} below. In its proof we use more or less
the same arguments (basically, integration by parts) but some additional
tricks are needed since the Stokes formula cannot be applied in this case.

Let $\hat\Delta$ be the reduced equation of $\partial\Omega$, i.e., the generator
of the ideal $I(\partial\Omega)$
(see Definition~\ref{def.red.eq}). We have $\Delta\in I(\partial\Omega)$,
hence $\hat\Delta$ divides $\Delta$. So, we proved condition~\eqref{thm.gal.de}.

To prove \eqref{thm.gal.bord}, notice that \eqref{eq.thm.gal.2}, which holds when
$0<\rho(x_0)<\infty$, combined with Proposition~\ref{rho=infty} imply that, for any $i$,
the left hand side of~\eqref{eq.thm.gal} identically vanishes on $\partial\Omega$,
i.e., belongs to $I(\partial\Omega)$. Hence it is equal to $\hat\Delta S^i$
for some polynomial $S^i$. By comparing the degrees, we conclude that $S^i\in\cP_1^d$.

\medskip
{\bf Sufficiency.}

{\it Step 1.}
Suppose that conditions~(1) -- (\LastConditionInThmGal) hold.
Let us write $\rho$ as in~\eqref{measure.form1} but with the
fraction $Q/(\dots)$ replaced by its reduced form $R/\prod\Delta_k^{n_k}$
where $\Delta_k$ does not divide $R$ unless $n_k=0$.
Since none of $\Delta_k$ vanishes in $\Omega$, we may assume that
$0<\Delta_k<1$ on $\Omega$ for each $k$.

Assume first that $\rho$ extends up to a continuous function on the closure
of $\Omega$ which vanishes at any smooth point of $\Omega$.
This means that for each $k$ such that $\Delta_k$ is a factor of $\hat\Delta$,
we have  $\alpha_k>0$ when $n_k=0$, and%
\footnote{%
However at some ``corners" of $\Omega$ where some $\Delta_q$
not included in $\hat\Delta$ vanishes, {\it a priori\/} $\rho$  
may be discontinuous if $n_q>0$ and $R|_\Omega>0$ somewhere near this ``corner".}
\beq\label{pf.thm.gal}
  \text{$R|_\Omega<0$ near $\{\Delta_k=0\}$
  when $n_k>0$.}
\eeq
In this case the result immediately follows from Proposition~\ref{prop.gal1.1}
combined with Lemmas~\ref{lem.eq.IPP3} and~\ref{lem.stokes}.

{\it Step 2.}
Now we turn to the general case.
According to Proposition~\ref{prop.gal1.1}, it is enough to prove
that equation~\ref{IPP1} holds for any two polynomials $P_1$ and $P_2$.
So, we fix $P_1$ and $P_2$ and we are going to vary the coefficients $\alpha_k$
in~\eqref{measure.form1}. Namely, up to renumbering the factors of $\Delta$,
we may assume that $\hat\Delta=\Delta_1\dots\Delta_r$. So,
for any $a=(a_1,\dots,a_r)\in\bC^r$,
we define $\rho_a$ by formula~\eqref{measure.form1} where $\alpha_1,\dots,\alpha_r$
are replaced with $a_1,\dots,a_r$. We set also $\alpha=(\alpha_1,\dots,\alpha_r)$.
Define $\LL_a$ by~\eqref{rep.L.sym} with $\rho_a$ standing for $\rho$.
Condition~\eqref{thm.gal.bord} ensures that each $\LL_a$ has the
form~\eqref{defdiffusion} with some $b^i_a$ standing for the $b_i$ and,
moreover, $a\mapsto b^i_a$, $i=1,\dots,d$, are affine linear functions on $\bC^r$. Indeed,
by equations~\eqref{eq.rho.sym} and~\eqref{eq.thm.gal.1},
for any $i$ we have:
$$
      b^i_a-b^i_\alpha = \sum_{i} g^{ij}\partial_j
     \log\frac{\rho_a}{\rho}
       = \sum_{i,k}(a_k-\alpha_k) g^{ij}
        \frac{\partial_j\Delta_k}{\Delta_k}
       = \sum_{k}(a_k-\alpha_k)S_k^i\,.
$$
A similar computation shows that $g$ and $\rho_a$ satisfy
condition~\eqref{thm.gal.bi}.

For fixed polynomials $P_1$ and $P_2$, we set
$$
    F(a) = \int_\Omega \Big(P_1\,\LL_a(P_2) - P_1\,\LL_a(P_2)\Big)\rho_a\,dx\,.
$$
Let $U_0=\{a\in\bC^r\mid \Re a_k \ge \alpha_k$ for all $k=1,\dots,r\,\}$.
Since $\rho_a=\rho\prod_{k\le r}\Delta_k^{a_k-\alpha_k}\le\rho$ on $U_0$
(recall that $0<\Delta_k<1$ on $\Omega$),
the function $F$ is defined (and is finite) 
in some domain $U$ containing $U_0$.
Moreover, it has the form
$$
    F(a) = \int_\Omega G(x,a)\,dx
$$
where $G(x,a)$ is an integrable function on $\Omega\times U$ which is complex
analytic with respect to $a$ for any fixed $x$. Hence $F$ is an analytic
function of $a$. Indeed, if we fix all variables except some $a_k$ and
let $a_k$ vary in the half-plane $\Re z>\alpha_k-\varepsilon$, $0<\varepsilon\ll1$,
then the integral of this function along each closed path is zero
by Fubini theorem.

The above arguments show that $(\Omega,g,\rho_a)$ satisfies all the
conditions~(1) -- (\LastConditionInThmGal)
Furthermore, the condition \eqref{pf.thm.gal} is also satisfied because
otherwise $\int_\Omega\rho\,dx$ would not be finite.
Therefore, by the result of Step 1, we have $F(a)=0$ when $a_k>0$ for all
$k=1,\dots,r$.
Hence $F=0$ on the whole $U$ which completes the proof.
\epf

\medskip
\bcor\label{cor.gal}
Let $\Omega$ be a natural bounded domain and $g$ a smooth (co)-metric in it.
A solution of the DOP problem with given $\Omega$ and $g$ exists if
and only if Conditions (1)--(3) of Theorem~\ref{thm.gal.bord.Omega} hold
for some reduced factor $\hat\Delta$ of $\det(g)$.

In this case one can choose any measure $\mu=\rho\,dx$ of the form
$\Delta_1^{a_1}\dots\Delta_p^{a_p}$, where $\Delta_1\dots\Delta_p=\hat\Delta$,
under condition that $\mu(\Omega)<\infty$, for example, one can choose for the
$a_k$ any non-negative real numbers.
\ecor
%%%%%%%%%%%%%%%%%%%%%%%% END OF MAIN THEOREM %%%%%%%%%%%%%%%%%%%%%%%%%%

\medskip
\bcor\label{cor.gal.2}
Let $(\Omega,g,\rho)$ be a solution of the SDOP problem in $\bR^d$,
and let $\Delta=\det(g)$.
If $\deg\Delta=2d$ and $\Delta$ does not have multiple factors, then
$\Omega$ is bounded.
\ecor

\bpf
By Theorem~\ref{thm.gal.bord.Omega}, the boundary of $\Omega$ is contained in
an algebraic hypersurface. On the other hand, the condition that
$\Delta$ is square-free and $\deg\Delta=2d$ combined with
Proposition~\ref{prop.measure} imply that
$\rho=\Delta_1^{\alpha_1}\dots\Delta_p^{\alpha_p}$
with polynomials $\Delta_1,\dots,\Delta_p$. Thus the unboundedness of $\Omega$
contradicts the integrability condition for polynomials. 
\epf

%%%%%%%%%%%%%%%%%%%%%%%% END OF MAIN THEOREM %%%%%%%%%%%%%%%%%%%%%%%%%%

\bigskip
The following lemma is needed in the proof of Theorem~\ref{thm.gal.bord.Omega}
only in the case when it is not a priori assumed that the boundary of $\Omega$
is piecewise smooth.

\medskip

\blem\label{lem.bad.bndry}
Let $(\Omega,g,\rho)$ be a solution of the SDOP problem.
Let $p_0$ be a point on $\partial\Omega$ such that $\rho(p_0)\ne 0$.
Then $\Delta(p_0)=0$.
\elem

\bpf
Suppose that $\Delta(p_0)\ne 0$, i.e., $g(p_0)$ is non-degenerate.
Let us choose coordinates so that
% $p_0$ is the origin and
$g^{ij}(p_0)=\delta^{ij}$.
For a unit vector $v$, we consider the linear function $l_v:x\mapsto v\cdot x$
and the derivation $f\mapsto \Gamma(l_v,f)$ which we denote by $\partial_v$.
% Let $h_v$ be the solution of the first order PDE $\partial_v u=1$
% with the initial coditions $u=0$ on the hyperplane $v^\perp$.
% Then $d h_v(p_0)=l_v$.
% and $h_v$ is $C^1$-close to $l_v$ in small neighbourhoods of $p_0$.

Using standard properties of submersions, it is easy to show that
there is a sufficiently small ball $U$ centered at the origin such that
for any two points $p,q\in U$, there exists a unit vector $v$ such that
$q$ lies on the trajectory of the vector field $\partial_v$ starting at $p$.
Let us fix a ball $U$ with this properties and choose points $p\in U\cap\Omega$
and $q\in U\cap\text{Int}(\bR^d\setminus\Omega)$
(this is possible because $\Omega$ coincides with the interior of its closure).
Let us choose $v$ as explained above (so that $q$ is on the
trajectory of $\partial_v$ starting at $p$).
Let us choose curvilinear coordinates
$(x_1,\dots,x_d)$ in $U$ so that $\partial_v=\partial_1:=\partial/\partial x_1$,
$p=(a,0,\dots,0)$, and $q=(b,0,\dots,0)$ with $a<b$.
We may assume that $U$ is small enough, so that $\Delta|_U\ne0$, hence by
Proposition~\ref{prop.measure} we may extend $\rho$ to a non-zero analytic function
in $U$.

\begin{figure}[ht]

\centering \includegraphics[width=.8\linewidth]{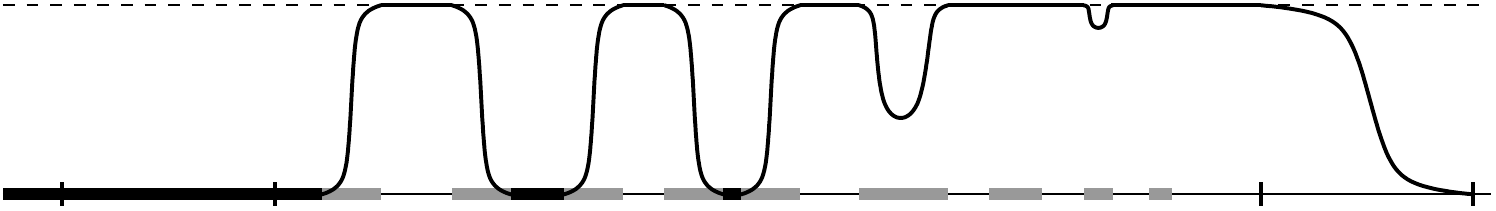}
\vskip-3pt
\centerline{\hskip15pt $a_1$ \hskip 25pt $a+r$\hskip 200pt $b_1$ \hskip25pt $b+r$}
\caption{The graph of $f_\varepsilon(t,0)$. Bold black: $\Omega_{2\varepsilon}$;
grey: $\Omega\setminus\Omega_{2\varepsilon}$.}
		\label{fig:BadDomain}
\end{figure}

Let $r>0$ be such that $B_r^1(a)\times B_r^{d-1}\subset\Omega$ and
$B_r^1(b)\times B_r^{d-1}\subset\bR^d\setminus\Omega$.
Let $\varphi$ (a smoothing kernel)
be a smooth non-negative function supported in the unit ball
such that $\int\varphi\,dx=1$. For $\varepsilon>0$ we set
$\varphi_\varepsilon(x) = \varphi(x/\varepsilon)/\varepsilon^d$, i.e.
$\text{supp}(\varphi_\varepsilon)=B_\varepsilon^d$ and
$\int\varphi_\varepsilon\,dx=1$.
Let $\Omega_\varepsilon = \{x\in\Omega\mid\text{dist}(x,\partial\Omega)
>\varepsilon\}$ and let
$h_\varepsilon = (1-\indic_{\Omega_\varepsilon})*\varphi_\varepsilon$ where
$\indic_{\Omega_\varepsilon}$ is the characteristic function of
$\Omega_\varepsilon$ and $*$ denotes the convolution.
Let $F_0$ be the finitely supported function such that
$\partial_1 F_0(x)=\varphi_r(x-p) - \varphi_r(x-q)$
%(the support of $F_0$ is the convex hull of $B_r(p)\cup B_r(q)$)
and, finally, we set $f_\varepsilon=F_0h_\varepsilon$ (see Figure~\ref{fig:BadDomain}).
% (see Figure ... where we show the graphs of the
% functions $F_0$, $\partial_1F_0$, $f_\varepsilon$ and
% $\partial_1 f_\varepsilon$ restricted on the line $(pq)$; the thick
% segments represent $(pq)\cap\Omega$). 

Observe that $\Omega\cap\text{supp}(\partial_1f_\varepsilon)
    = V\cap\text{supp}(\partial_1f_\varepsilon)$
% $f_\varepsilon$ vanishes outside $V\cap\Omega$
where $V=[a_1,b_1]\times B_r^{d-1}$ with $a_1=a-r$ and $b_1=b-r$.
Notice also that for any $y=(x_2,\dots,x_d)$
we have
$$
    f_\varepsilon(a_1,y)=0, \qquad
    f_\varepsilon(b_1,y)=F_0(b_1,y)=\int_{\bR} \varphi_r(t,y)\,dt,
$$
and by the choice of coordinates we have $\Gamma(l_v,f)=\partial_1f$. Hence
$$
\begin{aligned}
    \int_\Omega\Gamma(l_v,f_\varepsilon)\rho\,dx
     &  %= \int_\Omega (\partial_1 f_\varepsilon)\rho\,dx
     = \int_V (\partial_1 f_\varepsilon)\rho\,dx
     = \int_{B_r^{d-1}}dy\int_{a_1}^{b_1}
            \partial_1 f_\varepsilon(t,y)\rho(t,y)\,dt
\\&
     = \int_{B_r^{d-1}}\Bigg(f_\varepsilon(b_1,y)\rho(b_1,y)-\int_{a_1}^{b_1}
         f_\varepsilon(t,y)\partial_1\rho(t,y)\,dt\Bigg)\,dy
\\&
     = \Bigg(\int_{B_r^{d-1}}\rho(b_1,y)\,dy\int_{\bR}\varphi_r(t,y)\,dt\Bigg)
          - \int_V f_\varepsilon \partial_1\rho\,dx
\\&
     \ge \min_{|y|<r}\rho(b_1,y)
        -\int_V f_\varepsilon \partial_1\rho\,dx.
\end{aligned}
$$
On the other hand, we have $0\le f_\varepsilon\le r^{-d}$,
$\Omega\cap\text{supp}(f_\varepsilon)\subset(\Omega\setminus\Omega_{2\varepsilon})\cap V$,
and the Lebesgue measure of this set tends to $0$ as
$\varepsilon\to 0$. Hence
$$
    \left|\int_\Omega f_\varepsilon\big(\LL(l_v)-\partial_1\rho\big)\,dx\right|
     \le r^{-d}\int_{(\Omega\setminus\Omega_{2\varepsilon})\cap V}
            \big|\LL(l_v)-\partial_1\rho\big|\,dx
    \to 0 \quad\text{as $\varepsilon\to 0$},
$$
thus, setting $C=\min_{|y|<r}\rho(b_1,y)$, we obtain
$$
   \int_\Omega\Big( f_\varepsilon\LL(l_v) + \Ga(f_\varepsilon,l_v)\Big)\rho\,dx
    \ge \int_\Omega f_\varepsilon\LL(l_v)\rho\,dx + C
     - \int_\Omega f_\varepsilon\partial_1\rho\,dx
    \underset{\varepsilon\to0}\to C
$$
which contradicts Lemma~\ref{lem.eq.IPP3} because $C>0$.
\epf

\bigskip

 \brmq\label{rmq.compute.g}
The equation $\hat \Delta=0$ of the boundary being given,   the problem of finding a symmetric matrix $(g^{ij})(x)$ formed with second degree polynomials and first degree polynomial $S^i$ such that
$$
   \sum_j g^{ij} \partial_j \hat \Delta= S^i\hat\Delta
$$
is a linear problem in the coefficients of  $g^{ij}$ and $S^i $
(there are $d(d+1)^2(d+2)/4+ d(d+1)$ such coefficients)
which can be easily solved for small $d$.

In the case when each $\Delta_k=0$ is a rational hypersurface, i.e., it can be
parametrized by rational functions (this is the case in all known so far solutions
of the DOP problem), it could be more convenient to compute the coefficients of
the $g^{ij}$ from the boundary condition rewritten in the form
$$
  \sum_j g^{ij}\, dx_j^*=0 \qquad\text{on}\qquad \{\hat\Delta_k=0\}
$$
(cf.~Remark \ref{rmq.thm.gal}). This is also a system of linear equations on the coefficients
of the $g^{ij}$. For example, when $d=2$ and $x_1=\xi_1(t)$, $x_2=\xi_2(t)$ is
a rational parametrization of $\Delta_k=0$, we need to equate to zero
the coefficients of all powers of $t$ in the numerators of
$$
   g^{11}\dot\xi_2(t) - g^{12}\dot\xi_1(t)
\qquad\text{and}\qquad
   g^{21}\dot\xi_2(t) - g^{22}\dot\xi_1(t).
$$
\ermq

\medskip

 \brmq\label{rmq.compute.rho}
As soon as the matrix $(g^{ij})$ is known, all the admissible measure densities $\rho$
can be found as follows.The conditions \eqref{eq.rho.1} yield a system of
linear equations for the unknown parameters in \eqref{log.rho} (the
coefficients of $Q$ and the numbers $\gamma_q$). Then it remains to
select those solutions which satisfy the integrability conditions.
In dimension $2$, if $\Omega$ and $\rho$ are given in some local curvilinear
coordinates by $x^p>y^2$ and $(x^p-y^2)^\alpha f(x,y)$ (with $f(0,0)\ne0$)
respectively, then
the integrability condition reads $\alpha>-\tfrac{1}{2}-\tfrac{1}{p}$.
If $\Omega$ and $\rho$ are given locally by $0<y<x^p$ and
$y^\alpha(x^p-y)^\beta f(x,y)$ respectively, then it
reads $\alpha+\beta > -1-\frac{1}{p}$. Note that only these singularities
occur in our classification in dimension $2$. 
 \ermq

\medskip

\brmq\label{rmq.soukhanov}
When the boundary equation has maximal degree $2d$, then it is proportional to the determinant of the metric $\Delta$.  In this case, if $\Delta^{-1/2}$ is integrable on the domain , then the Laplace-Beltrami  operator associated with the co-metric $g$ is a solution of the DOP problem on $\Omega$.  It turns out  that in any example where  it is the case, the associated curvature (in  dimension $2$ the scalar curvature) is constant, and even either $0$ either positive.
 Lev Soukhanov recently proved that in the general case, whenever the boundary has
 maximal degree $2d$, the associated metric is the product of Einstein
 metrics~\cite{Soukhanov2014}, and it is locally homogeneous, i.~e., any two points
 have isometric neighbourhoods~\cite{Soukhanov2017}.
 The latter fact is proven in~\cite{Soukhanov2017} when polynomials are
 ordered by any weighted degree.
\ermq

%%%%%%%%%%%%%%%%%%%%%%%%%%%%%%%%%%%%%%%%%%%%%%%%%%%%%%%%%%%%%%%%%%%%%%%%
%%%%%%%%%%%%%%%%%%%%%%%%%%%%%%%%%%%%%%%%%%%%%%%%%%%%%%%%%%%%%%%%%%%%%%%%
%%%%%%%%%%%%%%%%%%%%%%%%%%%%%%%%%%%%%%%%%%%%%%%%%%%%%%%%%%%%%%%%%%%%%%%%
%%%%%%%%%%%%%%%%%%%%%%%%%%%%%%%%%%%%%%%%%%%%%%%%%%%%%%%%%%%%%%%%%%%%%%%%
%%%%%%%%%%%%%%%%%%%%%%%%%%%%%%%%%%%%%%%%%%%%%%%%%%%%%%%%%%%%%%%%%%%%%%%%
%%%%%%%%%%%%%%%%%%%%%%%%%%%%%%%%%%%%%%%%%%%%%%%%%%%%%%%%%%%%%%%%%%%%%%%%
%%%%%%%%%%%%%%%%%%%%%%%%%%%%%%%%%%%%%%%%%%%%%%%%%%%%%%%%%%%%%%%%%%%%%%%%
%%%%%%%%%%%%%%%%%%%%%%%%%%%%%%%%%%%%%%%%%%%%%%%%%%%%%%%%%%%%%%%%%%%%%%%%
%%%%%%%%%%%%%%%%%%%%%%%%%%%%%%%%%%%%%%%%%%%%%%%%%%%%%%%%%%%%%%%%%%%%%%%%
%%%%%%%%%%%%%%%%%%%%%%%%%%%%%%%%%%%%%%%%%%%%%%%%%%%%%%%%%%%%%%%%%%%%%%%%
%%%%%%%%%%%%%%%%%%%%%%%%%%%%%%%%%%%%%%%%%%%%%%%%%%%%%%%%%%%%%%%%%%%%%%%%

\section{The bounded solutions in dimension 2\label{dim2}}

In this section, we concentrate on the DOP problem in dimension 2 for bounded domains.
The central result of this section is the following

\bthm\label{thm.central}
 In $\bR^2$, up to affine transformations, there are exactly
10 relatively compact sets
and a one-parameter family
for  which there exists a solution for the DOP problem: the triangle, the square, the disk, and the areas bounded by two co-axial parabolas, by one parabola and two tangent lines, by one parabola, its axis, and a tangent line,
% the axis of the parabola,
by the nodal cubic $y^2=x^2+x^3$, by the cuspidal cubic $y^2=x^3$ and one tangent,
by the cuspidal cubic $y^2=x^3$ and the
 vertical line $x=1$,
%(that is a line passing through the infinite point of the cubic),
by a swallow tail,
or by a deltoid curve (see Section~\ref{comments.models} for more details).
\ethm

This theorem is an immediate consequence from
Propositions~\ref{RAlgDOP.quartic}, \ref{CAlgDOP.cubic},
\ref{AlgDOP.circle}, \ref{AlgDOP.parab}, and \ref{prop.only.lines}.
Since we look at bounded domains, we may therefore reduce to the SDOP problem, and  we solve the  algebraic problem described in Section~\ref{cas.general} in the particular case of dimension 2. For  basic references on plane algebraic curves
and their singularities, see \cite{Brieskorn.Knoerrer},
\cite[Ch.~I, \S3]{GLS}. \cite{walk}.

In the following definition we restrict ourselves by dimension 2, but
it can be obviously extended to any dimension.

\bdefi[AlgDOP problem]\label{gal.pb.compact} Let $\bK$ be $\bR$ or $\bC$ and
let $a,b,c,\hat\Delta$ be polynomials in $\bK[x,y]$.
We say that $(a,b,c,\hat\Delta)$ is a solution of the Algebraic counterpart
of the DOP problem over $\bK$ ($\,\bK$-AlgDOP problem for short), if
$a$, $b$, and $c$ are of degree at most $2$, the polynomial $\Delta:=ac-b^2$
is not identically zero,
and $\hat\Delta$ is a square-free
polynomial which divides each of the following three polynomials:
$$
   \Delta,  \qquad
   a\partial_1\hat\Delta+b\partial_2\hat\Delta, \qquad
   b\partial_1\hat\Delta+c\partial_2\hat\Delta.
$$
\edefi

Due to Theorem \ref{thm.gal.bord.Omega}, if $(\Omega,g,\rho)$ is a solution
to the DOP problem and with a bounded $\Omega$, then
$(g^{11},g^{12},g^{22},\hat\Delta)$ is a solution to the $\bR$-AlgDOP problem
where $\hat\Delta=0$ is the minimal equation of $\partial\Omega$.
So, our strategy is to find all solutions to the $\bC$-AlgDOP problem up
to affine linear transformations of $\bC^2$, then to find all solutions
to the $\bR$-AlgDOP problem such that $\bR^2\setminus\{\hat\Delta=0\}$
has a bounded component, and eventually to find all possible mesures $\rho$.

It is clear that the condition that $\hat\Delta$ divides
$a\partial_1\hat\Delta+b\partial_2\hat\Delta$ and
$b\partial_1\hat\Delta+c\partial_2\hat\Delta$ is equivalent to
the condition that for each irreducible factor $\Delta_1$ of $\hat\Delta$
one has
\begin{eqnarray}
\label{eqbord1}a\partial_1\Delta_1+b\partial_2\Delta_1=L_1\Delta_1\\
\label{eqbord2}b\partial_1\Delta_1+c\partial_2\Delta_1=L_2\Delta_1
\end{eqnarray}
where $\deg L_i\leq1, i=1,2$.
The equations \eqref{eqbord1}--\eqref{eqbord2}, in their turn, being equivalent
to
\beq
    \label{eqbordparam}
    a(\xi,\eta)\dot\eta = b(\xi,\eta)\dot\xi, \qquad
    b(\xi,\eta)\dot\eta = c(\xi,\eta)\dot\xi
\eeq
for any local analytic branch $x=\xi(t)$, $y=\eta(t)$ of the curve $\Delta_1=0$
(since $\Delta_1$ is irreducible, the equalities \eqref{eqbordparam}
for an arbitrary local branch of $\Delta_1=0$ imply the same equalities for
all local branches of $\Delta_1=0$).

The proof of Theorem \ref{thm.central}
is divided in many parts.  In Section~\ref{subsec.newton}, we prove that the curves
$\{\hat \Delta=0\}$ may  have flex or planar points at infinity only 
(Lemma~\ref{lemNoFlex}), unless $\Delta$ is reducible. We also describe the various
singularities which may occur at finite distance (Corollary~\ref{corBranch}) and 
the behavior at infinity (Lemma~\ref{lemTwoBranch}).  Section~\ref{dual} studies
the case where $\Delta$ is irreducible of degree $4$,
while the Sections~\ref{subsec.cubic}--\ref{subsec.linear}
concentrate on the reducible case.

%%%%%%%%%%%%%%%%%%%%%%%%%%%%%%%%%%%%%%%%%%%%%%%%%%%%%%%%%%%%%
%%%%%%%%%%%%%%%%%%%%%%%%%%%%%%%%%%%%%%%%%%%%%%%%%%%%%%%%%%%%%

\subsection{A preliminary study of Newton polygons of $a,b,c$ and $\Delta$\label{subsec.newton}}

Let
$(a,b,c,\hat\Delta)$ be a solution to the $\bC$-AlgDOP problem, $\Delta=ac-b^2$,
and $\Delta_1$ be an irreducible factor of $\hat\Delta$
which is not a common factor of $a$, $b$, and $c$. Note that the last condition
is always satisfied when $\Delta$ has an irreducible component of degree $\ge3$.

We shall use projective coordinates $(X:Y:Z)$ such that
$x=X/Z, y=Y/Z$ and denote $L_{\infty}$ the line $Z=0$ in $\CP^2$.
Let $\gamma$ be an {\it analytic branch} of the curve $\Delta_1=0$ at
some finite or infinite point, i.e.
$\gamma$ is a germ at $0$ of a non-constant meromorphic mapping $\bC\to\bC^2$,
$t\mapsto(\xi(t),\eta(t))$ such that $\Delta_1(\xi(t),\eta(t))=0$.
Let $v_\gamma:\bC[x,y]\to\bZ\cup\{\infty\}$ be the corresponding valuation,
i.e. $v_\gamma(f) = \ord_t f\big(\xi(t),\eta(t)\big)$ where
$$
	\ord_t u(t) =\begin{cases} n
        &\text{if $u(t)=\sum_{k\ge n}u_k t^k$ and
					$u_n\ne0$,}\\
			\infty	&\text{if $u(t)=0$}
	\end{cases}
$$
We denote $p=v_\gamma(x)=\ord\xi$ and $q=v_\gamma(y)=\ord\eta$.

\blem~\label{lemValu}\\ (a) Suppose that none of $\xi(t)$, $\eta(t)$ is constant. 
Then
\beq\label{eqVabc}
	v_\gamma(a)-v_\gamma(b) = v_\gamma(b)-v_\gamma(c)
		= \ord_t\dot\xi - \ord_t\dot\eta.
							\eeq
(b) Suppose that $\eta(t)$ is constant. Then $v_\gamma(b)=v_\gamma(c)=\infty$, i.e.,
$b$ and $c$ vanish identically on $\gamma$.
\elem

\bpf
\par
(a) By \eqref{eqbordparam},
 both $(a,b)$ and $(b,c)$ are proportional to $(\dot \xi, \dot \eta)$.
 Then, let us show  that no one of the coefficients $a,b$ and $c$ vanishes
 identically along $\gamma$.  
 Indeed if one vanishes then so will do the other ones because of this
 proportionality. Then $\Delta_1$ divides $a$, $b$, and $c$
 which contradicts our assumption about $\Delta_1$.

Then, again by \eqref{eqbordparam}, we have
$$
  v_\gamma(a)+\ord\dot\eta=v_\gamma(b)+\ord\dot\xi, \quad
  v_\gamma(b)+\ord\dot\eta=v_\gamma(c)+\ord\dot\xi.
$$

\smallskip
(b) Straightforward from the proportionality of $(a,b)$ and $(b,c)$ to $(\dot \xi, 0)$.
\epf

\medskip\par
As usually, for a polynomial $u=\sum u_{kl}x^ky^l$, we define its
Newton polygon $\cN(u)$ as the convex hull in $\bR^2$ of the set
$\{(k,l)\,|\,u_{kl}\ne0\}$.

Recall that $p=v_\gamma(x)$, $q=v_\gamma(y)$.
We have $v_\gamma(x^ky^l) = L_\gamma(k,l)$ where 
$L_\gamma$ is the linear
form $L_\gamma(r,s) = pr+qs$.
Thus, for any polynomial $u(x,y)$ we have 
$v_\gamma(u)\ge \min_{\cN(u)} L_\gamma$
and if the minimum of $L_\gamma$ is attained at a single vertex of $\cN(u)$,
then then  $v_\gamma(u) = \min_{\cN(u)} L_\gamma$

The notation of the style $b=\NND101100$ (any combination of 
$\circ$ and $\bullet$) means that $b$ is a linear combination
of monomials corresponding to the $\bullet$'s.
For example, $b=\NND101100$ means that $b_{01}=b_{10}=b_{20}=0$
(the coefficients of $y$, $x$, and $x^2$) and the other coefficients
may or may not be zero.

In the following lemma, we look for restrictions on Newton polygons of
$a$, $b$, and $c$ imposed by the fact that $(\xi,\eta)$ has a given valuation $(p,q)$.
The cases $p$ or $q$ negative correspond to points at infinity.

\par

\blem~\label{lemBranch}\\
(a)
If $(p,q)=(1,2)$, then 
$b=\NND111011$ and $c=\NND111001$.\\
(b)
If $(p,q)=(1,3)$, then
$b=\NND111001$ and $c=\NND101000$, in particular, \newline
$\mult_{(0,0)}\Delta\ge2$.
\smallskip
\\
(c) $(p,q)=(1,4)$ is impossible.
\\
(d)
If $(p,q)=(-1,0)$, then $b=\NND111110$, $c=\NND111110$, and $c_{10}+c_{11}=0$.
\\
(e)
If $(p,q)=(-1,1)$,  then $b=\NND111100$, $c=\NND111100$, and $c_{00}+c_{01}=0$.
\\
(f).
$p=-1$ and $2\le q<\infty$ is impossible.
\\
(g)
If $(p,q)\in\{(-2,-1),\, (-3,-2),\, (-4,-3)\}$, then 
$b=\NND111110$ and $c=\NND110110$.
\\
(h)
If $(p,q)=(-2,1)$, then
$b=\NND111100$ and $c=\NND100000$.
\\
(i) $(p,q)=(3,4)$ is impossible.
 
 \elem

\bpf
(a)
If $(p,q)=(1,2)$,  then $v_\gamma(\dot\xi)=0$ and $v_\gamma(\dot\eta)=1$.
Hence, by \eqref{eqVabc} we have 
$v_\gamma(b)=v_\gamma(a)+1\ge1$ and $v_\gamma(c)=v_\gamma(a)+2\ge2$
and the result follows from the fact that $v_\gamma(1)=0$, $v_\gamma(x)=1$,
and $v_\gamma(x^ky^l)\ge 2$ when $(k,l)\not\in\{(0,0),(0,1)\}$.

\smallskip
(b)
If $(p,q)=(1,3)$, then $v_\gamma(\dot\xi)=0$ and $v_\gamma(\dot\eta)=2$.
Hence, by \eqref{eqVabc} we have 
$v_\gamma(b)=v_\gamma(a)+2\ge2$ and
$v_\gamma(c)=v_\gamma(b)+2=v_\gamma(a)+4\ge4$.
The values of $v_\gamma$ on the monomials of degree $\le2$ are:
\beq\label{eqMonom}
	v_\gamma(1)=0,\;\;
	v_\gamma(x)=1,\;\;
	v_\gamma(x^2)=2,\;\;
	v_\gamma(y)=3,\;\;
	v_\gamma(xy)=4,\;\;
	v_\gamma(y^2)=6.				
\eeq
Thus, $v_\gamma(b)\ge2$ implies $b_{00}=b_{01}=0$
and $v_\gamma(c)\ge4$ implies $c_{00}=c_{10}=c_{20}=c_{01}=0$.
In particular, $\mult_{(0,0)}b\ge 1$ and $\mult_{(0,0)}c\ge2$.
Hence, $\mult_{(0,0)}(b^2-ac)\ge2$

\smallskip
(c) 
If $(p,q)=(1,4)$, then $v_\gamma(\dot\xi)=0$ and $v_\gamma(\dot\eta)=3$.
Hence, by \eqref{eqVabc} we have 
\beq\label{eqBrProof}
	v_\gamma(c)-v_\gamma(b) = v_\gamma(b)-v_\gamma(a) = 3
\eeq
The values of $v_\gamma$ on the monomials of degree $\le2$ are:
$$
	v_\gamma(1)=0,\;\;
	v_\gamma(x)=1,\;\;
	v_\gamma(x^2)=2,\;\;
	v_\gamma(y)=4,\;\;
	v_\gamma(xy)=5,\;\;
	v_\gamma(y^2)=8.
$$
Hence, we have $\{v_\gamma(a),v_\gamma(b),v_\gamma(c)\}
\subset\{0,1,2,4,5,8\}$.
Under this condition, \eqref{eqBrProof} is possible only for
$v_\gamma(a)=2$, $v_\gamma(b)=5$, $v_\gamma(c)=8$, hence
$a=\NND111001$, $b=\NND101000$, and $c=\NND100000$.
It follows that $\Delta=y^2 f(x,y)$.
This is impossible because $\gamma$ cannot be
a branch of a polynomial of degree $\le2$.

\smallskip
(d) $(p,q)=(-1,0)$. 
If $\dot\eta=0$, we use Lemma \ref{lemValu}(b).
Otherwise the proof is similar to (a)--(c). Indeed,
we have $\ord_t\dot\xi=-2$ and $\ord_t\dot\eta\ge 0$, thus
$v_\gamma(c)-v_\gamma(b)=v_\gamma(b)-v_\gamma(a)\ge 2$ by
\eqref{eqBrProof}. We have $v_\gamma(x^ky^l)=-k$, thus $v_\gamma(a)\ge-2$, hence
$v_\gamma(c)>v_\gamma(b)=v_\gamma(a)+2\ge 0$.
Therefore $b_{20}=c_{20}=0$
 (otherwise $v_\gamma(b)$ or $v_\gamma(c)$ would be $-2$)
and $b_{10}+b_{11}=c_{10}+c_{11}=0$
 (otherwise $v_\gamma(b)$ or $v_\gamma(c)$ would be $-1$).

\smallskip
(e) 
If $(p,q)=(-1,1)$, then by \eqref{eqVabc} we have
$v_\gamma(c)-v_\gamma(b)=v_\gamma(b)-v_\gamma(a)=2$, hence
$v_\gamma(a)\ge -2$, $v_\gamma(b)\ge 0$, $v_\gamma(c)\ge 2$ and
the result follows (as in point (d), $c_{00}+c_{11}=0$ because
otherwise we would have $v_\gamma(c)=0$).

\smallskip
(f) 
We have 
$v_\gamma(c)-v_\gamma(b) = v_\gamma(b)-v_\gamma(a) = q+1$
and
$v_\gamma(
x^2, x,1, xy,y,y^2)=
(-2,-1,0,q-1,q, 2q)$.
Thus, $v_\gamma(a,b,c)=(-2,q-1,2q)$, i.e.
$b=\NND111000$ and $c=\NND100000$,
Therefore,
$\Delta=y^2 f(x,y)$. This is impossible because $\gamma$ cannot be
a branch of a polynomial of degree $\le2$.

\smallskip
(g,h)  The proof is similar to the previous cases.

\smallskip
(i) 
We have $v_\gamma\{a,b,c\}\subset
v_\gamma\{1,x,y,x^2,xy,y^2\}=\{0,3,4,6,7,8\}$.
Combining this with
$v_\gamma(c)-v_\gamma(b)=v_\gamma(b)-v_\gamma(a)
=\ord_t\dot\eta-\ord_t\dot\xi=3-2=1$, we obtain
$v_\gamma(a)=6$, $v_\gamma(b)=7$, $v_\gamma(c)=8$, i.e.
$a=\NND101001$, $b=\NND101000$, $c=\NND100000$.
Thus, $\mult_0(\Delta)=4$, i.e., $\Delta=0$ is a union of
four lines which contradicts the condition $(p,q)=(3,4)$.
\epf
\medskip

According to the standard terminology (see, e.g.~\cite[Ch.~I, \S2.4]{GLS}),
we say that an analytic branch $\beta$
of an algebraic curve in $\CP^2$ is {\it generic}, {\it flex}, {\it planar},
or has singularity of type $A_2$ (called also {\it cusp}) or $E_6$ if
there exists an affine coordinate chart $(u,v)$ such that the pair
$(p,q):=v_\gamma(u,v)$ is as in the second column of Table~1 in
Section~\ref{dual}.

\bcor\label{corBranch}~\\
(a). $\Delta$ cannot have a singularity of type $E_6$
at a finite point.\\
(b). Suppose that $\gamma$ is a singular branch of $\Delta_1$
of type $A_2$ at a point $P\in L_\infty$ and
$L_\infty$ is not tangent to $\gamma$ at $P$. Then
there is another branch of $\Delta$ at $P$, or $\deg\Delta=3$.
\ecor

\bpf
\par
The point (a) follows from Lemma \ref{lemBranch}(i).
\par
Point (b) corresponds to $(p,q)=(-2,1)$
for a suitable choice of the coordinates
(whereas $(-3,-1)$
corresponds to a cusp on $L_\infty$ tangent to $L_\infty$).
We are therefore in case Lemma \ref{lemBranch}(h). 
Hence $b(x,0)=b_{00}$ (a constant) and $c(x,0)=0$ whence
$\Delta(x,0)=b_{00}^2$. This means that the local intersection of
the line $\{y=0\}$ with $\{\Delta=0\}$ at $P$ is equal to $\deg\Delta$.
On the other hand, the local intersection of this line with $\gamma$ at $P$ is $3$,
thus either $\Delta$ has another branch at $P$ or $\deg\Delta=3$.
\epf

\smallskip

\blem\label{lemNoFlex} 
Let $\gamma$ be a flex or planar branch of $\Delta_1$ at $P$. Then\\
(a) if $P\in L_\infty$, then 
$\gamma$ is tangent to $L_\infty$. \\
(b)
if $P\not\in L_\infty$, then $\gamma$ is not planar and
$\mult_P(\Delta)>\mult_P(\hat\Delta)$,
in particular, this is impossible when $\Delta$ is irreducible.\\
\elem

\bpf
(a) Follows from Lemma \ref{lemBranch}(f). 

\smallskip
(b)
The fact that $\gamma$ is not planar follows from Lemma \ref{lemBranch}(c).
Let us choose affine coordinates so that $P$
is the origin and the axis $y=0$ is tangent to $\gamma$.
Thus, $\ord_t\gamma=(1,3)$. By Lemma \ref{lemBranch}(b),
we have $\mult_{(0,0)}\Delta\ge2$, i.e.,
there is another branch $\beta$ of $\Delta$ passing through the origin.
Since $\deg_x\Delta(x,0)\le4$ the multiplicities of the intersection 
of the axis $y=0$ with $\gamma$ and $\beta$ are $3$ (i.e., $q=3$) and $1$
respectively.

It remains to prove that $\beta$ cannot be a branch of $\hat\Delta$. Suppose it is.
Let us choose coordinates so that the axis $x=0$ is tangent to $\beta$.
Then Lemmas~\ref{lemValu}(b) and \ref{lemBranch}(a,b) applied to $\beta$ imply
\beq
        \label{eqNoFProofOne}
	a_{00}=a_{01}=0.		
\eeq
(we swap $x\leftrightarrow y$ and $a\leftrightarrow  c$ in
Lemma \ref{lemBranch}). Hence, $v_\gamma(a)\ge 1$.
By (\ref{eqVabc}), we have 
\beq
        \label{eqNoFProofTwo}
	v_\gamma(c)-v_\gamma(b) = v_\gamma(b)-v_\gamma(a) = 2
\eeq
(see the proof of Lemma \ref{lemBranch}(b)).
Recall that the values of $v_\gamma$ on monomials
are given by (\ref{eqMonom}).
Hence, $\{v_\gamma(a),v_\gamma(b),v_\gamma(c)\}\subset\{0,1,2,3,4,6\}$.
Combining this with (\ref{eqNoFProofTwo}) and $v_\gamma(a)\ge 1$, we obtain
$v_\gamma(a)=2$, $v_\gamma(b)=4$, $v_\gamma(c)=6$.
By \eqref{eqNoFProofOne}, this implies 
$\mult_{(0,0)}(a)=\mult_{(0,0)}(b)=\mult_{(0,0)}(c)=2$, hence
$\mult_{(0,0)}(\Delta)=4$.
This means that $\Delta$ is a union of four lines passing through the origin.
Contradiction.
\epf

\blem\label{lemTwoBranch}
Let $\gamma$ be a smooth
branch of $\Delta_1$ at $P\in L_\infty$.
Suppose that there
exists a line $L$ passing through $P$ which
is tangent to a branch $\beta$ of $\hat\Delta$ at a finite point $Q$. 
Suppose also that $\beta,\gamma\not\subset L$.
Then\\
\smallskip
(a) $\beta$ is smooth at $Q$.\\
\smallskip
(b) $\mult_P(\Delta)\ge 2$ or $\deg\Delta\le 3$.

\elem

\bpf
Let us choose coordinates so that $L$ is the axis $y=0$ and 
and $\beta$ is tangent to $L$ at the origin.
Then all possibilities for $\gamma$ are covered by
Lemma \ref{lemBranch}(d)--(g) and in all these cases we have
$b,c=\NND111110$, i.e., $b_{20} = c_{20}=0$.

\smallskip
(a)
Let $\beta=(\xi,\eta)$ and $(p,q)=\ord_t(\xi,\eta)$.
Suppose that $\beta$ is singular. Then $\min(p,q)\ge2$.
We have also
$q>p$ (because $L$ is tangent to $\beta$) and
$q=(L.\beta)\le3$ (because $(L.\beta)+(L.\gamma)\le 4$). Thus,
$(p,q)=(2,3)$ hence, by (\ref{eqVabc}), we have
$v_\beta(c)-v_\beta(b)=v_\beta(b)-v_\beta(a)=1$.
Combining this fact with
$v_\beta(1,x,y,x^2,xy,y^2)=(0,2,3,4,5,6)$ and
$b_{20}=c_{20}=0$, we obtain $v_\beta(a,b,c)=(4,5,6)$, i.e.,
$a=\NND101001$, $b=\NND101000$, $c=\NND100000$.
Thus, $\Delta$ is homogeneous. A contradiction.

\smallskip
(b)
Combining $b_{20} = c_{20} = 0$ with Lemma \ref{lemBranch}(a) applied to $\beta$,
we obtain $b=\NND111010$, $c=\NND111000$.
Thus, it is enough to show that $c_{11}=0$.
Indeed, if $\gamma$ is tangent to $L_\infty$, this is already proven
in Lemma \ref{lemBranch}(g).
Otherwise by Lemma \ref{lemBranch}(d,e,f) we have
$c_{00}+c_{11}=0$ or $c_{10}+c_{11}=0$ and we know that $c_{00}=c_{10}=0$.
\epf

\subsection{The duals of quartic curves}\label{dual}

%We denote the {\it dual projective plane} by $\check\mP^2$.
Let $C$ be an irreducible algebraic curve in $\mathbb P^2$ of degree $d\geq 2$. 
Let $\check C$ be the dual curve in $\check\mP^2$ is the set of all
lines in $\mP^2$ endowed with the natural structure of the projective
plane, and $\check C$ is the set of all lines in $\mP^2$ which are
tangent to $C$. 

If $t\to\gamma(t)$ is a local analytic branch of $C$, then we
denote the {\it dual branch} of $\check C$ by $\check\gamma$.
It is defined by $t\mapsto\check\gamma(t)$ where $\check\gamma(t)$
is the line which is tangent to $C$ at $\gamma(t)$.

Let $\gamma$ be a local branch of $C$. Let us choose affine coordinates
$(X,Y)$ so that $\gamma$ is given by $X=\xi(t)$, $Y=\eta(t)$,
$\xi(0)=\eta(0)=0$. Then the equation of the line $\check\gamma(t)$
is $(X-\xi)\dot\eta - (Y-\eta)\dot\xi = 0$.
Thus, in the standard homogeneous
coordinates on $\check\mP^2$ corresponding to the coordinate chart $(X,Y)$,
the dual branch $\check\gamma$ has a parametrization of the form
\beq\label{eqDualSing}
    t\mapsto (\dot\eta\,:\, -\dot\xi\,:\, \dot\xi\eta - \xi\dot\eta)
								\eeq
and we obtain the following fact.

\blem\label{lemDualBranch}
Let $\gamma$ be a local branch of $C$ 
and $\check\gamma$ the dual branch
of $\check C$. Let $(X,Y)$ be an affine chart such that
$\gamma$ has the form $X=\xi(t)$, $Y=\eta(t)$ with $0<p<q$ where
$p=\ord_t\xi$ and $q=\ord_t\eta$. Then, in suitable affine coordinates
$(\check X,\check Y)$ on $\check\mP^2$, the branch $\check\gamma$
has the form $\check X=\check\xi(t)$, $\check Y=\check\eta(t)$
with $\ord_t\check\xi=q-p$ and $\ord_t\check\eta=q$. 
\elem

For a point $P\in C$, we denote the delta-invariant of $(C,P)$
by $\delta_P$ or $\delta_P(C)$ (see \cite[p.~206]{GLS}).
Informally speaking, $\delta_P$ is the number
of double points of $C$ concentrated in $P$. We have 
(see \cite[Thm.~10.2]{milnor} or \cite[Ch.~I, Prop.~3.34]{GLS})
$$
	2\delta_P = \mu+r-1 = \sum m_i(m_i-1)
$$
where $\mu$ is the Milnor number and $r$ is the number
of local branches of $C$ at $P$, and
$\mathbf{m}=[m_1,m_2,\dots]$ is the sequence of the multiplicities of all infinitely
near points of $P$.
If $P$ is a non-singular point of $C$, then $\delta_P=0$.
It easily follows from the definition
(see also \cite[Ch.~I, Lem.~3.32]{GLS}) that 
\beq\label{eqDeltaZ}
	\delta_P=\delta_P^0 + \sum_{i=1}^r \delta(\gamma_i)
	\qquad\text{where}\quad \delta_P^0 = \sum_{1\le i<j\le r}
		(\gamma_i\cdot\gamma_j),			
\eeq
$\gamma_1,\dots,\gamma_r$ are local branches of $C$ at $P$, and
$(\gamma_i\cdot\gamma_j)$ is the intersection number of $\gamma_i$ and $\gamma_j$
at $P$.

Let $g$ be the genus of $C$. By the genus formula
(see \cite[p.~624, Thm.~7]{Brieskorn.Knoerrer} or \cite[\S10, Eq.~(1)]{milnor}),
we have
\beq\label{eqGenus}
	2g + 2\sum_{P\in C}\delta_P = (d-1)(d-2).		\eeq
Combining (\ref{eqDeltaZ}) and (\ref{eqGenus}), we obtain
\beq\label{eqGenusZ}
	2g + 2n + 2\sum_\gamma\delta(\gamma) = (d-1)(d-2)	
\eeq
where $\gamma$ runs over all local branches of $C$ at all points
and  $n=\sum_{P\in C}\delta^0_P$
(only a finite number of terms in the both sums are non-zero).

For a local branch $\gamma$ of a curve $C$ at a point $P$, we
denote the multiplicity of $\gamma$ at $P$ by $m(\gamma)$.
If $\gamma$ is parametrized by $X=\xi(t)$, $Y=\eta(t)$
in some local coordinates $X,Y$, then $m(\gamma)=\min(\ord_t\xi,\ord_t\eta)$.
We set also $\bEps(\gamma) = 2\delta(\gamma)+m(\gamma)-1$.
Let $\check d$ be the degree of $\check C$.
In this notation, the first Pl\"ucker formula (the class formula)
takes the form (see \cite[Thm.~1.3]{Kulikov2016})
\beq\label{eqPlu}
	\check d = d(d-1) - 2n - \sum_\gamma \bEps(\gamma)	
\eeq
and the second Pl\"ucker
formula (the Riemann-Hurwitz formula for a generic
projection of $\check C$ onto a line) is
\beq\label{eqRH}
	2-2g = 2\check d - d - \sum_\gamma (m(\check\gamma) - 1).	
\eeq
In the both formulas $\gamma$ runs over all local branches of $C$.

If $d=4$, then $\sum\delta(\gamma)\le3$ by (\ref{eqGenusZ}) which is possible
for the sequences of multiplicities $[2]$, $[2,2]$, $[2,2,2]$, and $[3]$ only,
hence all singular branches
are of the types $A_2$, $A_4$, $A_6$ and $E_6$ (recall that $A_k$ and $E_6$
are given by $v^2=u^{k+1}$ and $v^3=u^4$ is suitable curvilinear local coordinates).
In Table 1, we list all types of local branches $\gamma(t)=(\xi(t),\eta(t))$,
$\ord_t\xi=p$, $\ord_t\eta=q$, $p<q$,
and their invariants contributing to (\ref{eqGenusZ}), (\ref{eqPlu}), and (\ref{eqRH})
(we use  Lemma \ref{lemDualBranch} to compute $\check p=\ord_t\check\xi$ and
$\check q=\ord_t\check\eta$).
\label{ListBranches}
\def\TL#1#2#3#4#5#6#7{\noindent
	\hbox to 35mm{\hfill#1\hfill}
	\hbox to 15mm{\hfill#2\hfill}
	\hbox to 12mm{\hfill#3\hfill}
	\hbox to 12mm{\hfill#4\hfill}
	\hbox to 12mm{\hfill#5\hfill}
	\hbox to 12mm{\hfill#6\hfill}
	\hbox to 15mm{\hfill#7\hfill}
}

\hbox to\hsize{\hfill Table 1\;\;\qquad\quad}

\medskip
\noindent\hrule
\smallskip
\TL{}{$\mathbf{m}$}{$(p,q)$}{$(\check p,\check q)$}
	{$\delta(\gamma)$}{$\bEps(\gamma)$}{$m(\check\gamma)-1$}

\smallskip
\noindent\hrule
\smallskip

\TL{generic point}	- {(1,2)}{(1,2)}{0}{0}{0}

\TL{flex point}		- {(1,3)}{(2,3)}{0}{0}{1}

\TL{planar point}	- {(1,4)}{(3,4)}{0}{0}{2}

\TL{$A_2$}		    {[2]}{(2,3)}{(1,3)}{1}{3}{0}

\TL{$A_4$}	          {[2,2]}{(2,4)}{(2,4)}{2}{5}{1}

\TL{$A_6$}		{[2,2,2]}{(2,4)}{(2,4)}{3}{7}{1}

\TL{$E_6$}		{[3]}    {(3,4)}{(1,4)}{3}{8}{0}
\smallskip
\noindent\hrule
\bigskip

Thus, denoting the number of branches of the respective types by 
$f$ (flex), $p$ (planar), $a_2$, $a_4$, $a_6$, and $e_6$, we rewrite
(\ref{eqGenusZ}) -- (\ref{eqRH}) as
$$
\begin{matrix}
	g + n + a_2 + 2a_4 + 3a_6 + 3e_6 = 3,\\
	\check d = 12 - 2n - 3a_2 - 5a_4 - 7a_6 - 8e_6, \\
	2 - 2g = 2\check d - 4 - f - 2p - a_4 - a_6.
\end{matrix}
$$
Eliminating $g$ and $\check d$, we obtain
\beq\label{eqPluDegFour}
	f+2p = 24 - 8 a_2 - 15 a_4 - 21 a_6 - 22 e_6 - 6 n.
							\eeq
Since all the ingredients (including $g$) are non-negative, we
obtain the following fact.
\blem\label{lemOneFlex}
Suppose that $C$ is an irreducible quartic curve in $\mP^2$ which
has at most one smooth non-generic {\rm(}i.e., flex or planar{\rm)}
local branch. Then $C$ is rational {\rm(}i.e., $g=0${\rm)}
and one of the following cases occurs:

\bitem
\item$(i)$ {\rm(tricuspidal quartic)}
	$C$ has three singular points of type $A_2$
	and no smooth non-generic branches {\rm(}i.e., $f=p=0${\rm)}.
	The dual curve $\check C$ is a nodal cubic.
\item$(ii)$ {\rm(swallow tail)}
	$C$ has two singular points of type $A_2$ and
	one planar point, and one ordinary double point
        {\rm(}i.e., $f=0$, $p=n=1${\rm)}.
	The degree of $\check C$ is $4$, it has one singular point
	of type $E_6$ and two flex points. The equation of $\check C$
	in suitable affine coordinates is $y=x^4-x^2$.
\item$(iii)$
	Each of $C$ and $\check C$
	has two singular points of types $A_2$ and $A_4$ and
	one flex point {\rm(}i.e., $f=1$, $p=0${\rm)},
	the degree of $\check C$ is $4$,
\item$(iv)$
 	Each of $C$ and $\check C$ has one singular point of type $E_6$ and
	one planar point {\rm(}i.e., $f=0$, $p=1${\rm)}.
	The degree of $\check C$ is $4$. The equation of $\check C$
	in suitable affine coordinates is $y=x^4$.
\eitem
In each of the cases $(i)$--$(iv)$ the formulated conditions
uniquely determine the curve $C$ up to
automorphism of $\CP^2$.
\elem

\bpf

By (\ref{eqGenusZ}) we have $g+n+a_2+2a_4 + 3a_6 + 3e_6 = 3$.
 Substituting each nonnegative solution of this equation
into (\ref{eqPluDegFour}), we see that the
only cases when $f+p\le 1$ are:
\beqnas
  &(i)~\;\;\check d=3,\;a_2=3;\qquad\qquad\quad
 		&(iii)\;\;\check d=4,\;a_2=a_4=1,\; f=1,\\
 	 &(ii)~\;\;\check d=4,\;a_2=2,\; p=n=1,	&(iv)\;\;\check d=4,\;e_6=p=1.
\eeqnas
Let us show that these cases are uniquely realizable.
In cases $(ii)$ and $(iv)$ this follows from the fact
that $\check C$ has the singularity $E_6$, hence 
it has the equation $y=f(x)$, $\deg_xf=4$, in 
suitable coordinates. By affine changes of coordinates,
this equation reduces to $y=x^4$ or $y=x^4-x^2$.

In case $(i)$, the dual curve is a nodal cubic.
It is unique up to projective transformation, thus $C$
is also unique.

In case $(iii)$, let us choose homogeneous coordinates $(X:Y:Z)$
so that $A_2$ and $A_4$ are at $(0:0:1)$ and $(0:1:0)$
respectively and the lines $Y=0$ and $Z=0$ are tangent to $C$
at these points. (These are two distinct lines. Indeed, the local
intersection of $C$ with the tangent lines at $A_2$ and $A_4$ is
is $3$ and $4$ respectively, thus it cannot be a single line
by Bezout's theorem.)
Let $F(X,Y,Z)=0$ be the equation of $C$. Let us consider the Newton
polygon of the polynomial $F(X,Y,1)$.
The choice of the coordinates near $A_2$ ensures
that it is placed above the segment $[(0,2),(3,0)]$.
The choice of the coordinates near $A_4$ ensures that
the segment $[(0,2),(4,0)]$ is an edge of the Newton polygon.
Hence $F=u_{30}X^3Z + G$ where $G=u_{40}X^4+u_{21}X^2YZ+u_{02}Y^2Z^2$.
Moreover, the fact that $F$ has a single branch at $(0:1:0)$
implies that $G$ is a complete square. Hence, rescaling the coordinate,
we can obtain $u_{30}=u_{40}=u_{02}=1$, $u_{11}=2$, whence the uniqueness
up to projective change of coordinates.

\epf

\bcor \label{quartic}
Suppose that $C$ is an irreducible quartic curve in $\bC^2$ which
satisfies the restrictions imposed by Lemmas \ref{lemBranch}(h,i),
\ref{lemNoFlex}\ and
\ref{lemTwoBranch}, i.e.:
\bitem
\item
	any smooth non-generic branch of $C$
	is tangent to the infinite line $L_\infty$;
\item
	if $C$ meets $L_\infty$ transversally at a point $P$
        and $C$ is smooth at $P$, 
        then there is no line through $P$ {\rm(}except, 
	maybe, $L_\infty${\rm)} which is tangent to
	$C$ at a smooth or singular point;
\item
	if $C$ has a cusp $A_2$ at a point $P\in L_\infty$, 
	and $L_\infty$ is not the tangent to $C$ at $P$,
	then $C$ has another branch through $P$;
\item
	$C$ does not have a singularity of type $E_6$ at a finite point.
\eitem

\noindent
Then one of the cases $(i)$ or $(ii)$ of Lemma \ref{lemOneFlex}  occurs
and the position of $C$ with respect to the infinite line $L_\infty$
is one of:
\bitem
\item$(i_1)$
	{\rm(deltoid} or {\rm$(1,3)$-hypocycloid)}
	$L_\infty$ is the bitangent of $C$.

\item$(i_2)$
	$L_\infty$ is the tangent at a cusp.
	
\item$(ii)$
	{\rm(swallow tail)}
	$L_\infty$ is the tangent at the planar point.
\eitem
In each of the cases $(i_1)$ and $(ii)$
the affine curve $C$ is unique up to affine transformation of $\bC^2$.
In suitable affine coordinates, $C$ is parametrized by
\bitem
\item$(i_1)$
	$x=2\cos\theta + \cos2\theta$, $y=2\sin\theta - \sin2\theta$;
\item$(i_2)$
	$x=t^3-2t+t^{-1}$, $y=3t-t^{-1}$;
\item$(ii)$
        $x=2t(2t^2-1)$, $y=t^2(3t^2-1)$.
\eitem
\ecor

\begin{figure}[ht]

\centering \includegraphics[width=.4\linewidth]{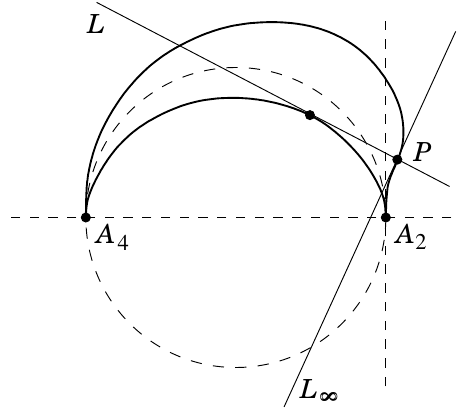}
\caption{The real quartic with $A_2$ and $A_4$}
		\label{fig:A2A4}
\end{figure}

\bpf 
Since $\deg C=4$, there is no room for more than one
non-generic tangency with $L_\infty$. Thus one of the cases
$(i)$--$(iv)$ of Lemma \ref{lemOneFlex}\ occurs.
We consider them separately.

\smallskip
$(i)$
Let $P$ be a smooth point of $C$.
Riemann-Hurwitz formula for the projection from $P$ implies that
there exists a unique line
$L_P$ through $P$ tangent to $C$ at another (smooth or singular) point.

Suppose that $(i_2)$ does not hold.
Then by Corollary \ref{corBranch}(b) 
all infinite points of $C$ are smooth.
Let $P$ be one of them. Let $Q$ be the point where $L_P$ is tangent to $C$.
Lemma \ref{lemTwoBranch} implies $L_P=L_\infty$, i.e., $Q\in L_\infty$.
Then, again by Lemma \ref{lemTwoBranch}, we have $L_Q=L_\infty$,
thus $(i_1)$ takes place.

In Case $(i_2)$, there are three cusps. However different choices of the
one tangent to $L_\infty$ lead to the same result because the cusps
are interchangeable by a projective automorphism of $\CP^2$
(one can see it, e.g., from the trigonometric parametrization).

\smallskip
$(ii)$ No other choice for $L_\infty$.

\smallskip
$(iii)$
Let $P$ be the flex point. Then $L_\infty$ is tangent to $C$ at $P$
by Lemma \ref{lemNoFlex}.
Applying Riemann-Hurwitz formula to the projection from $P$,
we see that there exists a line $L$ through $P$ which is tangent to $C$ at
a  smooth point. Contradiction with Lemma \ref{lemTwoBranch}.

{\it Remark.} 
The existence of such $L$ can be also derived from the
uniqueness of $C$ up to projective transformations. Indeed, we can
realize $C$ as a real curve in $\bR^2$ obtained by a small perturbation
of a double circle: $(x^2+y^2)^2 = \varepsilon y^3(x+1)$, $0<\varepsilon\ll1$.
Then $L$ is clearly visible in Figure \ref{fig:A2A4}.

\smallskip\pn
$(iv)$
Impossible by Lemma \ref{lemNoFlex}\ and Lemma \ref{corBranch}(a).
\epf 

\medskip
Thus there are only three candidates for solutions of the $\bC$-AlgDOP problem.
It remains to check that the linear equations for the metric (see Remark
\ref{rmq.compute.g}) have non-zero solutions, and then to select the real forms corresponding
to bounded domains.

\medskip
\bprop\label{CAlgDOP.quartic}
Up to affine transformations of $\bC^2$ and rescaling of $(a,b,c)$,
there are exactly three solutions to the $\bC$-AlgDOP problem
under condition that $\hat\Delta$ is irreducible
of degree $4$.
The curve $C=\{\hat\Delta=0\}$ is as in Corollary \ref{quartic}.
In Cases $(i_1)$ and $(ii)$ the formulas for $(a,b,c)$ are given in
Sections \ref{sect.swallow.tail} and \ref{sec.A2} respectively.
In Case $(i_2)$,
$a=9x^2+8xy$, $b=2y^2+3xy-8$, $c=y^2-12$.
\eprop

\medskip

\bprop\label{RAlgDOP.quartic}
Up to affine transformations of $\bR^2$ and rescaling of $(a,b,c)$,
there are exactly six solutions to the $\bR$-AlgDOP problem
under condition that $\hat\Delta$ is irreducible
of degree $4$: the three solutions
given in Proposition \ref{CAlgDOP.quartic} and those obtained from them
by the change of coordinates $(x,y)\mapsto(ix,y)$.
Only two among these solutions (those discussed in
Sections~\ref{sect.swallow.tail} and \ref{sec.A2}) correspond to bounded
domains in $\bR^2$ and thus provide a solution to the DOP problem.
\eprop

\bpf
First let us show that each projective curve $(i)$ and $(ii)$ of
Lemma \ref{quartic} has two real forms. It is easier to check this
fact for the dual curves. Indeed, for the nodal cubic (Case $(i)$) these are
$y^2=x^3\pm x$, and in Case $(ii)$ these are $y=x^4\pm x^2$.

The choice of the line at infinity is unique in all the six cases.
Indeed, for $(i_1)$ and $(ii)$ it is unique even over $\bC$.
In Case $(i_2)$ the curve $C$ has one or three real cusps, and if it
has three cusps, they are interchangeable by an automorphism of $\RP^2$.
Let us show that $\{\hat\Delta\ne0\}$ does not have bounded components
except the two cases.

\begin{figure}[ht]

\centering \includegraphics[width=.15\linewidth]{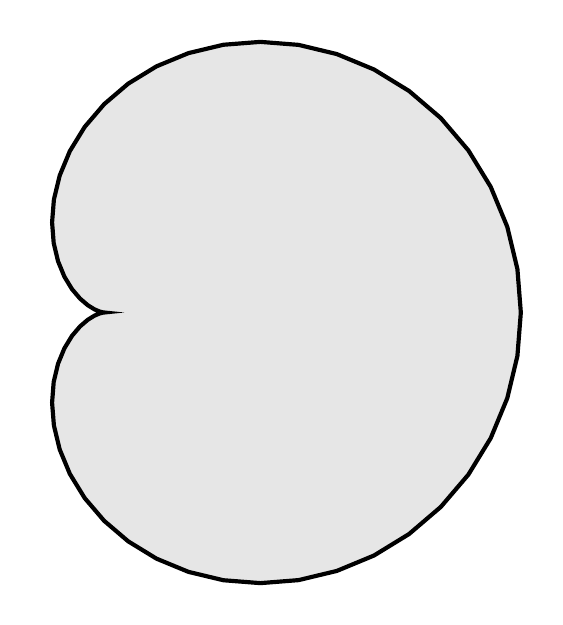}
\caption{$(1,1)$-hypercycloid: an irrelevant solution to AlgDOP problem}
		\label{fig:1-hypercyc}
\end{figure}

\smallskip
$(i)$ If $C$ has three real cusps, it looks as shown in Section
\ref{sec.A2}. So, it is clear that a tangent at a cusp is adjacent to
each component of the complement of $C$.
If $C$ has one real cusp, it can be realized as the $(1,1)$-hypercycloid
$x=2\cos\theta+\cos2\theta$, $y=2\sin\theta+\sin2\theta$ in some affine chart
(see Figure \ref{fig:1-hypercyc}).
So, in both cases $(i_1)$ and $(i_2)$, the line $L_\infty$
meets the closure of each component of the complement of $C$.

\smallskip
$(ii)$. If $\check C$ is $y=x^4+x^2$, it is convex in some affine chart,
hence so is $C$. Since $L_\infty$ is tangent to $C$, the affine part of $C$
is homeomorphic to a line dividing $\bR^2$ into two unbounded components.

\epf

%%%%%%%%%%%%%%%%%%%%%%%%%%%%%%%%%%%%%%%%%%%%%%%%%%%%%%%%%%%%%%%%%%%
%%%%%%%%%%%%%%%%%%%%%%%%%%%%%%%%%%%%%%%%%%%%%%%%%%%%%%%%%%%%%%%%%%%
%%%%%%%%%%%%%%%%%%%%%%%%%%%%%%%%%%%%%%%%%%%%%%%%%%%%%%%%%%%%%%%%%%%
%%%%%%%%%%%%%%%%%%%%%%%%%%%%%%%%%%%%%%%%%%%%%%%%%%%%%%%%%%%%%%%%%%%

\subsection{ Cubic factor of $\hat\Delta$ \label{subsec.cubic}}

In this section we suppose that
$\Delta=\Delta_3\Delta_1$ where
$\Delta_3$ is an irreducible cubic factor of  $\hat \Delta$
(As above, $(a,b,c,\hat\Delta)$ is a solution to
the $\bC$-AlgDOP problem  $\Delta=ac-b^2$).
By the genus formula~\eqref{eqGenusZ}, an irreducible cubic curve in $\CP^2$
is either smooth of genus one (and then depends of one parameter up to
projective transformations), or rational with a single singularity of
type $A_1$ (node) or $A_2$ (cusp).
In the latter case the curve is projectively rigid.

Let $C$ be the quartic curve defined by $\Delta=0$ and
let $C_3$ and $C_1$ be the respective irreducible components of $C$
(if $\deg\Delta=3$, then $C_1=L_\infty$).

\blem\label{lemCuOne} $C_3$ is rational.
\elem

\bpf Otherwise $C_3$ has nine flex points. They cannot
all be on $C_1\cup L_\infty$. So, this contradicts to
Lemma \ref{lemNoFlex}.
\epf 
\medskip

By an isomorphism of $\CP^2$, 
any rational cubic can by identified either with
the {\it nodal cubic} $y^2=x^3-x^2$ or with the 
{\it cuspidal cubic} $y^2=x^3$.
The nodal cubic has three flex points lying on the same line
and interchangeable by automorphisms of $\CP^2$.
The cuspidal cubic has a single flex point.

\blem\label{lemCuNodal}
Suppose that $C_3$ is a nodal cubic.
Then $\hat \Delta=\Delta_3$, the line
$L_\infty$ is tangent to $C_3$ at a flex point, and $C_1$
is the line passing through all the three flex points of $C_3$.
\elem
\bpf 
Let $L_0$ be the line passing through all the flex points of $C_3$.
Then $L_0\ne L_\infty$ by Lemma \ref{lemNoFlex}(a).
Thus, at least two flex points are not on $L_\infty$, hence
Lemma \ref{lemNoFlex}(b) implies that
a non-trivial component of $\Delta/\hat \Delta$
passes through them. Hence, $C_1=L_0$. and $\hat \Delta=\Delta_3$.

Suppose that $C_3$ has more than one point at infinity.
Then there is a point $P$ such that $(C_3.L_\infty)_P=1$.
Then $P$ is not a flex point by Lemma \ref{lemNoFlex}.
Hence, Riemann-Hurwitz formula for the projection from $P$
implies that there exists a line $L$ through $P$ which is tangent to
$C$ at some other point $Q$. If $Q$ were finite, then Lemma
\ref{lemTwoBranch}(b) would imply that $C_1$ passes through $P$.
This is impossible by Bezout's theorem because $C_1$ has already
three intersections with $C_3$ at the flex points. Thus, 
$Q\in L_\infty$. Applying the same arguments to $Q$, we obtain
a contradiction.

Thus, $C_3$ has a single point $P$ at the infinity.
It remains to show that $P$ is not the node of $C_3$.
Suppose it is. Choose coordinates $(X:Y:Z)$, $x=X/Z$,$y=Y/Z$,
so that $P=(1:0:0)$,
the axis $X=0$ is the tangent at a flex point,
and the tangents  at $P$ are $L_\infty$ and the axis $Y=0$.
Then,
up to rescaling of the coordinates,
$C_3$ admits a parametrization $x=\xi(t)=(t-1)^3/t$, $y=\eta(t)=t$.

So, we have an explicit parametric equation of a component of $\hat\Delta$.
As we pointed out in Remark \ref{rmq.compute.g}, then we have a system of
linear equations on the coefficients of $(a,b,c)$.
The rest of the proof is just checking by hand that this system does not have
any nonzero solution.

Applying Lemma \ref{lemBranch}(e,g) to the branches of $C_3$ at $P$,
we obtain $b=\NND111100$ and $c=\NND110100$.
Let $\gamma$ be the branch of $\hat\Delta$ at $P$ tangent to the axis $Y=0$.
We have $v_\gamma(x,y)=(-1,1)$, hence
by Lemma \ref{lemValu},
$v_\gamma(c)-v_\gamma(b)=2$. The values of $v_\gamma$ on the
monomials involved in $b$ and $c$ are $v_\gamma(1,xy,y,y^2)=(0,0,1,2)$.
Hence $c_{00}=c_{01}=0$, i.e., $c=c_{02}y^2$.
It follows that $c_{02}\ne0$ (otherwise $\Delta$ would be equal to
$b^2$), so we can assume that $c_{02}=1$.

Thus, the identity $b(\xi,\eta)\dot\eta=c(\xi,\eta)\dot\xi$ takes the form
$$
   b_{00} + b_{01}t + b_{02}t^2 + b_{11}(t^3-3t^2+3t-1)
   = t^2 (2t - 3 + t^{-2}).
$$
Equating the coefficients of $t^3,t^2,t,1$, we find 
$b_{11}=2$,
$b_{02}=3$,
$b_{01}=-6$,
$b_{00}=3$, i.e., $b=3(y-1)^2+2xy$. and hence
$b(\xi,\eta)=2t^3 - 3t^2 + 1$.
Substituting all these into $a(\xi,\eta)\dot\eta=b(\xi,\eta)\dot\xi$,
we obtain
$$
	a_{20}(t^4+\dots+t^{-2})=(2t^3-3t^2+1)(2t-3+t^{-2})=
		4t^4 + \dots + t^{-2}.
$$
A contradiction.
\epf 

\medskip

\blem\label{lemCuCuspidal}
Suppose that $C_3$ is a cuspidal cubic.
Then $L_\infty$ is tangent to $C_3$ at some point $P$.
Let $F$ be the flex point of $C_3$. Then:
\\
(a) If $P$ is the cusp, then $\hat\Delta=\Delta_3$ and $F\in C_1$.
\\
(b) If $P=F$, then $C_1$ is any line.
If, moreover, $\hat\Delta = \Delta$, then either $F\in C_1$ or
$C_1$ is tangent to $C_3$.
\\
(c) If $P$ is not as above, then $\hat\Delta=\Delta_3$ and $C_1$ is the
line $(PF)$.
\elem

\bpf
Let us prove that $L_\infty$ is tangent to $C_3$.
Suppose, it is not. Let us show that in this case 
$C_3\cap L_\infty\subset C_1$. Indeed,
let $Q\in C_3\cap L_\infty$. If $Q$ is the cusp of $C_3$, then
$Q\in C_1$ by Corollary \ref{corBranch}(b). If
$Q$ is a smooth point of $C_3$, then $Q\ne F$ by  Lemma \ref{lemNoFlex}(a)
and Riemann-Hurwitz formula for the projection from $Q$
implies that there is a line through $Q$ tangent to $C_3$,
hence $Q\in C_1$ by Lemma \ref{lemTwoBranch}(b).
Thus, we have shown that $C_3\cap L_\infty\subset C_1$.
Then, since  $C_3\cap L_\infty$
contains at least two points, we conclude that
 $C_1= L_\infty$.
However this is impossible because $F\not\in L_\infty$ by
Lemma \ref{lemNoFlex}(a) and then $F\in C_1$ by Lemma \ref{lemNoFlex}(b).
The obtained contradiction shows that $C_3$ is tangent to $L_\infty$.
So, let $P$ be the point where $C_3$ is tangent to $L_\infty$.

\smallskip
(a) Follows from Lemma \ref{lemNoFlex}(b).

\smallskip
(b) Suppose that $\hat\Delta = \Delta$ and $F\not\in C_1$.
Let $Q=C_1\cap L_\infty$. Let $L$ be a line through $Q$
tangent to $C_3$ at a finite point. Then $L=C_1$ by
Lemma \ref{lemTwoBranch}(b).

\smallskip
(c) By Lemma \ref{lemNoFlex}(b), we have $F\in C_1$ and $\hat\Delta\ne\Delta_3$.
Moreover, Riemann-Hurwitz formula for the projection from $P$
implies that there is a line through $P$ tangent to $C_3$,
hence $P\in C_1$ by Lemma \ref{lemTwoBranch}(b).
\epf 

\medskip
By combining Lemmas \ref{lemCuNodal} and \ref{lemCuCuspidal} and
computing $(a,b,c)$ from linear equations (see Remark \ref{rmq.compute.g}),
we summarize as follows.

\bprop\label{CAlgDOP.cubic}
Each solution of the $\bC$-AlgDOP problem where $\hat\Delta$ has an
irreducible cubic factor
is determined by $\hat\Delta$ and $\Delta:=ac-b^2$
up to rescaling of $(a,b,c)$ except the case $(ii_5)$.
Up to affine transformation of $\bC^2$, all realizable pairs
$(\Delta,\hat\Delta)$ are:

\smallskip
      $(i)$ (nodal cubic; see Section \ref{sec:DlePtCub})
      $\hat\Delta=x^3+x^2-y^2$, $\Delta=(3x+4)\hat\Delta$;

\smallskip
      $(ii_1)$ (see Section \ref{sec:CuCuSec})
      $\hat\Delta=\Delta=(x-1)(y^2-x^3)$;

\smallskip
      $(ii_2)$ (see Section \ref{sec:CuCuTan})
      $\hat\Delta=\Delta=(2y-3x+1)(y^2-x^3)$;

\smallskip
      $(ii_3)$ $\hat\Delta=\Delta=y(y^2-x^3)$;

\smallskip
      $(ii_4)$ $\hat\Delta=\Delta=y^2-x^3$;

\smallskip
      $(ii_5)$ $\hat\Delta=y^2-x^3$,
      $\Delta=(\alpha x+\beta y+\gamma)\hat\Delta$,
      $(\alpha,\beta,\gamma)\ne(0,0,0)$;

\smallskip
      $(ii_6)$ (Lemma \ref{lemCuCuspidal}(c))
      $\hat\Delta=\Delta=(1 + 2x - 2y)(x^3 - y^2 - 3xy^2 + 2y^3)$;

\smallskip
      $(ii_7)$ $\hat\Delta=y-x^3$, $\Delta=(\alpha x+\beta y)\hat\Delta$,
      $(\alpha,\beta)\ne(0,0)$.

\smallskip
\noindent
The solution $(i)$ has two real forms. Only one of them provides a bounded
solution to the $DOP$ problem. Each of the other solutions has one
real form. Only $(ii_1)$ and $(ii_2)$ provide bounded
solutions to the $DOP$ problem.
\eprop

\brmq The cusp at infinity (Case $(ii_7)$ leads to a non compact domain. Moreover,  because of the form of the measure, it is not possible even
in the non compact case.\ermq

%%%%%%%%%%%%%%%%%%%%%%%%%%%%%%%%%%%%%%%%%%%%%%%%
%%%%%%%%%%%%%%%%%%%%%%%%%%%%%%%%%%%%%%%%%%%%%%%%
%%%%%%%%%%%%%%%%%%%%%%%%%%%%%%%%%%%%%%%%%%%%%%%%

\subsection{ Quadratic factor of $\hat\Delta$ \label{subsec.quadratic}}

In this section we suppose that
$\Delta=\Delta_2\tilde\Delta_2$ where
$\Delta_2$ is an irreducible quadratic factor of  $\hat \Delta$.
As above, $(a,b,c,\hat\Delta)$ is a solution to
the $\bC$-AlgDOP problem  $\Delta=\det g = ac-b^2$.
Let $C$, $C_2$ and $\tilde C_2$ be the corresponding curves in $\CP^2$,
Up to an affine linear transformation of $\bC^2$
there are two cases for $C_2$:
a hyperbola $xy-1$ and a parabola $y-x^2$.

\bprop\label{AlgDOP.hyperb}
Let $\Delta=ac-b^2=\Delta_2\tilde\Delta_2$ with $\Delta_2=xy-1$
and let $\Delta_2$ be a factor of $\hat\Delta$.
Then $(a,b,c,\hat\Delta)$
is a solution of the $\bC$-AlgDOP problem if and only if
% $\hat\Delta=\Delta_2$ (up to a constant factor) and
\[
     \bpm a&b\\b&c\epm
    =\Delta_2\bpm \alpha&\beta\\\beta&\gamma \epm
     +r\bpm x^2&-xy\\-xy&y^2\epm,
\]
with $(r,\alpha\gamma-\beta^2)\ne(0,0)$ and $(\alpha,\beta,\gamma)\ne(0,0,0)$,
and one of the following cases occurs
up to rescaling and exchange of $x$ and $y$:
\benum
 \item $\hat\Delta=\Delta_2$.
 \item $\alpha=\beta=0$ and $\hat\Delta=x\Delta_2$; in this case
       $\Delta=\gamma rx^2\Delta_2$.
\eenum
Furthermore, $\deg\Delta=2$ if and only if $\alpha=\gamma=0$ and $\beta=2r$
{\rm(then $(a,b,c)=r(x^2,xy-2,y^2)$)}.
Otherwise $\deg\Delta=4$.
\eprop

\bpf
By solving the linear equations~\eqref{eqbordparam} for the coefficients of
$a$, $b$, and $c$, we find the announced form of $g$. Then we have
$$
   \tilde\Delta_2 = (\alpha\gamma-\beta^2)\Delta_2
             + r(\gamma x^2 + 2\beta xy +\alpha y^2)
$$
whence the non-vanishing conditions
$(r,\alpha\gamma-\beta^2)\ne(0,0)$ and $(\alpha,\beta,\gamma)\ne(0,0,0)$.
We also see that $\deg\Delta=4$ unless $\alpha=\gamma=0$ and $\beta=2r$ 
in which case $\deg\Delta=2$.

If $r=0$, then $\Delta=\Delta_2^2$, hence $\hat\Delta=\Delta_2$
(recall that $\hat\Delta$ does not have multiple factors).
So, from now on we assume that $r\ne 0$.
%If $r\ne 0$, we may assume that $r=1$. 
%Then up to exchange of $\alpha$ and $\gamma$ we have the following cases.

\smallskip
Case 1. $\alpha\gamma-\beta^2\ne 0$ and $(\gamma,\alpha)\ne(0,0)$. Then 
Then $\tilde C_2$ is a nonsingular conic such that
$C_2\cap\tilde C_2\cap L_\infty=\varnothing$,
thus $\hat\Delta=\Delta_2$ by Lemma~\ref{lemTwoBranch}
applied to the tangent to $C_2$ passing through a point from
$\tilde C_2\cap L_\infty$.

\smallskip
Case 2. $\alpha=\gamma=0$ (and hence $\beta\ne 0$).
Then we have $\hat\Delta_2=\beta(\beta-(\beta-2r)xy)$.
If $\beta=2r$, we have $\Delta_2=\text{const}$, hence $\hat\Delta=\Delta_2$.
Otherwise, plugging the parametrization
$x=\xi(t)=t$,
$y=\eta(t)=\beta/((\beta-2r)t)$
of $\tilde C_2$ into $b\dot\eta-c\dot\xi$,
we obtain $-2\beta r/(\beta-2r)$, i.e., the relation~\eqref{eqbordparam}
is not satisfied for $\tilde C_2$.
Thus $\tilde\Delta_2$ is not a factor of $\hat\Delta$.

\smallskip
Case 3. $\alpha\gamma-\beta^2=0$. Up to exchange of $x$ and $y$,
we may set $\gamma=1$, $\alpha=\beta^2$. Then
$\tilde\Delta_2=r(x+\beta y)^2$.
Plugging $\xi(t)=\beta t$, $\eta(t)=-t$
into \eqref{eqbordparam}, we obtain
$a\dot\eta-b\dot\xi=\beta(b\dot\eta-c\dot\xi)
=-2\beta^2(1+(\beta-r)t)^2$. Thus $x+\beta y$ can be a factor of $\hat\Delta$
if and only if $\beta=0$.
\epf
\medskip

Up to affine change of coordinates there are three real forms of
$xy-1$: a circle (the only one which gives a bounded solution to the DOP
problem; see Section \ref{circle}), a hyperbola, and a purely imaginary
conic $1+x^2+y^2$. The latter real solution will be used in the study of
the SDOP problem in the case $\Omega=\bR^2$ (Section \ref{no.boundary}).

\bprop\label{AlgDOP.circle}
Let $\Delta=ac-b^2=\Delta_2\tilde\Delta_2$ with $\Delta_2=1-x^2-y^2$
and let $\Delta_2$ be a factor of $\hat\Delta$.
Then $(a,b,c,\hat\Delta)$
is a solution of the $\bR$-AlgDOP problem if and only if
$\hat\Delta=\Delta_2$ (up to a constant factor) and, up to rotation,
\[
     \bpm a&b\\b&c\epm
    =\Delta_2\bpm \alpha&0\\0&\gamma \epm+r\bpm1-x^2&-xy\\-xy&1-y^2\epm.
\]
where at most one of the numbers $r$, $\alpha+r$, $\gamma+r$ is zero.
Furthermore, $\deg\Delta=2$ if and only if $\alpha=\gamma=0$ and $r\ne0$.
Otherwise $\deg\Delta=4$.
\eprop

\medskip
\bpf
By solving the linear equations~\eqref{eqbordparam} for the coefficients of
$a$, $b$, and $c$, we find the announced form of $g$ but with a matrix
$\left(\begin{smallmatrix}\alpha&\beta\\\beta&\gamma\end{smallmatrix}\right)$.
However by rotation we may reduce to the case $\beta=0$.
The rest can be easily derived from Proposition~\ref{AlgDOP.hyperb} after a
suitable change of variables.
\epf

\medskip

\bprop\label{AlgDOP.parab}
Let $\hat\Delta=\Delta_2\tilde\Delta_2$ where $\Delta_2=y-x^2$.
Then a solution of $\bC$-AlgDOP
problem with this $\hat\Delta$ exists if and only if one of the following cases
occurs up to affine linear change of coordinates:
\benum
\item (coaxial parabolas; see Section~\ref{parab1})
   $\tilde\Delta_2=y-\alpha x^2-1$, $\alpha\ne 1$;
\item (parabola with a tangent and the axis; see Section~\ref{parab2})
   $\tilde\Delta_2=y(x-1)$;
\item (parabola with two tangents; see Section~ \ref{sec.B2})
   $\tilde\Delta_2=(y+1)^2-4x^2$;
\item
   $\tilde\Delta_2$ is $1$, $x$, or $y$.
\eenum
\eprop

The real forms are evident (in Case (3) the tangents may be either real or
complex conjugate). The only bounded solutions to the DOP problem
are the first three cases 

\medskip
\bpf
If $\deg\tilde\Delta_2<2$, everything is realizable.
Indeed, by a change of coordinates preserving the parabola, any
line can be transformed to one of $x=0$, $y=0$, or $y=1$ 
(which is (3) with $\alpha=0$), and the
system of equations can be easily solved in all the three cases.

Let then $\deg\tilde\Delta_2=2$.
By Lemma \ref{lemTwoBranch} we know that if
$P\in\tilde C_2\cap L_\infty$ but $P\ne\hbox{(0\,:\,1\,:\,0)}$, then $\tilde C_2$
contains the line passing through $P$ and tangent to $C_2$.
Thus, if $\tilde\Delta_2$ is not as required, it can be transformed
to one of:
$y(y-1)$, $x(x-1)$, $y-x^2-1$, $y-2x^2$. One easily checks that
there are no solutions in these cases.

\epf

The following proposition is not needed for the classification of
compact solutions but it will be helpful in the study of the
non-compact case in Section~\ref{non.compact.with.boundary}.

\bprop\label{AlgDOP.parab.2}
Let $(a,b,c,\hat\Delta)$ be
a solution of the $\bR$-AlgDOP problem.
Suppose that $y-x^2$ is a factor of $\hat\Delta$.
Then, for some $\alpha,\beta,\gamma,r,\mu,\nu\in\bR$,
\beq\label{eq.parab.2}
  g
   =(y-x^2)\bpm \alpha&\beta\\\beta&\gamma \epm
     +r\bpm x&2y\\2y&4xy\epm
     +(\lambda+\mu y)\bpm 1&2x\\2x&4y\epm
\eeq

Up to change of variables, we may suppose that
either $r=0$, or $\lambda=\mu=0$.
Moreover, up to scaling and change of variables, we have:

\benum
\item[\rm(1)] $\Delta=C(y-x^2)^2$ if and only if   one of the following cases occurs
  \benum
  \item[\rm(1i)] $r=\lambda=\mu=0$, $\alpha\gamma-\beta^2\ne 0$;
  \item[\rm(1ii)] $\alpha=\gamma=\lambda=\mu=0$, $r=-\beta=1$.
  \item[\rm(1iii)] $r=\beta=0$, $-\alpha=\mu=\pm1$, $-\gamma=4\lambda=4$.
  \eenum
In these cases $g$ is, respectively,
$$
    (y-x^2)\bpm\alpha&\beta\\\beta&\gamma\epm, \quad
    \bpm x&y+x^2\\y+x^2&4xy\epm,\quad
    \bpm 1 & 2x\\ 2x & 4x^2\epm\pm\bpm x^2 & 2xy\\ 2xy& 4y^2\epm.
$$

\item[\rm(2)] $\deg\Delta=3$ if and only if one of the following cases occurs
   \benum
   \item[\rm(2i)] $\beta=\gamma=r=0$, $\mu=-\alpha=1$, $\lambda=\pm1$;
   \item[\rm(2ii)] $\alpha=\beta=\lambda=\mu=0$, $r=1$, $\gamma\in\{0,1\}$;
   \item[\rm(2iii)] $\beta=\gamma=\mu=r=0$, $\lambda=1$, $\alpha=\pm1$.
   \eenum
\item[\rm(3)] $\Delta=y-x^2$
     if and only if $\alpha=\beta=r=\mu=0$, $\lambda=1$, $\gamma\ne-4$;
\eenum
\eprop

\bpf
By solving the system of linear equations \eqref{eqbordparam}, we find
\eqref{eq.parab.2}.
The change of variables $x'=x+q$, $y'=y+2qx+q^2$ transforms the parameters
$(r,\lambda,\beta)$ into
\beq
   \label{eq.parab.3}
      r' = r-2q\mu,\qquad
      \lambda' = \lambda-qr+q^2\mu,\qquad
      \beta'=\beta+(2\alpha+4\mu)q.
\eeq
 Thus
if $\mu\ne 0$, we may assume that $r=0$, and if $\mu=0$ but $r\ne0$, we may assume
that $\lambda=0$. The rest of the proof is a straightforward case by case computation
(in case (2iii), when $r=\mu=0$ and $\alpha\ne0$, we kill $\beta$ using
\eqref{eq.parab.3}).
At the final stage we use the variable change $(x,y)\mapsto(px,p^2y)$ to
normalize the coefficients.
\epf

%%%%%%%%%%%%%%%%%%%%%%%%%%%%%%%%%%%%%%%%%%%%%%%%
%%%%%%%%%%%%%%%%%%%%%%%%%%%%%%%%%%%%%%%%%%%%%%%%
%%%%%%%%%%%%%%%%%%%%%%%%%%%%%%%%%%%%%%%%%%%%%%%%

\subsection{ All factors of $\hat\Delta$ are linear\label{subsec.linear}}

\blem\label{lem.ParLines}
Let $(a,b,c,\hat\Delta)$. $\Delta=ac-b^2$, be a solution of the $\bC$-AlgDOP problem
such that $\hat\Delta=(1-x^2)\tilde\Delta_2$. Then, up to an
affine linear transformation, we have
$$
      g=\bpm a&b\\b&c\epm = \bpm \alpha(1-x^2) & \beta(1-x^2)\\
                                 \beta (1-x^2) & c(x,y) \epm
$$
where either $\alpha=0$ and then $\Delta=\beta^2(1-x^2)^2$, or
$\beta=0$ and then $\Delta=\alpha(1-x^2)c(x,y)$.

Moreover, $\hat\Delta$ cannot have a factor $x-x_0$ with $x_0\ne\pm1$.
\elem

\bpf
By solving the linear equations~\eqref{eqbordparam} for the coefficients of
$a$, $b$, and $c$, we find the announced form of $g$ but with arbitrary
$\alpha$ and $\beta$. However, if $\alpha\ne 0$, then the change of
variables $x'=x$, $y'=\alpha y-\beta x$ kills $\beta$.

Suppose that $x-x_0$ is a factor of $\hat\Delta$. If $\alpha=0$, then
$x_0=\pm1$ because $\Delta=(1-x^2)^2$.
If $\beta=0$, then for $(\xi(t),\eta(t))=(x_0,t)$ we have
$a(\xi,\eta)\dot\eta-b(\xi,\eta)\dot\xi=1-x_0^2$ whence $x_0=\pm1$
by~\eqref{eqbordparam}.
\epf

\medskip
The following lemma is very similar to Proposition \ref{AlgDOP.circle}.
We shall see in Section \ref{circle} that there is a deep reason for this.
\medskip

\blem\label{lem.3lines}
Let $(a,b,c,\hat\Delta)$, $\Delta=ac-b^2$, be a solution of the $\bC$-AlgDOP problem
such that $\Delta_3:=xy(1-x-y)$ divides $\hat\Delta$.
Then $\hat\Delta=\Delta_3$ up to constant factor, and
\[
     g=\bpm a&b\\b&c\epm
    =(1-x-y)\bpm \alpha x&0\\0&\gamma y\epm+r\bpm x(1-x)&-xy\\-xy&y(1-y)\epm.
\]
where at most one of the numbers $r$, $\alpha+r$, $\gamma+r$ is zero.
Furthermore,
\benum
\item
 $\Delta= \Delta_3$ if and only if $\alpha=\gamma=0$ and $r=1$;
\item
 $\Delta=x\Delta_3$ if and only if $\gamma+r=0$ and $r(\alpha+r)=1$;
\item
 $\Delta=y\Delta_3$ if and only if $\alpha+r=0$ and $r(\gamma+r)=1$;
\item
 $\Delta=(1-x-y)\Delta_3$ if and only if $r=0$ and $\alpha\gamma=1$;
\eenum
\elem

\bpf
By solving~\eqref{eqbordparam} for parametrizations of the three lines
$x=0$, $y=0$, and $x+y=1$, we find the required form of $g$.
Hence $\Delta=\Delta_3\Delta_1$ where
$$
   \Delta_1 = (\alpha+r)(\gamma+r) - \gamma(\alpha+r)x - \alpha(\gamma+r)y
$$
and we easily derive the assertions (1)--(4) as well as
the fact that at most one of $r$, $\alpha+r$, $\gamma+r$ is zero.

Let us show that $\hat\Delta=\Delta_3$ up to constant factor.
If $\alpha=\gamma=0$ or one of $r$, $\alpha+r$, $\gamma+r$ is zero, this follows from
the assertions (1)--(4) (recall that $\hat\Delta$ does not have multiple factors).
So, we assume that
$r(\alpha+r)(\gamma+r)\ne 0$ and $(\alpha,\gamma)\ne(0,0)$.
Let us show that in this case a parametrization $(x,y)=(\xi(t),\eta(t))$ of
the line $\Delta_1=0$ does not satisfy condition \eqref{eqbordparam}.

Indeed, if $\alpha\ne0$ and $\gamma=0$ (the case $\alpha=0$ and $\gamma\ne0$ is
similar), then
$\xi=t$, $\eta=(\alpha+r)/\alpha$, hence
$a(\xi,\eta)\dot\eta-b(\xi,\eta)\dot\xi = \frac r\alpha(\alpha+r)t$
and the coefficient of $t$ in non-zero.
If $\alpha\ne0$ and $\gamma\ne0$, then
$\xi=(1+\frac r\gamma)t$, $\eta=(1+\frac r\alpha)(1-t)$ and
$$
    a(\xi,\eta)\dot\eta-b(\xi,\eta)\dot\xi
    = \alpha C\big((\gamma+r)t+(\alpha-\gamma)t^2\big)
$$
with $C=r(\alpha+r)(\gamma+r)/(\alpha\gamma^2)$ and again
the coefficient of $t$ is non-zero.
\epf

\medskip
\bprop\label{prop.only.lines}
Let $\hat\Delta$ is a product of linear factors.
Then a solution of $\bC$-AlgDOP
problem with this $\hat\Delta$ exists if and only if one of the following cases
occurs up to affine linear change of coordinates:
\benum
\item (square; see Section~\ref{square}) $\hat\Delta=(x^2-1)(y^2-1)$;
\item (triangle; see Section~\ref{triangle}) $\hat\Delta=xy(1-x-y)$;
\item $\hat\Delta$ is one of $xy(x+y)$, $y(x^2-1)$, $xy$, $x^2-1$, $x$, or $1$.
\eenum
\eprop

\bpf
Follows from Lemmas~\ref{lem.ParLines}--\ref{lem.3lines} (to exclude
four concurrent lines, one should slightly modify Lemma~\ref{lem.3lines}).
\epf

\medskip
By combining Propositions~\ref{RAlgDOP.quartic}, \ref{CAlgDOP.cubic},
\ref{AlgDOP.circle}, \ref{AlgDOP.parab}, and \ref{prop.only.lines}, we get a proof
of Theorem~\ref{thm.central}.

%%%%%%%%%%%%%%%%%%%%%%%%%%%%%%%%%%%%%%%%%%%%%%%%%%%%%%%%%%%%%%%%%%%%%
%%%%%%%%%%%%%%%%%%%%%%%%%%%%%%%%%%%%%%%%%%%%%%%%%%%%%%%%%%%%%%%%%%%%%
%%%%%%%%%%%%%%%%%%%%%%%%%%%%%%%%%%%%%%%%%%%%%%%%%%%%%%%%%%%%%%%%%%%%%
%%%%%%%%%%%%%%%%%%%%%%%%%%%%%%%%%%%%%%%%%%%%%%%%%%%%%%%%%%%%%%%%%%%%%
%%%%%%%%%%%%%%%%%%%%%%%%%%%%%%%%%%%%%%%%%%%%%%%%%%%%%%%%%%%%%%%%%%%%%
%%%%%%%%%%%%%%%%%%%%%%%%%%%%%%%%%%%%%%%%%%%%%%%%%%%%%%%%%%%%%%%%%%%%%
%%%%%%%%%%%%%%%%%%%%%%%%%%%%%%%%%%%%%%%%%%%%%%%%%%%%%%%%%%%%%%%%%%%%%
%%%%%%%%%%%%%%%%%%%%%%%%%%%%%%%%%%%%%%%%%%%%%%%%%%%%%%%%%%%%%%%%%%%%%
%%%%%%%%%%%%%%%%%%%%%%%%%%%%%%%%%%%%%%%%%%%%%%%%%%%%%%%%%%%%%%%%%%%%%
%%%%%%%%%%%%%%%%%%%%%%%%%%%%%%%%%%%%%%%%%%%%%%%%%%%%%%%%%%%%%%%%%%%%%

\section{The bounded $2$-dimensional models\label{comments.models}}

\subsection{Generalities\label{subsec.generalities}}

In this section, we will explore separately the various 2  dimensional compact  models. It turns out that for some values (in general half-integer) of the parameters appearing in the measure, one may produce a geometric interpretation, coming in general from Lie groups or symmetric spaces, as it is the case for the one dimensional Jacobi operator (Section~\ref{dim1}). We do not pretend to present all the possible origins of the various models, but  we provide some insight whenever they are at hand and relatively easy to produce. Moreover, these geometric interpretation may lead to natural higher dimensional models for the DOP problem.

Recall  that the boundary of $\Omega$ is an algebraic curve of degree at most $4$. When the degree is $4$, this boundary is $\{\Delta =0\}$  where  $\Delta$ is the determinant of the matrix $(g^{ij})$. Among the admissible measures, one may chose $\rho(x)= \Delta^{-1/2}$, which corresponds to the Laplace-Beltrami operator associated with the (co-)metric $g$. It turns out that in every such example, this Laplace-Beltrami operator has constant curvature, either $0$ or positive. We did not succeed in proving this fact in the general setting (and we do not even know if is is true in higher dimension;
see Remark \ref{rmq.soukhanov} where we also mention a result of
Soukhanov~\cite{Soukhanov2014} in this direction).
However, when the boundary has degree less than 4, it is not always true that the curvature is constant (see Section~\ref{sec:DlePtCub}). But even in this latter case, when the measure has density $\Delta_1^{-1/2}$, where $\Delta_1$ is the irreducible equation of the boundary (while in this model $\Delta$ has degree 4 and $\Delta_1$ degree 3), there exists a natural interpretation coming from a 4-dimensional sphere. 

Then, one may interpret the associated model as some quotient of the Euclidean or spherical Laplace operator through some  discrete or continuous  symmetry subgroups.  When the curvature is $0$ (Section~\ref{sec.B2} and Section~\ref{sec.A2}), this shows some relation with root systems and the associated Hall  polynomials \cite{Macdo}, with connection to Hecke algebras. 
See also Araki \cite{Araki} and Harish-Chandra \cite{Harish1,Harish2}. Many other natural geometric interpretations come from spherical functions on rank 2 symmetric spaces (see  Helgason \cite{Helgason}, Heckman et al. \cite{HeckOp,Heck, opdam, HeckSchil}). For references on Dunkl operators, we also refer to Dunkl \cite{Dunkl} or
to a more recent paper by
R\"osler \cite{Rosler}.
\par

\par
From a general point of view, there is a dictionary linking the angles of the reflection associated with the symmetries and the type of singularities of the boundary of $\Omega$: double points, cusps and tangency points
correspond respectively to angles $\pi/2$, $\pi/3$, $\pi/4$.

 It turns out that many of the models described  above have some nice geometric interpretation in terms of compact homogeneous spaces $M= G/H$: we try to interpret the given operator as the Laplace-Beltrami operator $\De_G$ on $G$ acting on some specific functions $(X,Y): G\mapsto \bR^2$.  

In this whole section, the identification will be made with the Laplace operator acting on the $n$-dimensional sphere, on the $n$-dimensional Euclidean space or on some classical Lie group such as $SO(n)$ and $SU(n)$. For the sake  of clarity, we recall here some well known formulas and facts on these operators. The general principle is the following. When $\LL$ is a Laplace-Beltrami operator (or more generally any second order differential operator with no 0-order term) on some model space $E$, recall that the associated carr\'e du champ is defined by
$$
    \Ga_{\LL}(f,g)= \frac{1}{2} \Big(\LL(fg)-f\LL(g)-g\LL(f)\Big)\, ,
$$
and $\LL$ satisfies the change of variable formula~\eqref{chgt.de.variables}. Then, we are looking for pairs $(X^1,X^2)$ of real functions $E\to \bR$
such that $\LL(X^i)= L^i(X^1, X^2)$, and $\Ga(X^i,X^j)= G^{ij}(X^1,X^2)$,
where $L^i$ are some
degree 1 polynomials and $G^{ij}$ are degree 2 polynomials in the two variables
$(X^1,X^2)$. Then, from the change of variable formula~\eqref{chgt.de.variables},
for any smooth function $\Phi: \bR^2\to \bR$, one has
$\LL\big(\Phi(X^1,X^2)\big) = \LL_1(\Phi)(X^1,X^2)$, where 
$$
   \LL_1(f)= \sum_{ij} G^{ij}(x)\partial_{ij} f+ \sum_i L^i(x)\partial_i f.
$$
We shall say that such an operator $\LL_1$ is the image measure of $\LL$ through $(X^1,X^2)$. Moreover, it is immediate  that, if $\mu$ is the reversible measure for
$\LL$, then $\LL_1$ has reversible measure the image of $\mu$ through $(X^1,X^2)$.

Indeed, what is immediate from the study of the various models is the knowledge of
$\Ga(X^i,X^j)= g^{ij}$  and  the density measure $\rho$.
From~\eqref{sym.diffusion}, it is then immediate that
\beq
      \label{eq.drift}
      L^i(x) = \sum_j \partial_j g^{ij} + g^{ij}\partial_j \log \rho.
\eeq

Through an affine change of coordinates, one is reduced to find two eigenvectors $X^1$ and $X^2$ of $\LL$  for which $\Ga(X^i,X^j)$ satisfy a quadratic relation $\Ga(X^i,X^j)= G^{ij}(X^i,X^j)$. In this respect, similar problems are studied (although mainly in dimension 1) in the study of isoparametric surfaces (see Cartan \cite{cartan1, cartan2,cartan3,cartan4}).

 It can be quite hard to find from which model space a given model comes from. Spectral analysis can be useful: indeed,  for any polynomial $P(x,y)$ and 
 whenever $P(X^1,X^2)\in \cL^2(\mu)$, the spectrum of $\LL_1$ is embedded in the discrete spectrum of $\LL$. But, as it happens in Section~\ref{sec.B2} and Section~\ref{sec.A2}, it could be that the reversible measure $\mu$  for $\LL$ has infinite mass, and that $X^1$ and $X^2$ are eigenvectors for $\LL$ which are not in $\cL^2(\mu)$. However, whenever $\LL_1$ is the image of some geometric operator $\LL$ on some compact  model space $E$, the spectrum of $\LL_1$ is imbedded in the spectrum of $\LL$. Nevertheless, this could be misleading in some specific situation. For example, 
on the unit sphere $\bS^n$ imbedded in $\bR^{n+1}$ with the induced Riemannian metric,
the spectrum of the associated Laplace-Beltrami operator $\De_{\bS^n}$ is $\{-k(k+n-1), k\in \bN\}$. But then, for any integer $p$, the spectrum of $p^2\De_{\bS^n}$ is included into $\{-k(k+p(n-1), k\in \bN\}$, which is the spectrum of a sphere of dimension $p(n-1)+1$, and therefore, since in general we know $\LL_1$ only up to some scaling factor, we are not even able to determine the dimension of the sphere it may come from (if ever). We already saw this phenomenon in the case of Jacobi operators with parameter $(p,p)$ for which  we have two distinct geometric interpretation, one coming from $\bS^{p}$ and another one coming from $\bS^{2p-1}$ (Section~\ref{dim1}).

As mentioned above, for the purpose of the description of our 11 models, we shall mainly use a few model spaces, namely Euclidean, spheres, $SO(n)$ and $SU(n)$. In order to be able to carry the identification described above from these models, it is worth to describe the Laplace-Beltrami (or Casimir) operators for those models, in a simple way leading to the further interpretations. In some cases, it is also useful to extend the operators $\LL$ and $\Ga$ to complex valued functions, and we shall do that without further notice.

\par 
If $E_n$ is an $n$-dimensional Euclidean space, and $x^i$ the coordinates in some orthonormal basis, one has
for the Euclidean Laplace operator $\De_{E_n}$:
$$
      \De_{E_n}(x^i)= 0,\;\; ~\Ga_{E_n}(x^i, x^j)= \delta^{ij},
$$
and the dimension does not appear in these relations, hence we can omit the subscript $n$. These descriptions of course do no depend of the chosen orthonormal basis in $E_n$.

On the unit sphere $\bS^n$ imbedded  into $\bR^{n+1}$,  and for the restriction to $\bS^n$ of the same Euclidean coordinates $x^i$, one has
for the Laplace operator $\De_{\bS^n}$ 
\beq
     \label{eq.lapl.spheres}
     \De_{\bS^n}(x^i)= -nx^i,\qquad \Ga_{\bS^n}(x^i,x^j)= \delta^{ij}-x^ix^j.
\eeq
For complex coordinates $z^j=x^j+ix^{j+n}$, $j=1,\dots,n$, on $\bR^{2n}$
(here $i=\sqrt{-1}$),
one easily derives from \eqref{eq.lapl.spheres} that
\beq
     \label{eq.lapl.spheres.compl}
      \Ga_{\bS^{2n-1}}(z^j,z^k)=-z^j z^k,\qquad
      \Ga_{\bS^{2n-1}}(z^j,\bar z^k)=2\delta^{jk}-z^j\bar z^k.
\eeq

As previously, we can write $\Ga_\bS$ since it does not depend on the dimension.
When $F$ is the restriction to the unit sphere in $\bR^{n+1}$ of some smooth
function in the Euclidean space,  $\De_{\bS^n} (F)$  and $\Ga_{\bS}(F)$ may be computed
from the related quantities $\De_{\bE}$ and $\Ga_{\bE}$ in the ambient Euclidean space,
since they are the restriction to the sphere  of the quantities
\beq
    \label{lapl.sph.2}
\begin{aligned}
    \De_{\bS^n}(F)&=\De_{\bE}(F) -(r\partial_r)^2 F -(n-1) r\partial_r F, \\
    \Ga_{\bS}(F,G) &= \Ga_{\bE}(F,G)- (r\partial_r F)(r\partial_r G)
\end{aligned}
\eeq
where $r\partial_r F= \sum_i x^i \partial_i F$. Hence, if $F$ and $G$ are
homogeneous of degree $a$ and $b$ respectively, we have
\beq
    \label{lapl.sph.homog}
    \De_{\bS^n}(F)=\De_{\bE}(F) - a(a+n-1) F, \quad
    \Ga_{\bS}(F,G) = \Ga_{\bE}(F,G)- abFG.
\eeq

As mentioned above the spectrum of $-\De_{\bS^n}$ is $\{k(k+n-1), k\in \bN\}$.
The eigenspace associated with $-k(k+n-1)$ consists of the restriction to the sphere of degree $k$ harmonic homogeneous  polynomials (see Stein and Weiss \cite{SteinWeiss}).

Beyond the case of spheres, we shall also use Casimir operators on the semi-simple groups $SU(n)$ and $SO(n)$.  Once again, in order to describe them, we consider the entries $x^{ij}$ ($SO(n)$ case) and $z^{ij}$ ($SU(n)$ case) as functions on the group (complex valued in the latter case) and describe the operators through the action of $\LL$ and $\Ga$ on them.

For $SO(n)$, up to some constant, we have,
\beq
    \label{lapl.so(n)}
    \De_{SO(n)}(x^{ij})=-(n-1) x^{ij},\qquad
    ~\Ga_{SO(n)}( x^{kl},x^{pq})= \delta^{kp}\delta^{lq}- x^{kq}x^{pl}.
\eeq
For $SU(n)$ the formulas are similar:
\beq
    \label{lapl.su(n)1}
    \De_{SU(n)}(z^{ij})=-2\frac{(n-1)(n+1)}{n}  z^{ij},
\eeq
and also 
\beq
    \label{lapl.su(n)2}
\begin{aligned}
    &\Ga_{SU(n)}( z^{kl},z^{pq})=-2z^{kq}z^{pl}+\tfrac{2}{n}z^{kl}z^{pq},\\
    &\Ga_{SU(n)}( z^{kl},\bar z^{pq})
    = \;\;2\delta^{kp}\delta^{lq} -\tfrac{2}{n}z^{kl}\bar z^{pq}.
\end{aligned}
\eeq
 
In order not to get confused in the notation, we shall use upper case letters
$(X,Y)$ instead of $(X^1,X^2)$ for the coordinate system in the different
$2$-dimensional models $\Omega$, and  lower case letters $(x_i)$ for the
coordinates on the geometric model it comes from
(we switch also to lower indices because we will not use
much summation over repeated indices whereas polynomial expressions we will be
widely used).

\medskip
To conclude this subsection, we give an example of computation of
the coefficients of $\De_{S^n}$, $n=p+q-1$, pushed down
to $\bR$ through the function $X=x_1^2+\dots+x_p^2$.
We shall obtain the Jacobi operator (Section~\ref{dim1}).

By \eqref{defGamma} we have
$$
   \De_{S^n}(x_i^2) = 2x_i\De_{S^n}(x_i) + 2\Ga(x_i,x_i),
$$
hence $\De_{S^n}(x_i^2) = 2-2(p+q)x_i^2$ (by \eqref{eq.lapl.spheres})
whence $\De_{S^n}(X)=2p-2(p+q)X$ by linearity.
Further, $\Ga$ is a derivation with respect to each argument, i.e.,
$\Ga(fg,h)=\Ga(h,fg)=f\Ga(g,h)+g\Ga(f,h)$ whence
$$
  \Ga_{\bS}(x_i^2,x_j^2)=4x_ix_j\Ga_{\bS}(x_i,x_j)=4x_ix_j(\delta_{ij}-x_ix_j)
$$
(see \eqref{eq.lapl.spheres}). Hence
\beq
   \label{expl.jacobi}
   \Ga_{\bS}(X,X)=\sum_{i,j}\Ga_{\bS}(x_i^2,x_j^2)
                 =4\sum_i x_i^2 - 4\sum_{i,j}x_i^2x_j^2
                 = 4X-4X^2
\eeq
This means that the obtained one-dimensional operator $\LL$ reads
$$
   \LL(f(X)) = 4X(1-X)f''(X) + 2p-2(p+q)Xf'(X).
$$
After the change of
variable $Y=2X-1$ we obtain $4 J_{q/2,p/2}$.
All the projections of Laplace operators presented
in this section are computed by this scheme.

%%%%%%%%%%%%%%%%%%%%%%%%%%%%%%%%%%%%%%%%%%%%%%%%%%
%%%%%%%%%%%%%%%%%%%%%%%%%%%%%%%%%%%%%%%%%%%%%%%%%%
%%%%%%%%%%%%%%%%%%%%%%%%%%%%%%%%%%%%%%%%%%%%%%%%%%
%%%%%%%%%%%%%%%%%%%%%%%%%%%%%%%%%%%%%%%%%%%%%%%%%%

%                  S Q U A R E

\subsection{ The square or rectangle\label{square}}

\begin{figure}[ht!]
\begin{center}
\includegraphics[width=.3\textwidth]{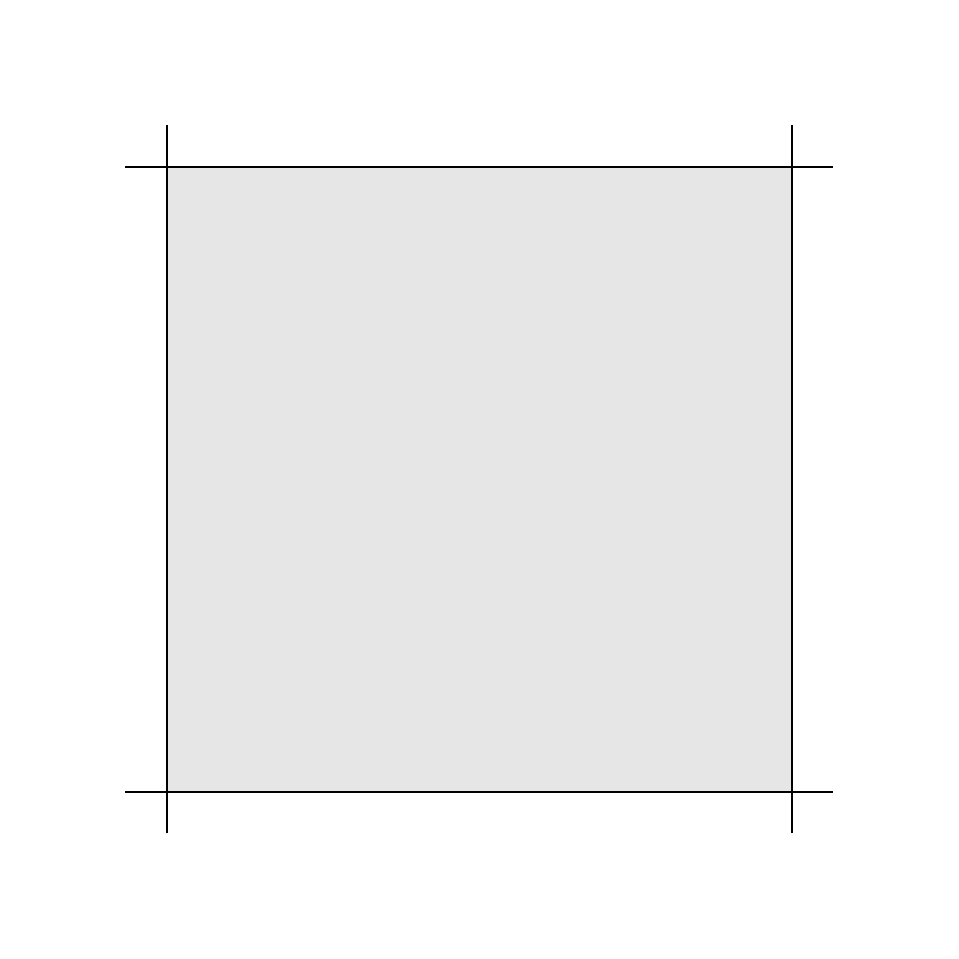}
\end{center}
\end{figure}

This is the simplest model. By affine transformation, we may chose the square  to be $[-1,1]\times [-1,1]$. The metric is
$$
   G=\bpm 1-X^2&0\\0&1-Y^2\epm
$$
and the density of the measure is 
$$
   \rho(X,Y)=C(1-X)^{a-1}(1+X)^{b-1}(1-Y)^{c-1}(1+Y)^{d-1}.
$$
This corresponds to the products of dimension $1$ Jacobi polynomials.
We recall that for any positive half-integer $p$ and $q$,
the one-dimensional Jacobi operator  $J_{p,q}$ with reversible  measure 
$(1-X)^{p-1}(1+X)^{q-1}$ can be realized on a $(2p+2q-1)$-dimensional sphere
$\{ x_1^2+ \ldots + x_{p+q}^2=1\}$
through the function 
$X= 2\big(x_1^2+ \ldots + x_{2p}^2\big)-1$. Hence we have
$$
   \De_{\bS^{2p+2q-1}}\big(h(X)\big)= 4J_{p,q}(h)(X).
$$
Since the boundary is degree four, among the admissible density measures is
$\det(G)^{-1/2}$, and the metric is then the Euclidean metric,
through the change of coordinates
$X= \cos(x_1), ~Y= \cos(x_2)$. Then, the operator is nothing else than the Laplace operator on $\bR^2$, acting on functions which are invariant under
the reflections with respect to the lines $\{x_1=k\pi\}$, 
 $\{x_2= k'\pi\}$, which is the square lattice in $\bR^2$.
Of course, this is covered by the first case since when $p=q=1/2$, the sphere is nothing else than the 1-dimensional torus. 

 Therefore, this square model for half integer values of the coefficients $(a,b,c,d)$
may be seen as images of products of spheres.  We already mentioned also the various
interpretations coming from compact  rank 1 symmetric  spaces.

 %%%%%%%%%%%%%%%%%%%%%%%%%%%%%%%%%%%%%%%%%%%%%%%%%%%%%%%%%%%
 %%%%%%%%%%%%%%%%%%%%%%%%%%%%%%%%%%%%%%%%%%%%%%%%%%%%%%%%%%%
 %%%%%%%%%%%%%%%%%%%%%%%%%%%%%%%%%%%%%%%%%%%%%%%%%%%%%%%%%%%

 %                  C I R C L E

\subsection{The circle\label{circle}}

\begin{figure}[ht!]
\begin{center}
\includegraphics[width=.3\textwidth]{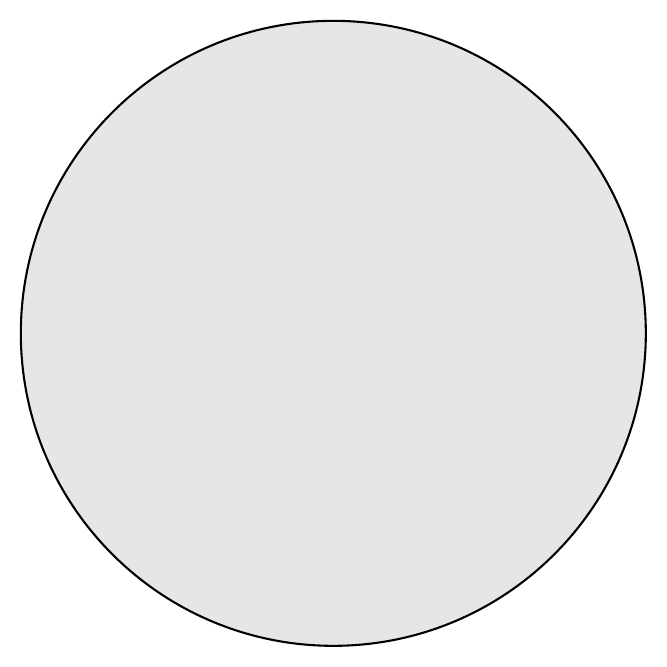}
\end{center}
\end{figure}
We may chose $\Omega$ to be the unit disk in $\bR^2$. In this case, the metric is not unique, and, up to scaling, depends  on 2 free parameters. Up to some rotation,
$$
   G_{a,b,c}= (1-X^2-Y^2) \bpm a&0\\0&b\epm + c\bpm1-X^2&-XY\\-XY&1-Y^2\epm.
$$
Ellipticity imposes $c\geq 0$ (otherwise $G_{a,b,c}^{11}(0,1)<0$), and whenever $c\neq 0$,
we may reduce by homogeneity to $c=1$.  We concentrate only on this  case.

Ellipticity condition also imposes $a>-1$, $b>-1$. When $a,b\neq 0$, the determinant
of $G_{a,b,c}$ writes
$(1-X^2-Y^2)P_2(X,Y)$, where $P_2$ has degree $2$ and is irreducible (and is constant whenever $a=b=0$). Comparing Proposition~\ref{prop.measure} and formula~\eqref{eq.rho.1}  with the value of the determinant, it is easily seen (see Remark~\ref{rmq.compute.rho})
that the only admissible measures have density
$$
   \rho_p(X,Y)= C(1-X^2-Y^2)^{p-1}\,.
$$
This remains the case even when $ab=0$ (or $c=0$),
although in these cases $P_2$ is real reducible (resp. $P_2=1-X^2-Y^2$).
In complex notation, with $Z= X+iY$, the operator associated with $c=1$ and measure  with density $C(1-X^2-Y^2)^{p-1}$ may be described from 
$$
\begin{aligned}
   &\LL_{p,a,b,1}(X) = -(1+2p+2ap)X,\quad
    \LL_{p,a,b,1}(Y) = -(1+2p+2bp)Y,\\
   &\LL_{p,a,b,1}(Z) = -\big(1+2p + (a+b)p \big) Z  - (a-b)p \bar Z,\\
   &\Ga_{a,b,1}(Z,Z) = (a-b)(1-Z\bar Z) -Z^2,\\
   &\Ga_{a,b,1}(Z, \bar Z) = (a+b+2)-(a+b+1) Z\bar Z,
\end{aligned}
$$
with of course the conjugate values for $\LL_{p,a,b,1}(\bar Z)$
and $\Ga_{a,b,1}(\bar Z, \bar Z)$.

When $a=b=0$, the metric has constant curvature $2c$. This model is well known.
For $p=1/2$, the operator corresponds to  the Laplace operator on
$\bS^2= \{x_1^2+x_2^2+x_3^2=1\}$, acting on functions which are invariant under the symmetry $x_3\mapsto -x_3$.
If one considers the unit disk as a local chart for the upper half-sphere, this is nothing else than the Laplace operator acting on functions of $(x_1,x_2)$.
The spectrum is then the spectrum of the sphere.
The eigenvalues are $\lambda_k= -k(k+1)$.

When $p= (n-1)/2$, $n\in \bN, n\geq 3$, $a=b=0$, this still corresponds to a Laplace
operator on the sphere $\bS^n$. More precisely, if one considers the Laplace operator
$\De_{\bS^n}$ on the unit sphere $\bS^n= \{x_1^2+ \ldots+ x_{n+1}^2=1\}$, and  a function depending only on $(x_1, x_2)= (X,Y)$, one gets 
$$
    \De_{\bS^n}\big(f(x_1, x_2)\big)= \LL_{(n-1)/2, 0,0,1}(f)(x_1 ,x_2).
$$
It is therefore the image of $\De_{\bS^{n}}$
through the projection $x\mapsto (x_1,x_2)$.

In this case, one may also get some other interpretation, from spheres
$\bS^{2n-1}$: on  the unit sphere $\bS^{2n-1}\subset \bR^{2n}$, $n\ge2$,
let $z_k$ be the complex functions 
which are the restrictions to the sphere  of the linear forms
$z_k(x)= x_k + ix_{n+k}$. Then, consider the complex
function $Z= z_1^2+\cdots+  z_n^2$. One can see that
$$
   \De_{\bS^{2n-1}}(Z)= -4n Z,\qquad
   \Ga_{\bS}(Z,Z)= -4Z^2,\qquad
   \Ga_{\bS}(Z, \bar Z)= 8-4Z\bar Z.
$$
Therefore, passing to the real forms $Z= X+iY$, one has
$$
  \Ga_{\bS}(X,X)= 4(1-X^2),\qquad
  \Ga_{\bS}(Y,Y)= 4(1-Y^2),\qquad
  \Ga_{\bS}(X,Y)= -4XY.
$$
This corresponds to the operator $4\LL_{(n-1)/2, 0,0,1}$.

One may also obtain similar forms using $Z= \sum_i z_iz'_i$, where $z_i$ and $z'_i$
are defined in a similar way but on the product of two spheres
$M=\bS^{2n-1}\times\bS^{2n-1}$
endowed with the product metric, which leads to $2\LL_{n-1, 0,0,1}$.
Indeed,
$$
   \De_M(Z)=(2-4n)Z,\qquad
   \Ga_M(Z,Z)=-2Z^2,\qquad
   \Ga_M(Z,\bar Z)=4-2Z\bar Z.
$$

For the other metrics, the situation is more complicated.  Still restricting to $c=1$, the condition for the metric to be non-negative on the disk is $a> -1$ and  $b> -1$. 
Even in the case $a=b$, 
the Laplace operator associated with the metric is no longer a solution to our problem,
and one may check that the metric has not constant curvature.

 If  we restrict our attention to the diagonal case $a=b$, then the equation simplifies.  Up to some scaling, the operator may be considered as the sum of the previous operator with $a=b=0$ and $\gamma(x\partial_y-y\partial_x)^2$, which corresponds to a circular Brownian motion in the plane. But we may construct this in a more geometric way as follows. 
For $-1<a<0$,  and density measure $(1-X^2-Y^2)^{p-1}$, one may consider a sphere
$\bS^n_r$ of radius $r$, where $a=-r^2/(1+r^2)$,  and dimension $n=2p+1$.
Then, we chose $e_1$ and $e_2$ two vectors in $\bR^{n+1}$ which are orthogonal
and of norm $1$, and consider  the complex  linear forms
$Z_1(x)= \lag e_1, x\rag +i \lag e_2, x\rag$, that we restrict on the sphere $\bS^n_r$.
It satisfies, for the Laplace operator on the sphere,
$$
   \De_{\bS^n_r}(Z_1)=-\frac{n}{r^2} Z_1,\;\;\;
   ~\Ga_{\bS_r}(Z_1,Z_1)= -\frac{1}{r^2} Z_1^2,\;\;\;
   \Ga_{\bS_r}(Z_1,\bar Z_1)= \frac{1}{r^2}(2- Z_1\bar Z_1).
$$

Consider now the product 
$\bS^1\times \bS^{n}_r$, with the product structure and Laplacian $\LL$.
With the function on $z= e^{i\theta}$ on $\bS^1$,  we look at the function 
$Z= zZ_1$. We have, for the product structure 
$$
   \LL (Z)= -\Big(\frac{n}{r^2}+1\Big) Z,\;\;\;
   ~ \Ga(Z,Z)=-\Big(1+\frac{1}{r^2}\Big)Z^2,\;\;\;
  ~\Ga(Z, \bar Z)= \frac{2}{r^2}+\Big(1-\frac{1}{r^2}\Big)Z\bar Z.
$$
Then, the image of $\dis\frac{r^2}{1+r^2}\LL$ through $Z$ is 
$\LL_{p, a,a,1}$ for $a=-r^2/(r^2+1)$ and $p=(n-1)/2$. 
However, the case $a\neq b$  remains mysterious.

This model has an immediate $d$-dimensional generalization. On the unit ball in $\bR^n$, one may consider the operator corresponding to the co-metric
\beq
    \label{Cercle-d-dim}
    (1-|x|^2) D(a_1, \ldots, a_n)+ (g_0^{ij}(x)),
\eeq
where $|x|$ denotes the Euclidean norm, 
$(g_0^{ij})= (\delta^{ij}-x_i x_j)$ corresponds to the projection of the spherical metric of $\bS^n$ onto a hyperplane, and $D(a_1, \ldots, a_n)$ is any diagonal matrix. In fact, it is quite easy to check with this co-metric the boundary condition
$\sum_jg^{ij}\partial_j P= -2(a_i+1)x_iP$, where $P= 1-|x|^2$ is the equation of the boundary.
The condition for the metric to be non negative on the unit ball
is again that  $a_i >-1$ for any $i$. Again, when $D(a_1, \ldots, a_n)=0$,
the choice of the measure $(1-|x|^2)^{(q-1)/2}$ corresponds to a Laplace operator
on the $(n+q)$-sphere.
 Furthermore, adding squares of infinitesimal rotations in various directions
with different coefficients provides a larger class of matrices $(g^{ij})$
solutions of the DOP problem for the boundary $|x|^2=1$.

%%%%%%%%%%%%%%%%%%%%%%%%%%%%%%%%%%%%%%%%%%%%%%%%%%%%%%%%%%%%%%%%%%%%%%%%%
%%%%%%%%%%%%%%%%%%%%%%%%%%%%%%%%%%%%%%%%%%%%%%%%%%%%%%%%%%%%%%%%%%%%%%%%%
%%%%%%%%%%%%%%%%%%%%%%%%%%%%%%%%%%%%%%%%%%%%%%%%%%%%%%%%%%%%%%%%%%%%%%%%%
%
%                T R I A N G L E

\subsection{The triangle\label{triangle}}
\begin{figure}[ht!]
\begin{center}
\includegraphics[width=.3\textwidth]{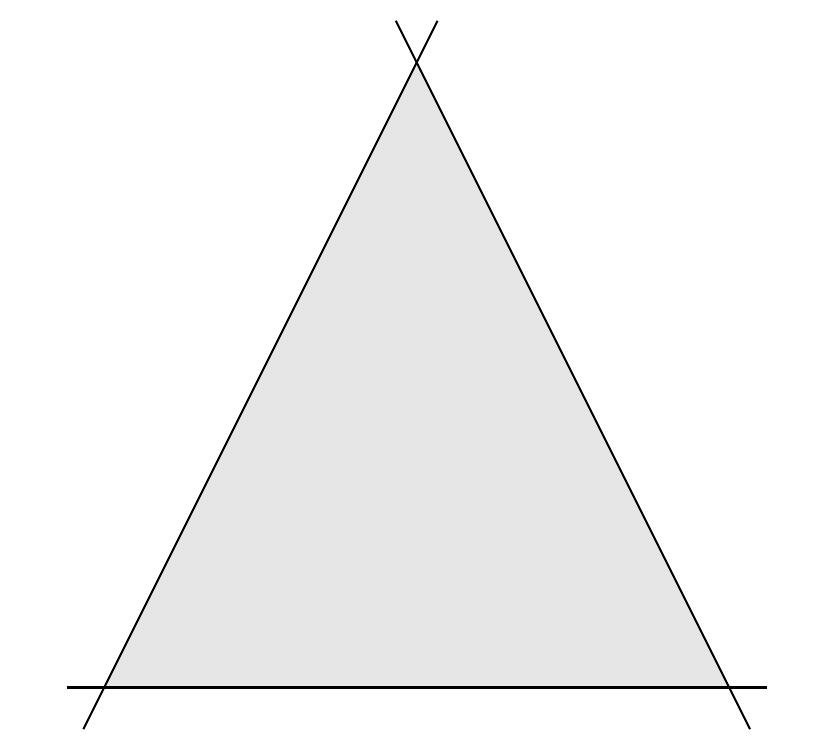}
\end{center}
\end{figure}
By affine transformation, one may reduce to the case where the triangle is delimited by the lines  $X=0$, $X+Y=1$, $Y=0$, such that the domain $\Omega$ is the 2-dimensional simplex $\{X\geq 0, Y\geq 0, X+Y\leq 1\}$. Then,the metric depends again on three  parameters
$$
   G_{a,b,c}=(1-X-Y)\bpm aX & 0\\0 & bY\epm +
            c \bpm  X( 1-X) & -XY\\ -XY & Y(1-Y)\epm.
$$
The density of the measure also depends on three parameters
$$
   \rho_{p,q,r}(X,Y)=C X^{p-1}Y^{q-1}(1-X-Y)^{r-1}\, ,
$$
which leads to a family of operators $\LL^{p,q,r}_{a,b,c}$, for which
$$
   \LL^{p,q,r}_{a,b,c}(X)=-\big((a + c)(r + p) + cq\big)X - apY + (a + c)p
$$
and a similar form for $Y$ but exchanging $X$ with $Y$,
$p$ with $q$, and $a$ with $b$.

One can check that the affine linear symmetries of the triangle correspond to simultaneous
permutations of $(p,q,r)$ and $(b+c,a+c,c)$ (cf.~Lemma~\ref{lem.3lines}), for
example, the mapping $(X,Y)\mapsto(1-X-Y,\,Y)$ transforms $\LL_{a,b,c}^{p,q,r}$
into $\LL_{a-b,c+b,-b}^{r,q,p}$ which is $\LL_{a_1,b_1,c_1}^{p_1,q_1,r_1}$ with $(p_1,q_1,r_1)=(r,q,p)$ and $(b_1+c_1,a_1+c_1,c_1)=(c,a+c,b+c)$.

This model is closely related to the circle one. We first observe that if we take
the circle model in $\bR^2$ with coordinates $(x,y)$,
and let the operator (divided by~$4$) act on  functions of $X=x^2$, $Y=y^2$,
we find the operator on the triangle acting on the variable $(X,Y)$
(the simplex is clearly the image of the disk under $(x,y)\mapsto(x^2,y^2)$).
We obtain in this way the complete family of metrics, but only the measures
$\rho_{1/2,\,1/2,\,r}$ which are the image measures of the measures on
the unit disk with density $(1-X^2-Y^2)^{r-1}$. 

For  other measures $\rho_{p,q,r}$, provided
$p$, $q$  are half-integer numbers,
one may use the $n$-dimensional model on the unit ball \eqref{Cercle-d-dim}.
As for the circle case, we restrict our study to $c=1$.
Setting $m=2p$ and $n=2p+2q$, 
consider the operator $\LL^{\bB}_{p,q,r,a,b}$
on the unit ball in $\bR^n$ given by the metric
\eqref{Cercle-d-dim} and the measure $(1-|x|^2)^{r-1}$ where
$a_1= \ldots= a_m=a$ and $a_{m+1}= \ldots= a_n=b$.
Let
$X= \sum_{i=1}^m x_i^2$, $Y= \sum_{i=m+1}^n x_i^2$. Then its image of
$\LL^{\bB}_{p,q,r,a,b}$
through $(X,Y)$ is $4\LL^{p,q,r}_{a,b,1}$,
as easily checked comparing for both cases
$\LL(X)$, $\LL(Y)$, $\Ga(X,X)$, $\Ga(X,Y)$, and $\Ga(Y,Y)$.  
 
Therefore, we see that the triangle case may be interpreted, at least for half integers values of the measure parameters, as images of the unit ball operators, in exactly the same way that one dimensional non symmetric Jacobi operators may be obtained from spheres.

Once again, those operator have an immediate $n$-dimensional extension on the $d$-dimensional simplex $\{x_i\geq 0, i= 1,\ldots, n, ~\sum_i x_i\leq 1\}$, with the (co)-metric
$$
   G^{ij}= \delta_{ij}\Big(\alpha_ix_i\big(1-x_1-\dots-x_n\big)+1\Big) - x_ix_j
$$
where $\alpha_i$ are constants.

For half-integers $p$, $q$, $r$, let us consider
the  unit sphere in $\bR^{n}$, $n=2p+2q+2r$,
where a point in $\bR^{n}$ is represented as
$(\xx_1,\xx_2,\xx_3)\in \bR^{2p}\times \bR^{2q}\times \bR^{2r}$.
Then, as we mentioned in Section~\ref{circle}, the operator
$\LL^{\bB}_{p,q,r,0,0}$ is the image of $\De_{\bS^{n-1}}$
under the projection on the first $2p+2q$
coordinates. By composing it with the mapping
of $\bB^{2p+2q}$ onto the triangle, we obtain an interpretation of
$4\LL^{p,q,r}_{0,0,1}$
as the image of $\De_{\bS^{n-1}}$ through
$(X,Y)=(\|\xx_1\|^2,\,\|\xx_2\|^2)$.
In particular, for $p=q=r=1/2$ this is the quotient of $\bS^2$ by the reflections
in the three coordinate planes. See also Remark~\ref{rmq.coax}
in Section~\ref{parab1}.

%========================= REMARK 4.1 =======================
\brmq
In the same way that we gave another interpretation on the circle coming from complex representations, one may give other interpretations on the triangle
in some particular case. For example, on $\bS^5$, consider the complex linear forms
$z_1= x_1+ ix_2$, $z_2= x_3+ i x_4$, $z_3= x_5+ix_6$
restricted to the sphere, and the function 
$$
   Z= z_1\bar z_2+ z_2\bar z_3+ z_3\bar z_1.\
$$
which maps $\bS^5$ onto the triangle in $\bC$ with vertices $1$, $\omega$, $\omega^2$
where $\omega=e^{2\pi i/3}$.
One may check that, for the Laplace operator $\De_{\bS^5}$ on the sphere, one has
$$
   \De_{\bS^5}(Z)= -12Z,\quad
   \Ga_{\bS}(Z,Z)= 4\bar Z- 4Z^2,\quad
   \Ga_{\bS}(Z, \bar Z)= 4-4Z\bar Z,
$$
which corresponds to the change of variables
$Z=1-\frac{3}{2} (X+Y)+ i \frac{\sqrt{3}}{2}(Y-X)$, i.e.,
$X=\frac13(1+\omega Z+\omega^2\bar Z)$,
$Y=\frac13(1+\omega^2 Z+\omega\bar Z)$.
Then the image of $\frac{1}{4}\De_{\bS^5}$ under $(X,Y)$ is
$\LL^{1,1,1}_{0,0,1}$.

\ermq

%%%%%%%%%%%%%%%%%%%%%%%%%%%%%%%%%%%%%%%%%%%%%%%
%%%%%%%%%%%%%%%%%%%%%%%%%%%%%%%%%%%%%%%%%%%%%%%
%%%%%%%%%%%%%%%%%%%%%%%%%%%%%%%%%%%%%%%%%%%%%%%

%         C O A X I A L     P A R A B O L A S

\subsection{The coaxial parabolas.\label{parab1}}
 
\begin{figure}[ht!]
\begin{center}
\includegraphics[width=.9\textwidth]{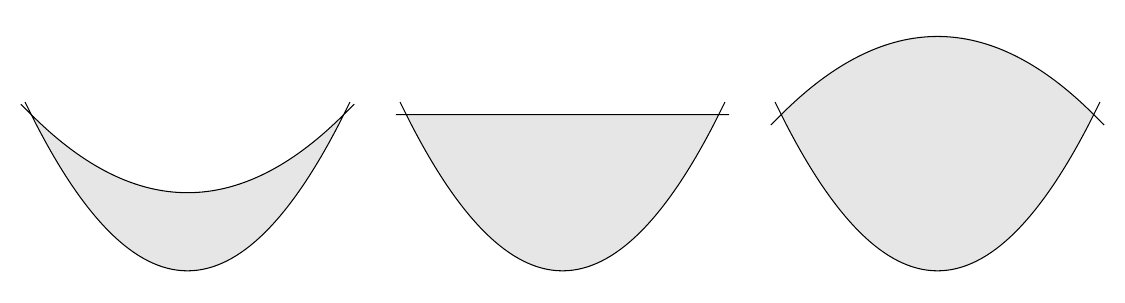}
\end{center}
\end{figure}
 
Up to affine transforms, the domain may be bounded by the two parabolas  
$Y=X^{2}-1$ and $Y=a(1-X^{2})$ with $a>-1$.
This forms a one-parameter family up to affine transformations,
but may be reduced to a single model via some non-linear transformation.

The (co)-metric is  
$$
  G_a= \bpm 1-X^{2} & -2XY     \\
         -2XY       & -4Y(1+Y) \epm
  +  4a\bpm0&0\\0&1-X^2+Y\epm.
$$
If $a\ne 0$, it is unique up to rescaling but if $a=0$, it extends to a one-parameter
family that we present at the end of this subsection.

When $a\neq 0$, the boundary has degree $4$, and therefore the Laplace operator
associated with the metric  is an admissible solution, corresponding to the measure
$\det(G_a)^{-1/2}$. It turns out that the associated metric has scalar curvature
equal to $2$, and therefore the operator is locally a spherical Laplacian. 

In fact, the (non-affine) change of coordinates
$X=X_1$, $Y= (a+1)Y_1+a(1-X_1^2)$,
allows us to reduce, up to a scaling parameter, to the case  $a=0$,
and then $\Omega=\{X^2-1<Y<0\}$ -- the domain is bounded by
the parabola $Y=X^2-1$ and the axis $Y=0$.

Even though the boundary has no longer degree 4 in this case,
the Laplace operator is still an admissible solution. In fact, the determinant
of the metric is still equal to the reduced equation of the boundary even in this case. This particular model is known in the literature as the parabolic biangle
(see Koornwinder and Schwartz \cite{KoornSchw}).

For symmetry reasons, we prefer  to consider the case $a=1$, in which case
$$
   G_1= \bpm 1-X^2& -2XY\\-2XY& 4(1-X^2-Y^2)\epm.
$$ 

When the density of the measure is $\det(G_1)^{-1/2}$, the operator may be directly seen
as the image of the Laplace operator on $\bS^2=\{x_1^2+x_2^2+x_3^2=1\}\subset \bR^3$
through $(X,Y)$, where $X= x_3$ and $Y= 2x_1x_2$.
Then, the associated operator is nothing else than the spherical Laplace operator
acting on functions which are invariant under the reflections in the hyperplanes
$\{x_1=x_2\}$ and $\{x_1=-x_2\}$
(the angle between them is $\pi/2$, which corresponds to the ordinary double
points of the boundary of the domain).

One easily checks (see Remark~\ref{rmq.compute.rho}) that all admissible
measure densities are
\beq
   \label{coax.rho}
   \rho= (1-X^2+Y)^{p-1}    (a(1-X^2)-Y)^{q-1}.
\eeq
Then we obtain an operator $\LL_{p,q,a}$ for which we have
$$
    \LL_{p,q,a}(X)= - 2(p+q)X,\qquad~
    \LL_{p,q,a}(Y)=-(2+ 4(p+q))Y+4(ap-q),
$$
When $p$ and $q$ are half-integer and $a=1$,
this operator is an image of the Laplace operator on a sphere
  $\bS^{n}\subset \bR^{n+1}$, $n=2p+2q$,
by the functions
$X = x_{n+1}$, $Y = x_1^2 + \dots x_m^2 - x_{m+1}^2 - \dots - x_n^2$ where $m=2p$.
Using formulae~\eqref{eq.lapl.spheres}, it is easily checked that they give
the required values for
$\LL_{p,q,1}(X), \LL_{p,q,1}(Y)$, $\Ga_1(X,X)$, $\Ga_1(X,Y)$, and  $\Ga_1(Y,Y)$.

Although the general case may be reduced to $a=0$,
this latter case offers a more general admissible family of (co)-metrics, namely
$$
   G'_{\alpha}= G_0 + \alpha\bpm 1-X^2+Y&0\\0&0\epm
$$
which is positive definite in the domain $X^2-1<Y<0$ for
$\alpha\geq -1$. In the special case $\alpha= -1$ we have 
$$
  G'_{-1}= -Y\bpm1 &2X\\ 2X& 4Y(1+Y)\epm.
$$
The curvature of the metric associated with $G'_{\alpha}$ is
$2(1+\alpha+2\alpha Y)/(1+\alpha+\alpha Y)^2$,
so it is non-constant for $\alpha\ne 0$.
For the measure density \eqref{coax.rho} (with $a=0$), we obtain the
operator $\LL'_{p,q,\alpha}$ with
$$
   \LL'_{p,q,\alpha}(X)=-2((1+\alpha)p+q)X,\qquad
   \LL'_{p,q,\alpha}(Y)=-(2+4(p+q))Y-4q.
$$

\brmq\label{rmq.coax2}
When $q=1/2$, the operator $\LL'_{p,\,1/2,\,\alpha}$
is the image of the operator $\LL_{p,\alpha,0,1}$
on the unit disk (see Section~\ref{circle}) under the mapping
$(X,Y)\mapsto(X,-Y^2)$.
However, as we already pointed out in Section~\ref{circle}, we do not
know any geometric interpretation for this operator when $\alpha\ne0$,
and when $\alpha=0$ we get the same models as above up to change of
coordinates.
\ermq

\brmq\label{rmq.coax}
The mapping $(X,Y)\mapsto(X^2,-Y)$ transforms $\LL'_{p,q,\alpha}$
to the operator $4L^{1/2,\,q,p}_{\alpha,0,1}$ on the triangle
(cf.~Section~\ref{triangle}).
\ermq

%%%%%%%%%%%%%%%%%%%%%%%%%%%%%%%%%%%%%%%%%%%%%%%%%%%%%%%%%%
%%%%%%%%%%%%%%%%%%%%%%%%%%%%%%%%%%%%%%%%%%%%%%%%%%%%%%%%%%
%%%%%%%%%%%%%%%%%%%%%%%%%%%%%%%%%%%%%%%%%%%%%%%%%%%%%%%%%%
%%%%%%%%%%%%%%%%%%%%%%%%%%%%%%%%%%%%%%%%%%%%%%%%%%%%%%%%%%
%%%%%%%%%%%%%%%%%%%%%%%%%%%%%%%%%%%%%%%%%%%%%%%%%%%%%%%%%%

%        Parabola + axis + tangent

\subsection{The parabola with the axis and a tangent.\label{parab2}}
\begin{figure}[ht!]
\begin{center}
\includegraphics[width=.3\textwidth]{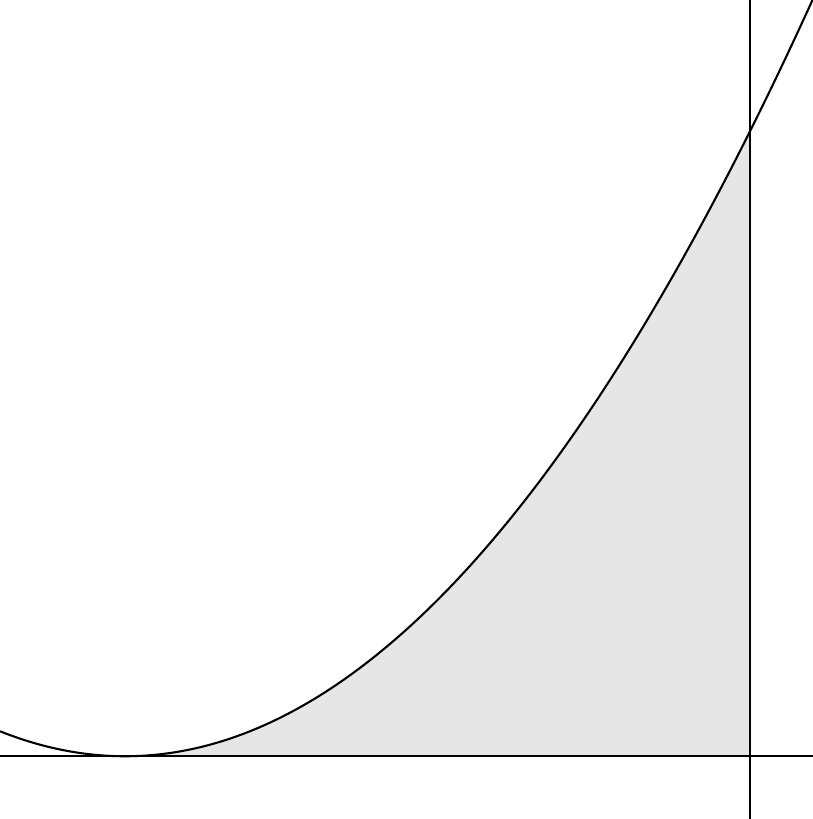}
\end{center}
\end{figure}

Notice that the axis cuts the line at infinity at the same point that the parabola.
Then, up to affine transformation, we may chose the domain $\Omega$
delimited by the curves
$$
   Y= X^2,\quad ~Y=0,\quad ~X=1.
$$
Up to scaling, there is just one (co)-metric which is a solution of the problem:
$$
    G=4\bpm X(1-X)&2Y(1-X)\\2Y(1-X)&4Y(X-Y)\epm.
$$
Once again, the boundary has degree $4$, and the Laplace  operator corresponding to the associated metric  is a solution, which corresponds to a metric with constant scalar curvature equal to 2 (that is why we chosed this normalization of the metric), and therefore  it may be realized on a unit sphere $\bS^2\subset \bR^3$.
In the general case the measure density is
 $(X^2-Y)^{p-1}Y^{q-1}(1-X)^{r-1}$, $p,q,r>0$, $p+q>1/2$.
It provides a family of operators $\LL_{p,q,r}$ for which
$$
\begin{aligned}
   \LL_{p,q,r}(X)&= 4(2p+2q-1)-4(2p+2q+r-1)X,\\
   \LL_{p,q,r}(Y)&= 16qX-8(2p+2q+r)Y.
\end{aligned}
$$
The Laplace operator  corresponds to
 $\LL_{ 1/2, 1/2, 1/2}$,
and is the image of $\De_{\bS^2}= \{x_1^2+x_2^2+x_3^2\}=1$
through $(X,Y)$ where  $X =x_{1}^{2}+ 
x_{2}^{2}$ and $Y= 4x_{1}^{2}x_{2}^{2}$.
These are functions on the sphere
invariant under the symmetries with respect to the hyperplanes
$\{x_1=0\}$, $\{x_2=0\}$, $\{x_3=0\}$, and $\{ x_1= \pm x_2\}$.
The fundamental domain on the sphere for this group action is a triangle with two $\pi/2$ angles, and one $\pi/4$ angle, which corresponds to the two nodes
and one tacnode of $\partial\Omega$.

One can check that $\LL_{1/2,\,q, r}$ is the image of the operator
$\frac4{a+c}\LL_{a,a,c}^{q,q,r}$ on the triangle (see Section~\ref{triangle})
under the mapping $(x,y)\mapsto(X,Y)=(x+y,4xy)$.
Thus each model for $\LL_{a,a,c}^{q,q,r}$ yields a model on $\Omega$ with $p=1/2$. 
In particular,
$\LL_{1/2,\,q,r}$
is the image of the Laplace orerator on the unit sphere
in $\bR^{2q}\times\bR^{2q}\times\bR^{2r}$
under the composition
$$
   (\xx_1,\xx_2,\xx_3)\mapsto\big(\|\xx_1\|^2,\|\xx_2\|^2\big)
     \mapsto\big(\|\xx_1\|^2+\|\xx_1\|^2,\;4\,\|\xx_1\|^2\,\|\xx_2\|^2\big).
$$

\def\cd{c}
For $m\geq 2$ and $\cd\in\{1,2,4,8\}$, we may construct the operator as an image
of a sphere of an appropriate dimension, namely,
of the unit sphere in $\bR^{2\cd m+2r}$.
For  those values of $\cd$, we may construct
$\cd$ orthogonal transformations $\ell_i$ on $\bR^{\cd}$ such that
$\ell_1(\uu), \dots, \ell_{\cd}(\uu)$ form an orthonormal basis
for any $\uu\in \bR^{\cd}$,  with $\|\uu\|=1$.
This is done through the complex, quaternionic or
octionionic multiplications (say from the left) by the basis elements of the algebra,
which provides orthonormal transformations of the space which satisfy the required
conditions (although in the octonionic case it is not just a simple application
of the algebra rule due to the non-associativity of the product), 
see Conway and Smith \cite{ConwaySmith}.
Indeed, this property fails for higher order Cayley-Dickson algebras.
For any $m$, the operators $\ell_i$ lift to
$\bR^{m}\otimes \bR^{\cd}$ into orthogonal transformations such that
$\ell_1(\xx), \dots, \ell_{\cd}(\xx)$ are pairwise orthogonal for any
$\xx\in\bR^m\otimes\bR^{\cd}$.

Let $n=2\cd m+2r$.
We consider a point in $\bR^n$ as a triple of vectors
$(\xx,\yy,\zz)$ with
$\xx,\yy\in\bR^{m}\otimes \bR^{\cd}$ and $\zz\in\bR^{2r}$.
Then we consider $X= \|\xx\|^2+\|\yy\|^2$ and 
$$
   Y=4\Big( \|\xx\|^2\|\yy\|^2 -
  \sum_{i=1}^{\cd} (\xx\cdot \ell_i(\yy))^2 \Big)
$$
where $\xx\cdot\yy$ denotes the usual scalar product in 
$\bR^{\cd m}$, and $\|\xx\|^2=\xx\cdot\xx$
(if $m=1$,then $Y=0$; this is why we imposed the restriction $m\ge 2$).

It may be checked that the restriction of the functions $X$ and $Y$ to the unit
sphere  in $\bR^n$ satisfy the relations required for
$\Ga(X,X)$, $\Ga(X,Y)$ and $\Ga(Y,Y)$. Indeed, once we have remarked that
$X$ and $Y$ are homogeneous with degree respectively $2$ and $4$ in
$\bR^{n}$,
and for this value of $X$,
by \eqref{lapl.sph.homog}
 everything boils down to verify that, for the Euclidean
operator $\Ga_{\bE}$ in $\bR^{n}$, one has $\Ga_{\bE}(Y,Y)= 16XY$,
which is quite easy to check. Then, one also checks that
$$
    \De_{\bS^{n-1}}(X)= 4\cd m-2nX,\quad
   ~\De_{\bS^{n-1}}(Y)= 8\cd(m-1)X-4(n+2)Y,
$$
wich corresponds to 
$\LL_{p,q,r}$ with $2p=\cd+1$ and $2q=\cd(m-1)$ (recall that $n=2\cd m+2r$).

%%%%%%%%%%%%%%%%%%%%%%%%%%%%%%%%%%%%%%%%%%%%%%%%%%%%%%%%%%%%%%
%%%%%%%%%%%%%%%%%%%%%%%%%%%%%%%%%%%%%%%%%%%%%%%%%%%%%%%%%%%%%%
%%%%%%%%%%%%%%%%%%%%%%%%%%%%%%%%%%%%%%%%%%%%%%%%%%%%%%%%%%%%%%

%   P A R A B O L A    W I T H    T W O    T A N G E N T S

  \subsection{The parabola with two tangents\label{sec.B2} }
  \begin{figure}[ht!]
\begin{center}
\includegraphics[width=.3\textwidth]{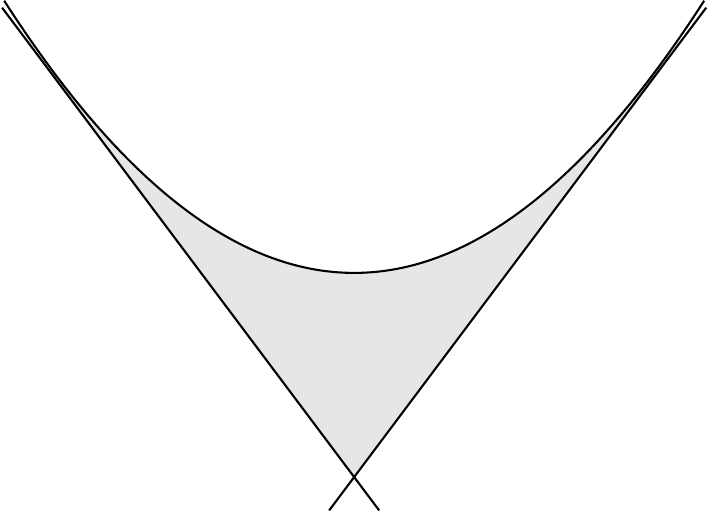}
\end{center}
\end{figure}
  
Here,  the domain $\Omega$ is delimited   by the equations
$$
      Y= X^2,\quad ~Y=2X-1,\quad ~Y=-2X-1.
$$
With this boundary, up to scaling, the (co)-metric is unique and is
$$
   G=\bpm Y+1-2X^2  &  2X(1-Y)      \\
            2X(1-Y)  &  4(2X^2-Y-Y^2)\epm.
$$
Once again, the boundary has degree $4$,  the Laplace  operator corresponding to this
(co)-metric  is a solution, and has constant  curvature $0$. The general density
measure is $(X^2-Y)^{p-1}(Y-2X+1)^{q-1}(Y+2X+1)^{r-1}$ with $p,q,r>0$, $p+q>1/2$,
$p+r>1/2$. For such measure we get an operator $\LL_{p,q,r}$ with
$$
\begin{aligned}
   \LL_{p,q,r}(X)&= 2(r-q)-2(2p+q+r-1)X,\\
   \LL_{p,q,r}(Y)&= -2(2p-1)+4(r-q)X-2(2p+2q+2r-1)Y.
\end{aligned}
$$

When $p=q=r= 1/2$, this corresponds to the image of a 2-dimensional
Euclidean Laplacian, constructed from the root system $B_2$ as follows.
Consider in $\bR^2$, with canonical basis $(e_1,e_2)$,  the 4 roots
$\lambda_j =\pm\sqrt{2} e_i$, and the 4 roots
$\mu_j= \pm\sqrt{2} e_i\pm \sqrt{2} e_j$ (the factor $\sqrt{2}$ is there
to fit with the final values of $X$ and $Y$). Then, let
$$
\begin{aligned}
  &X(x,y)=\frac{1}{4} \sum_{j= 1}^4 \exp(i \lambda_j\cdot(x,y))
  = \big(\cos(\sqrt{2}x)+\cos(\sqrt{2}y)\big)/2,
\\
  &Y(x,y)
= \frac{1}{4} \sum_{j=1}^4 \exp(i\mu_j\cdot(x,y))
= \cos(\sqrt{2} x)\cos(\sqrt{2}y).
\end{aligned}
$$
Then, it is directly checked that $\Ga_{\bR^2}(X,X)$, $\Ga_{\bR^2}(X,Y)$,
$\Ga_{\bR^2}(Y,Y)$, $\De_{\bR^2}(X)$, $\De_{\bR^2}(Y)$ satisfy the relations
required for $\LL_{1/2,\,1/2,\,1/2}$. This is just one example of the family of
Jack polynomials associated with root systems (see MacDonald \cite{Macdo1}).
Following Koornwinder \cite{Koorn1,Koorn2}, one may find other representations
for symmetric rank 2 spaces with restricted root systems $B_2$
(which include for example $SO(5)$ and \hbox{$SO(n+2)/SO(n)$}). For a reference on this model,
see also Sprinkhuizen-Kuyper \cite{Sprink}. For the sake  of completeness,
we give below some naive representations of those  models coming from the
Laplace-Beltrami operator on $SO(n)$ described in~\eqref{lapl.so(n)}.
One may find more complete descriptions of those models in Doumerc's thesis
\cite{YanDoumerc}. Moreover, this allows us to show how to deal in a convenient
way with matrix operators.

For a given operator on square matrices in dimension $n$, such as the one described
in~\eqref{lapl.so(n)} or~\eqref{lapl.su(n)1} and~\eqref{lapl.su(n)2}, one may consider
the image of the operator on the spectrum, determined by the coefficients
$a_0, \dots, a_{n-1}$ of the characteristic polynomial
$P(\lambda)= \det(M-\lambda\,\Id)$. Of course, for small values of $n$,
one may perform computations by hand, but it is perhaps worth to describe
general methods. 

The first task is to compute the various derivatives with respect to the entries
$M_{ij}$ of $M$ of the various coefficients of $P(\lambda)$. One may start from
the comatrix $\hat M= \hat M_{ij}$ for which  $\hat M^t= \det(M) M^{-1}$
(where $M^t$ is the transposed of $M$) and satisfies $\partial_{M_{ij}}\hat M_{ik}=0$.
Together with $\partial_{M_{ij}} M^{-1}_{kl} = -M^{-1}_{ki}M^{-1}_{jl}$,
we get  $\partial_{M_{ij}}\log \det(M)= M^{-1}_{ji}$
(which is valid on the dense domain where $\det(M)\neq 0$).  

Now, for an operator on matrices satisfying $\LL(M_{ij})= -\mu M_{ij}$ and
$\Ga(M_{ij}, M_{kl})= \delta_{ik}\delta_{jl}- M_{il}M_{kj}$, 
denoting $M(\lambda)$ the matrix $M-\lambda\,\Id$, the previous formulae
combined with the change of variable formula~\eqref{chgt.de.variables} leads to
$$
\begin{aligned}
 \LL\big(\log P(\lambda)\big)
  =&-\mu\,\tr\big(M(0)M(\lambda)^{-1}\big)
     -\tr\big(M(\lambda)^{-1}M^t(\lambda)^{-1}\big)\\
   &+
 \Big(\tr\big(M(0)M(\lambda)^{-1}\big)\Big)^2
\end{aligned}
$$
$$
\begin{aligned}
 \Ga\Big(\log\big(P(\lambda_1)\big),\, \log\big( P(\lambda_2)\big)\Big)
  =&\; \tr\Big(M^t(\lambda_1)^{-1}M(\lambda_2)^{-1}\Big)\\
   &\,-\tr\Big(M(\lambda_1)^{-1}M(0)M(\lambda_2)^{-1} M(0)\Big).
\end{aligned}
$$

For the special case of $SO(n)$  where $\mu= n-1$ and $M^t= M^{-1}$, denoting
$a_1,\dots,a_n$ the eigenvalues of $M$, one has
$$
\begin{aligned}
   &\tr\big(M(0)M(\lambda)^{-1}\big) = \sum\frac{a_i}{a_i-\lambda}
   =\sum\Big(1+\frac{\lambda}{a_i-\lambda}\Big)
      = n - \lambda\frac{P'(\lambda)}{P(\lambda)}\,,
\\
   &\tr\big(M(\lambda)^{-1}M^t(\lambda)^{-1}\big)
   =\sum\frac{1}{(a_i-\lambda)(a_i^{-1}-\lambda)}
\\
   &\qquad=\frac1{1-\lambda^2}\sum\Big(1+\frac{\lambda}{a_i-\lambda}
     +\frac{\lambda}{a_i^{-1}-\lambda}\Big)
    =\frac1{1-\lambda^2}\Big(n-2\lambda\frac{P'(\lambda)}{P(\lambda)}\Big)\,.
\end{aligned}
$$
Putting this into the previous formula, we obtain
$$
  \De_{SO(n)}\big(\log P(\lambda)\big)
  =-\frac{n \lambda^2}{1-\lambda^2}
   +\lambda\frac{P'(\lambda)}{P(\lambda)}\Big(\frac{1+\lambda^2}{1-\lambda^2}-n\Big)
   +\Big(\lambda \frac{P'(\lambda)}{P(\lambda)}\Big)^2\,.
$$
Similarly,
$$
\begin{aligned}
  \Ga_{SO(n)}\big(\log P(\lambda_1),\,\log P(\lambda_2)\big)
  =&\frac{1}{1-\lambda_1\lambda_2}\Big(n\lambda_1\lambda_2
     -\lambda_1\frac{P'(\lambda_1)}{P(\lambda_1)}
     -\lambda_2\frac{P'(\lambda_2)}{P(\lambda_2)}\Big)\\
   &+\frac{1}{\lambda_1-\lambda_2}\Big( \lambda_1^2\frac{P'(\lambda_1)}{P(\lambda_1)}
     -\lambda_2^2\frac{P'(\lambda_2)}{P(\lambda_2)}\Big)
\end{aligned}
$$
(if $\lambda_1=\lambda_2=\lambda$, the second term is
$\partial_\lambda(\lambda^2P'/P)$).
Using \eqref{chgt.de.variables} with $\Phi=\exp$,
this leads to the very simple formulas
\beq
\label{PaTanTan.SO}
\begin{aligned}
   \De_{SO(n)}\big(P(\lambda)\big)&= -(n-1) \lambda P' + \lambda^2P'',
\\
   \Ga_{SO(n)}\big(P(\lambda),P(\lambda)\big)
   &= \lambda^2\Big(\big(n P^2 - 2\lambda PP'\big)/(1-\lambda^2)
  + PP'' -(P')^2\Big). 
\end{aligned}
\eeq
For $n=4$, we write $P(\lambda)= \lambda^4+X\lambda^3+ Y\lambda^2 + X\lambda + 1$.
Plugging this expression to the right hand side of \eqref{PaTanTan.SO}
and comparing the coefficients of powers of $\lambda$ with those in
the expansions
$$
\begin{aligned}
   &\De(P)=\De(X)\lambda+\De(Y)\lambda^2+\dots\,,\\
   &\Ga(P,P)=\Ga(X,X)\lambda^2 + 2\Ga(X,Y)\lambda^3
         + \big(2\Ga(X,X)+\Ga(Y,Y)\big)\lambda^4
                  +\dots
\end{aligned}
$$
one gets
$$
   \De_{SO(4)}(X)= -3X, \qquad\De_{SO(4)}(Y)= -4Y.
$$
and 
$$
\begin{aligned}
   &\Ga_{SO(4)}(X,X)= 4-X^2+2Y,\\
   &\Ga_{SO(4)}(X,Y)=6X-XY,\\
   &\Ga_{SO(4)}(Y,Y)= 8+4X^2-2Y.
\end{aligned}
$$
From this, we see that $\LL_{3/2,\,1/2,\,1/2}$ is the image of $2\De_{SO(4)}$
through $(X_1,Y_1)$ where $4X_1=X$ and $4Y_1=Y-2$.

For $SO(5)$, setting $P(\lambda)= \lambda^5+ X\lambda^4+ Y\lambda^3+
Y\lambda^2+ X\lambda+1$, and $X= 4X_2+1$, $Y= 4X_2+4Y_2+2$, one sees
with the same method that the image of $2\De_{SO(5)}$ through
$(X_2,Y_2)$ is  $\LL_{3/2,\, 3/2,\, 1/2}$.

One may also project $\De_{SO(n)}$ on any $m\times s$ submatrix $\hat M$ (it is obvious from formulae~\eqref{lapl.so(n)} that the operator projects). It is less obvious 
a priori (but still easy to check using~\eqref{lapl.so(n)})
that it also projects on the square  $s\times s$  matrices $N=\hat M^t \hat M$, and  produces on the entries $N_{ij}$ of those matrices the operator defined by
$\De_{SO(n)}(N_{ij})= 2m \delta_{ij} - 2nN_{ij}$ and 
$$
  \Ga_{SO(n)}(N_{ij}, N_{kl})= N_{ik}\delta_{jl}+N_{il}\delta_{jk}+N_{jk}\delta_{il}
                              +N_{jl}\delta_{ik}-2\big(N_{ik}N_{jl}+N_{il}N_{jk}\big).
$$
Again, this projects on the spectrum of such matrices. In particular,
when $s=2$, $m=2r+1$, and $n=2q+2r+2$ for positive half-integers $p$ and $q$,
one may chose as variables $\tr(N)= X+1$ and $4\det(N)= Y+2X+1$, and then
the image of $\frac{1}{2}\De_{SO(n)}$ through
$(X,Y)$ is $\LL_{1,q,r}$. For $r=0$ (thus $m=1$), the image is obviously degenerate,
and concentrated on the boundary  $\{Y+2X+1=0\}$, while for $q=0$ (thus $n=m+1$),
it concentrates on $\{Y-2X+1=0\}$, as would do the image measure
when $r\to 0$ or $q\to 0$ respectively.

\brmq\label{rmq.CuCuTan}
The singularities of $\Omega$ correspond to the angles
$(\frac\pi2,\frac\pi4,\frac\pi4)$.
This is a Euclidean triangle which can be obtained by folding a square along the
diagonal. This corresponds to the mapping $[-1,1]^2\to\Omega$ given
by $(X,Y)\mapsto \big(\frac12(X+Y),XY\big)$ which maps the lines
$X\pm1=0$ and $Y\pm1=0$ to the line $Y\pm2X+1$ and the diagonal $X=Y$ to $Y=X^2$.
This mapping transforms the product of Jacobi operators
$J_{q,r}\times J_{q,r}$ (see Section~\ref{square}) to $\frac12\LL_{1/2,\,q,r}$.
In particular, $\LL_{1/2,\,q,r}$ can be interpreted in this way as an appropriate
projection of $\De_{\bS^n\times\bS^n}$.

Similarly, if we fold the triangle with angles $(\frac\pi2,\frac\pi4,\frac\pi4)$ along
its axis of symmetry, we obtain a triangle with the same angles.
This corresponds to the 2-to-1 mapping $\Omega\to\Omega$ given by
$(X,Y)\mapsto(Y,4X^2-2Y-1)$. Under this mapping,
the parabola $\{Y=X^2\}$ is mapped to the line $\{Y-2X+1=0\}$, the both lines
$\{Y\pm2X+1=0\}$ are mapped to the parabola $\{Y=X^2\}$, and the axis $\{X=0\}$
is mapped to the line $\{Y+2X+1=0\}$.
This mapping transforms $\LL_{p,q,q}$ into $2\LL_{q,p,1/2}$.
In particular, the above model on $SU(4)$ is transformed into $\LL_{1/2,\,3/2,\,1/2}$,
the model of $\LL_{1,q,q}$ via $SU(4q+2)$ is transformed into
$\LL_{q,\,1,\,1/2}$, and the model coming from $J_{p,p}\times J_{p,p}$ is transformed into
$\LL_{p,\,1/2,\,1/2}$.
\ermq

%%%%%%%%%%%%%%%%%%%%%%%%%%%%%%%%%%%%%%%%%%%%%%%%%%%%%%%%%%%%%%%%%%%%
%%%%%%%%%%%%%%%%%%%%%%%%%%%%%%%%%%%%%%%%%%%%%%%%%%%%%%%%%%%%%%%%%%%%
%%%%%%%%%%%%%%%%%%%%%%%%%%%%%%%%%%%%%%%%%%%%%%%%%%%%%%%%%%%%%%%%%%%%
%%%%%%%%%%%%%%%%%%%%%%%%%%%%%%%%%%%%%%%%%%%%%%%%%%%%%%%%%%%%%%%%%%%%

%                    N O D A L    C U B I C

\subsection{The nodal cubic\label{sec:DlePtCub}}

\begin{figure}[ht!]
\begin{center}
\includegraphics[width=.3\textwidth]{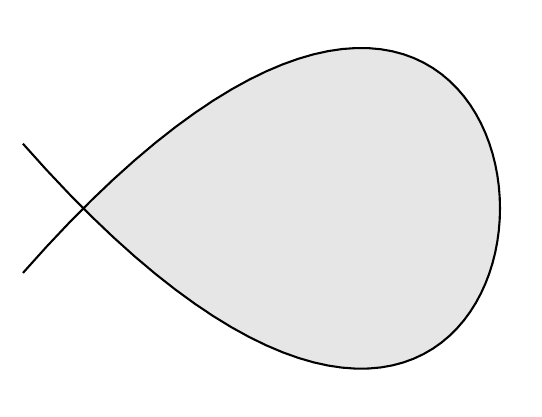}
\end{center}
\end{figure}
In this situation, we may choose the equation of the boundary  to be $Y^{2}= X^{2}(1-X)$. 
There  is a unique metric up to scaling
$$
    G= \bpm 4X(1- X)  &  2Y(2-3X)        \\
            2Y(2-3X)  & 4X-3X^{2}-9Y^{2} \epm.
$$
The boundary has degree 3, and in this situation the measure  density $\rho(x)= \det(G)^{-1/2}$ is not an admissible measure (it does not satisfy equationı \eqref{eq.rho.1}, as one may check directly).  Also, the metric has a non constant curvature. The general form of the density measure is 
$\rho_p(X,Y)= \big(X^2(1-X)-Y^2\big)^{p-1} $, for which we have
\begin{equation*}
\LL(X)=-2(6p+1)X + 8p\,,\qquad \LL(Y)=-6(3p+1)Y\,.
\end{equation*}

 It turns out that for $p=1/2$, the operator may be interpreted from a
$3$-dimensional sphere, through a projection which is very close to the
Hopf fibration. Indeed, on the unit sphere
$\bS^3= \{x_1^2+ x_2^2+ x_3^2+ x_4^2=1\}$,
consider the functions
$X=x_1^2+ x_2^2$ and
$Y=\big(x_1^2-x_2^2\big)x_3 + 2x_1 x_2 x_4$.
We may check directly that they satisfy the required equations on
$\De_{\bS^3}(X)$, $\De_{\bS^3}(Y)$, $\Ga_{\bS^3}(X,X)$, $\Ga_{\bS^3}(X,Y)$,
and $\Ga_{\bS^3}(Y,Y)$.

\begin{figure}[ht]

\centering \includegraphics[width=.40\linewidth]{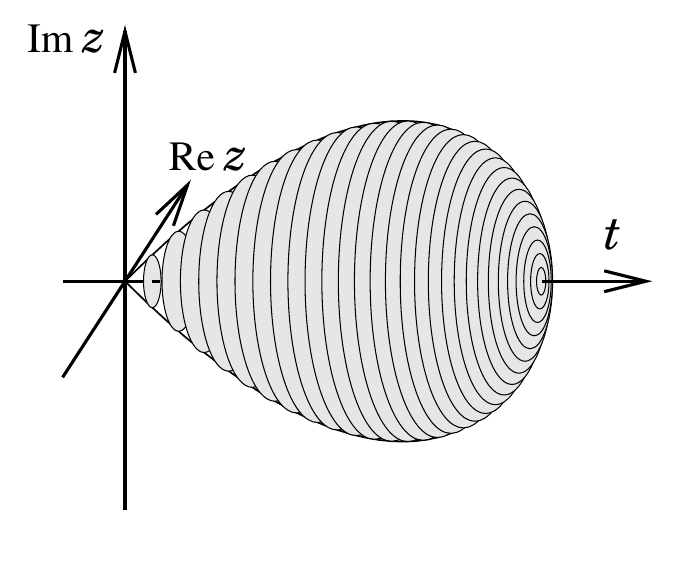}
\caption{The surface of revolution over the nodal cubic}
		\label{fig:revol}
\end{figure}

To understand which functions on the sphere are of the form $f(X,Y)$, one may represent
the sphere in complex notation as $\{|z_1|^2+ |z_2|^2=1\}$,
where $(z_1, z_2)\in \bC^2$, that we write in polar coordinates as
$z_j= \rho_j\exp(i\theta_j)$.  We then see that 
$$
   (X,Y)= \big(|z_1|^2,\, \Re(z_1^2\bar z_2)\big) =
   \big(\rho_1^2,\, \rho_1^2\rho_2\cos(2\theta_1-\theta_2)\big).
$$
Then $(X,Y)$ is invariant under 
$$
   (z_1,z_2)\mapsto (e^{i\theta}z_1, e^{2i\theta}z_2).
$$
Moreover, the quotient of the sphere under this action 
can be identified with the image $\Omega_1=F(\bS^3)$ of the sphere under
the mapping $F:\bS^3\to\bR_+\times\bC$,
$(z_1,z_2)\mapsto(|z_1|^2, z_1^2\bar z_2)= (t,z)$. 
The functions $t=\rho_1^2$ and $r=|z|=\rho_1^2\rho_1$ satisfy
the relation  $r^2= t^2(1-t)$, and therefore $\Omega_1$ is the surface of
revolution, with axis $\bR\times\{0\}\subset\bR\times\bC$, whose
meridional section is $\partial\Omega$ placed in the real half-plane
$\{(t,z)\mid\Im\,z=0\}$; see Figure~\ref{fig:revol}. So, $\Omega$ is the quotient
of $\Omega_1$ by the symmetry $(t,z)\mapsto(t,\bar z)$.

\smallskip

This construction admits the following generalization for
some other values of the measure parameter. Namely, for $p=\cd/2$
where $\cd\in\{1,2,4,8\}$ (with $\cd=1$ corresponding to the 
considered model on $\bS^3$). 
We shall use an interpretations similar to the one described
in Section~\ref{parab2}. Let us write a point in
$\bR^{3\cd+1}$ as $(\uu,\vv,\ww)$, with $\uu, \vv\in \bR^{\cd}$ and
$\ww= (w_0, w_1, \dots, w_{\cd})\in \bR^{\cd+1}$. Consider then
$X= \|\uu\|^2+ \|\vv\|^2$. For these values of $c$, as seen above, 
there exist $c$ orthogonal transformations $\ell_k$  in $\bR^{\cd}$ such that
$\{\ell_1(\vv), \ell_2(\vv), \dots, \ell_{\cd}(\vv)\}$
is  an orthonormal basis for any unit vector $\vv\in \bR^{\cd}$.
Then, on $\bR^{2\cd}$, one considers the bilinear functions
$B_k(\uu,\vv)= 2\uu\cdot \ell_k(\vv)$, $k=1,\dots,\cd$,
for which it is immediate that
$\sum_{j=1}^{\cd} B_j^2= 4\|\uu\|^2\|\vv\|^2$.
Let also $B_0(\uu,\vv)= \|\uu\|^2-\|\vv\|^2$. Then
$\sum_{i=0}^{\cd} B_i^2= X^2$.
For the Euclidean Laplacian on $\bR^{2\cd}$, one has
$$
  \De_\bE B_i=0,\qquad
  \Ga_\bE(B_i,B_j) = 4\delta_{ij}(\|\uu\|^2+\|\vv\|^2),\qquad
  i,j= 0, \dots, \cd. 
$$

We then consider the function $Y= \sum_{i=0}^{\cd} w_iB_i$.
For the Euclidean Laplace operator in $\bR^{3\cd+1}$, one easily checks that
$\Ga_\bE(X,Y)= 4Y$ and that
$$
    \Ga_\bE(Y,Y)= X^2+4X\|\ww\|^2=X^2+4X(1-X).
$$
The comparison~\eqref{lapl.sph.homog} of spherical Laplace operator and the Euclidean
one  shows that the restrictions of $X$, $Y$, and the Laplace operator on the unit
sphere in $\bR^{3\cd+1}$ satisfy the required relations for  $\LL_p$ with $2p=\cd$.

It is perhaps worth to observe that in the above construction,
the bilinear functions $B_0, B_1, \dots , B_{\cd}$, considered as functions on
$\bR^{2\cd}$ are harmonic and satisfy
$\Ga_\bE(B_i,B_i)= 4\delta_{ij} (\|\uu\|^2+ \|\vv\|^2)$,  and
$\sum_iB_i^2=(\|\uu\|^2 +\|\vv\|^2)^2$ . Their restriction to
$\bS^{2\cd-1}$ satisfy  then the same relations (up to some factor $4$) than the coordinates on a unit sphere  $\bS^{\cd}$.
Any construction performed on those spheres may be then  carried to $\bS^{2\cd-1}$,
just replacing $X_i$ by $B_i$.

%%%%%%%%%%%%%%%%%%%%%%%%%%%%%%%%%%%%%%%%%%%%%%%%%%%%%%%%%%%%%%%%%%%%%%%%%%%%%%%%%%%%%
%%%%%%%%%%%%%%%%%%%%%%%%%%%%%%%%%%%%%%%%%%%%%%%%%%%%%%%%%%%%%%%%%%%%%%%%%%%%%%%
%%%%%%%%%%%%%%%%%%%%%%%%%%%%%%%%%%%%%%%%%%%%%%%%%%%%%%%%%%%%%%%%%%%%%%%%%%%%%%%
%%%%%%%%%%%%%%%%%%%%%%%%%%%%%%%%%%%%%%%%%%%%%%%%%%%%%%%%%%%%%%%%%%%%%%%%%%%%%%%
%%%%%%%%%%%%%%%%%%%%%%%%%%%%%%%%%%%%%%%%%%%%%%%%%%%%%%%%%%%%%%%%%%%%%%%%%

%                C U B I C   +   S E C A N T

\subsection{The cuspidal cubic with one secant line\label{sec:CuCuSec}}

\begin{figure}[ht!]
\begin{center}
\includegraphics[width=.3\textwidth]{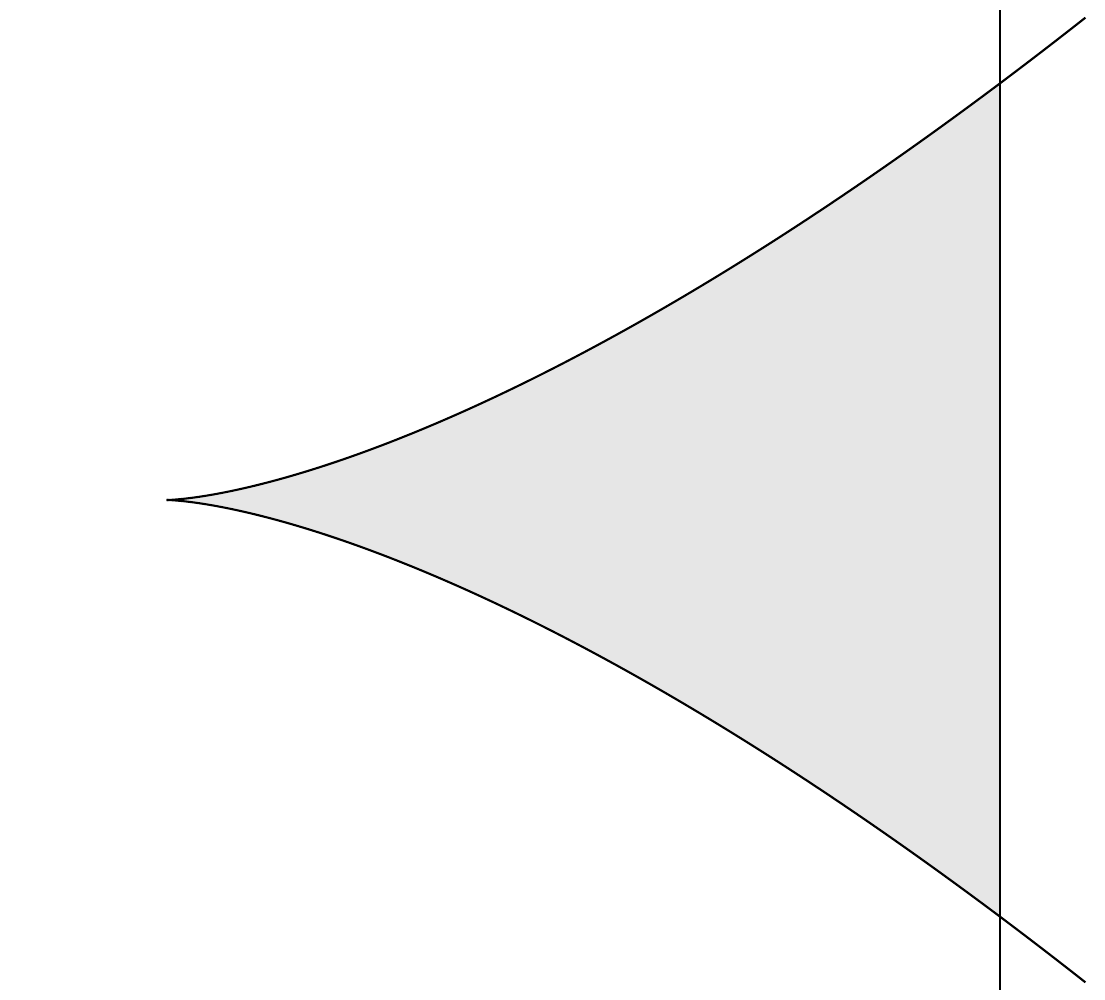}
\end{center}
\end{figure}
We may choose the boundary equation to be $(X^3-Y^2)(X-1)=0$.
Up to scaling, the associated metric is unique and we have
$$
   G= \bpm 4X(1-X)  &   6Y(1-X)        \\
           6Y(1-X)  &   9(X^{2}-Y^{2}) \epm.
$$
Since the boundary has degree 4, the Laplace operator associated with this metric
belongs to the admissible solutions and we may check that the associated metric has
constant scalar curvature $2$ and therefore may be realized from the unit sphere
$\bS^2$. 

The general density measure is $\rho_{p,q}= (X^3-Y^2)^{p-1}(1-X)^{q-1}$, $p>1/6$,
for which we have for the associated operator $\LL_{p,q}$
$$
   \LL_{p,q}(X)= -2(6p+2q-1)X + 2(6p-1) ,\qquad
   \LL_{p,q}(Y)= -3(6p+2q)Y.
$$
For the Laplacian case, $\LL_{1/2,\,1/2}$ is the image of $\De_{\bS^2}$ through 
$$
    X= x_1^{2}+x_2^{2}, \qquad
    Y= x_1\big(x_1^{2}-3x_2^{2}\big).
$$
 The functions $F(X,Y)$ are the functions on the unit sphere which are invariant under
 $x_3\mapsto -x_3$ and such that the projection $z= x_1+ix_2= \rho e^{i\theta}$ on the
 hyperplane $\{x_3=0\}$ depend only on $\rho$ and $\cos(3\theta)$. These are the
 functions which are invariant under symmetries through the hyperplanes $H=\{x_2=0\}$
 and the two hyperplanes having an angle $\pm\pi/3$ with $H$. The  fundamental domain
 for these symmetries  on the sphere is a triangle with angles
 $(\pi/3, \pi/2, \pi/2)$, which correspond to one cusp and two double points.

For the other density measures, we may
 consider the unit sphere in
 $\bR^{3\cd+2}\times \bR^{2q}$ where $2p=\cd+1$.
 For a point $(\uu,\vv)\in \bR^{3\cd+2}\times \bR^{2q}$, we set $X= \|\uu\|^2$
 and we chose for $Y$ some homogeneous degree 3 harmonic polynomial $P(\uu)$.
 Then, the required formulae for $\LL_{p,q}(X)$, $\LL_{p,q}(Y)$, $\Ga(X,X)$,
 $\Ga(X,Y)$ and $\Ga(Y,Y)$ are satisfied as soon as $\Ga_\bE(Y,Y)= 9\|\uu\|^4$,
 where $\Ga_\bE$ denotes the Euclidean operator $\Ga$. 

This problem has been studied by Cartan \cite{cartan1} where he proved that such
 polynomials exist only for $\cd=0,1,2,4,8$. Beyond the case $\cd=0$
 (the above example),
 this corresponds respectively to real, complex, quaternionic and octonionic structures.
 Such a function (for $\cd=1,2,4,8$ ) may be for example represented as follows:
 consider a Hermitian $3\times 3$ matrix with trace $0$ and respectively
 real, complex, quaternionic, or octonionic entries.
 On this space of matrices, one may consider the Euclidean
 structure given by $X=\|M\|^2= \tr(M^*M)$, and, for this structure, the function
 $Y: M\mapsto 3\sqrt6\det(M)$, satisfies $\|M\|^2=\tr(M^2)$,
 as one may check by direct computation.
 The case $p=0$
 corresponds to diagonal matrices. 

\brmq\label{rmq.DetH}
The determinant of a $3\times3$ Hermitian matrix over $\bH$ or $\bO$, in fact,
does not make much troubles.
Indeed, only two terms of its expansion depend on the order of multiplication:
$M_{12}M_{23}M_{31}$ and $M_{13}M_{32}M_{21}$. So, we can choose any order
for one of them and take the conjugate value for the other one.
\ermq

%%%%%%%%%%%%%%%%%%%%%%%%%%%%%%%%%%%%%%%%%%%%%%%%%%%%%%%%%%%%%%%%%%%%%%%%%%%%%%%%%%%%%
%%%%%%%%%%%%%%%%%%%%%%%%%%%%%%%%%%%%%%%%%%%%%%%%%%%%%%%%%%%%%%%%%%%%%%%%%%%%%%%
%%%%%%%%%%%%%%%%%%%%%%%%%%%%%%%%%%%%%%%%%%%%%%%%%%%%%%%%%%%%%%%%%%%%%%%%%%%%%%%
%%%%%%%%%%%%%%%%%%%%%%%%%%%%%%%%%%%%%%%%%%%%%%%%%%%%%%%%%%%%%%%%%%%%%%%%%%%%%%%
%%%%%%%%%%%%%%%%%%%%%%%%%%%%%%%%%%%%%%%%%%%%%%%%%%%%%%%%%%%%%%%%%%%%%%%%%

%                C U B I C    &    T A N G E N T

\subsection{The cuspidal cubic with one tangent\label{sec:CuCuTan}}

\begin{figure}[ht!]
\begin{center}
\includegraphics[width=.3\textwidth]{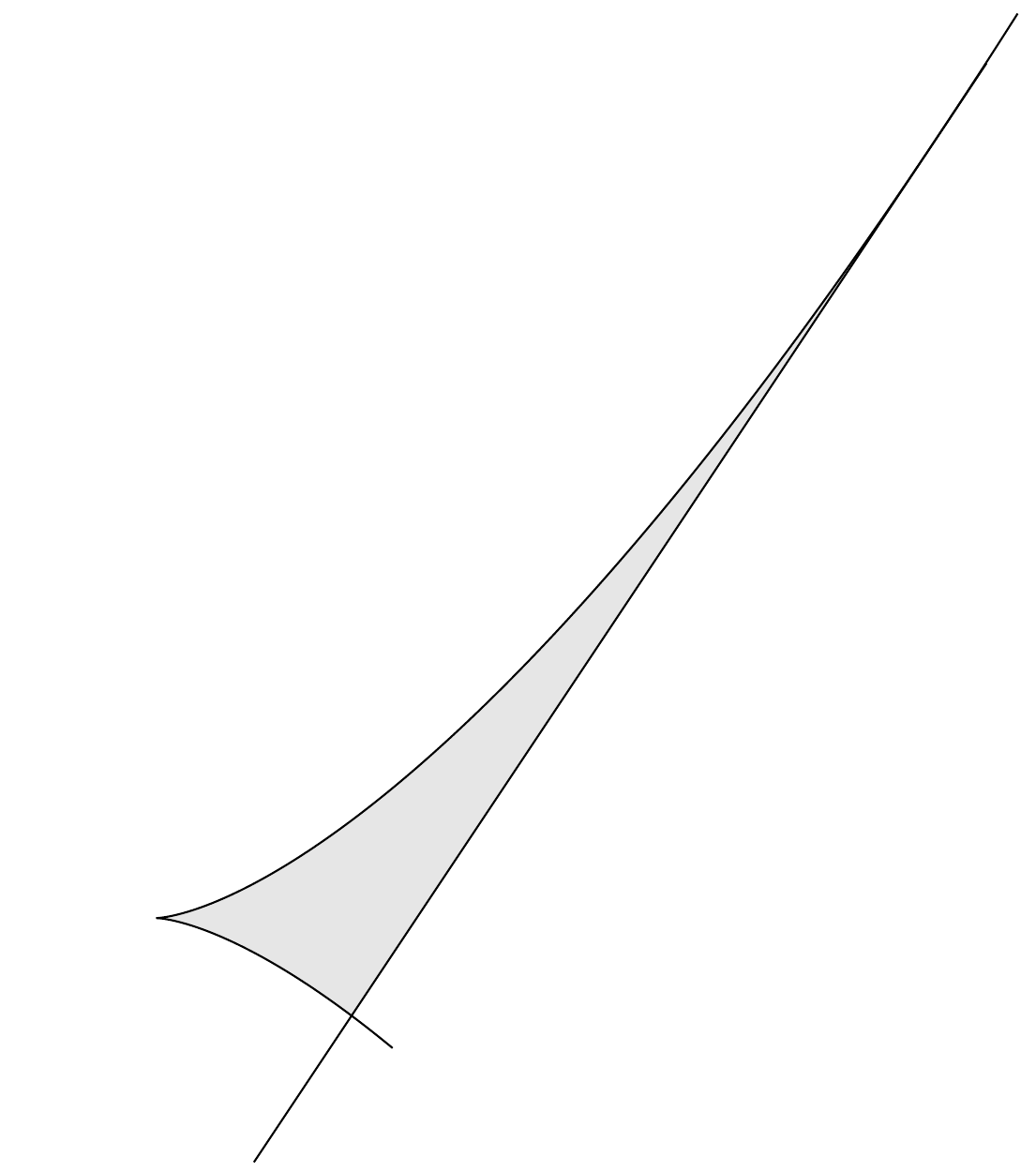}
\end{center}
\end{figure}
We may choose the boundary equation to be 
$$
   (X^3-Y^2)(2Y-3X+1))=0
$$
(one can write the second factor in the form $2(Y-1)-3(X-1)\,$).
Then, up to scaling, there is a unique solution
$$
  G= 2\bpm  4(X+Y-2X^2)   &  6(Y-2XY+X^2)   \\
            6(Y-2XY+X^2)  &  9(X-Y)(X+2Y)   \epm.
$$

The boundary having degree $4$, the density measure $\det(G)^{-1/2}$ belongs to the admissible solutions. Therefore, the Laplace  operator associated with this (co)-metric is an admissible solution. The scalar curvature is $2$, and therefore we may realize this Laplace operator as an image of the spherical Laplacian $\De_{\bS^2}$. 

The general density measure is $\rho= (X^3-Y^2)^{p-1}(2Y-3X+1)^{q-1}$,
$p>1/6$, $q>0$, $p+q>1/2$. For this measure we have
$$
\begin{aligned}
   &\LL_{p,q}(X)= \; -8(6p+3q-2)X + 4(6p-1), \\
   &\LL_{p,q}(Y)= -12(6p+3q-1)Y + 6(6p+1)X.
\end{aligned}
$$
In the case $p=q=1/2$, which corresponds to the Laplace operator, one may see
 that the operator is the image of a two-dimensional sphere, where $X$ is
 a degree $4$ polynomial and $Y$ has degree $6$.  
Indeed, it is worth to represent $X$ and $Y$  as
\beq
    \label{eq.CuCuTan}
    X=-\tfrac{1}{3}( t_1t_2+ t_2t_3+t_3t_1) \qquad\text{and}\qquad
    Y= \tfrac{1}{2}t_1t_2t_3
\eeq
 with  $t_1+t_2+t_3=0$,
 which reflects the fact that $X^3-Y^2$ (up to scaling) is the discriminant
 of the polynomial $T^3-3XT+2Y$.  

A solution is given by $t_i= 3x_i^2-1$, and one may check that all
 the relations concerning $\LL_{1/2,\,1/2}(X)$, $\LL_{1/2,\,1/2}(Y)$,
 $\Ga(X,X)$, $\Ga(X,Y)$ and $\Ga(Y,Y)$ are satisfied for this choice
(on the $2$-sphere $x_1^2+x_2^2+x_3^2=1$).

From this representation, it is clear that $X$ and $Y$ are invariant  under the symmetries through the hyperplanes $\{x_i=0\}$ and $\{x_i=x_j\}$. The fundamental domain for those reflexions is a triangle on the sphere, defined by the hyperplane coordinates, cut along its three medians, with angles $\pi/2, \pi/3, \pi/4$. This corresponds to one double point, one cusp and one tangency point.

In the case when $p=1/2$ and $q$ is a positive half-integer, 
 we may take   a unit sphere  in $\bR^n$, $n=6q$, whose elements we represent by
triples $(\xx_1,\xx_2,\xx_3)$ with $\xx_j\in\bR^{2q}$.
and consider $X$ and $Y$ given by \eqref{eq.CuCuTan} but
with $t_j=3\|\xx_j\|^2-1$, $j=1,2,3$.
Then, for the spherical Laplace operator on $\bS^{n-1}$ one can check that
$\Ga(X,X)$, $\Ga(X,Y)$ and $\Ga(Y,Y)$, $\LL_{p,q}(X)$ and $\LL_{p,q}(Y)$ satisfy the required equations.  It is certainly worth to mention that this model may also
be seen as the image of the triangle model (Section~\ref{triangle}) on the
triangle $\{(s_1,s_2,s_3)\in\bR^3\mid s_1+s_2+s_3=1, ~s_i\geq 0\}$
through the transformation $X= s_1s_2+s_2s_3+s_3s_1, ~Y= s_1s_2s_3$. 

For $p=1$ and a positive half-integer $q$, one may consider the following model.
For $(\xx_1,\xx_2,\xx_3)\in\bR^n$, $n=6+6q$,
$\xx_j\in\bR^{2+2q}$, we consider the $3\times 3$ symmetric matrix
$M_{ij}= (\xx_i\cdot\xx_j)-\frac13 \delta_{ij}$.
Then the restriction of this matrix to the
unit sphere in $\bR^n$ has trace $0$, and one considers its characteristic polynomial
$P(\lambda)= \det(\lambda\,\Id-M)$. Write
$P(\lambda)= \lambda^3-\frac13\lambda X -\frac2{27} Y$.
Then the image of the operator $\De_{\bS^{n-1}}$ is $\LL_{p,q}$.

This is case $\cd=1$ of a construction which works for $\cd=\dim_{\bR}\bK$
where $\bK$ is $\bR$, $\bC$, or $\bH$. In each of these three cases we
consider $(\xx_1,\xx_2,\xx_3)$ with $\xx_j\in\bK^{m}$, $m\ge 3$, and define
$M$, $X$, and $Y$ as above but with $\xx\cdot\yy$ understood as the Hermitiant
product $\sum_{j=1}^m x_i\bar y_i$ (see also Remark~\ref{rmq.DetH}).
Then the projection of $\De_{\bS^{n-1}}$, $n=3\cd m$, yields
$\LL_{p,q}$ with $2p=\cd+1$ and $2q =\cd(m-2)$ (recall that $\cd\in\{1,2,4\}$ and
$m\ge 3$).
It could be interesting to generalize this construction for the octonions.
A computation shows that literally the same formulas do not lead to the desired result.

 \par

%%%%%%%%%%%%%%%%%%%%%%%%%%%%%%%%%%%%%%%%%%%%%%%%%%%%%%%%%%%%%%
%%%%%%%%%%%%%%%%%%%%%%%%%%%%%%%%%%%%%%%%%%%%%%%%%%%%%%%%%%%%%%
%%%%%%%%%%%%%%%%%%%%%%%%%%%%%%%%%%%%%%%%%%%%%%%%%%%%%%%%%%%%%%
%%%%%%%%%%%%%%%%%%%%%%%%%%%%%%%%%%%%%%%%%%%%%%%%%%%%%%%%%%%%%%

%          S W A L L O W    T A I L

\subsection{ The swallow tail}\label{sect.swallow.tail}

\begin{figure}[ht!]
\begin{center}
\includegraphics[width=.3\textwidth]{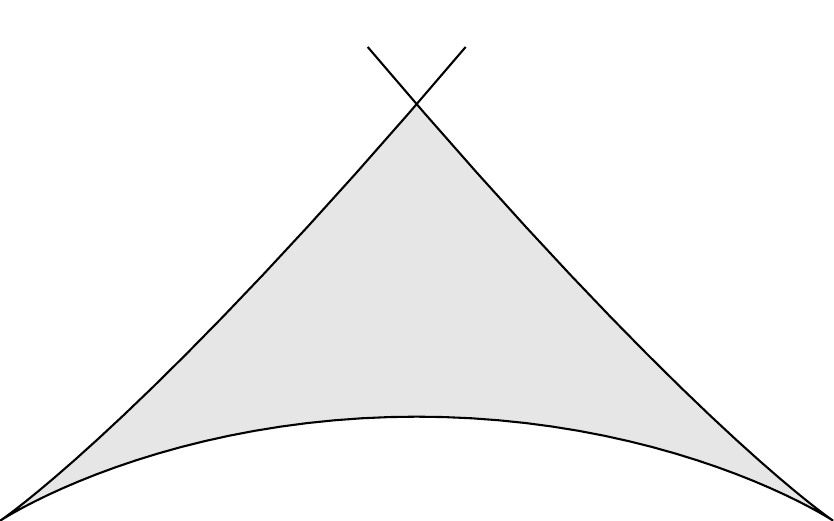}
\end{center}
\end{figure}
This is a degree $4$ algebraic curve, whose, up to affine transformations, we may chose the equation to be 
$$4\,X^{2}-27\,X^{4}+16\,Y-128\,Y^{2}-144
\,X^{2}Y+256\,Y^{3}=0.$$

This is the discriminant in $T$ of the polynomial $T^{4}-T^{2}+XT+Y$. Once again, the metric is unique up to scaling,  and we have
$$
   G= \bpm  2-8Y-9X^2  &  -X(12Y+1)                \\
            -X(12Y+1)  &  \frac{3}{2}X^2-16Y^2+4Y  \epm.
$$
The boundary having degree $4$, the measure density $\det(G)^{-1/2}$ is an admissible
solution, and for this measure, the corresponding Laplace operator has  constant
scalar curvature $2$, and therefore the operator may be represented on the unit
sphere $\bS^2$. 

The general measure density is $\rho= \det(G)^{p-1}$, $p>1/6$. For it we have
\[
    \LL_p(X)=  -6(6p-1)X, \qquad
    \LL_p(Y)= -4(12p-1)Y + 4p-1.
\]
In the Laplace-Beltrami case
$$
   \LL_{1/2}(X)= -12X,\qquad
   \LL_{1/2}Y= 1-20Y,
$$
which corresponds for $X$ to be an eigenvector of degree $3$
and $Y-1$ to be an eigenvector of degree $4$.

Taking in account that the boundary is a discriminant, we should look for 
$$
   -X= t_1t_2t_3+ t_2t_3t_4+ t_3t_4t_1+ t_4t_1t_2, \qquad
    Y= t_1t_2t_3t_4,
$$
on the variety given by
$$
   t_1+t_2+t_3+t_4= 0, \qquad\sum_{i< j} t_it_j=-1.
$$
Since $\sum t_i^2 = \big(\sum t_i\big)^2-2\sum t_i t_j$,
this variety is the intersection of a sphere $\bS^{3}$ of radius $\sqrt{2}$ with the
hyperplane $\{\sum t_i= 0\}$, which is again a sphere with radius $\sqrt{2}$.

To compute the image of $\De_\Sigma$ through $(X,Y)$, we introduce the following
orthogonal coordinates on the plane $\sum t_i=0$:
$$
\begin{matrix}
   t_1=( x_1+x_2+x_3)/\sqrt2\,,\quad &
   t_3=(-x_1+x_2-x_3)/\sqrt2\,,\\
   t_2=( x_1-x_2-x_3)/\sqrt2\,,\quad &
   t_4=(-x_1-x_2+x_3)/\sqrt2\,
\end{matrix}
$$
where the scaling factor $\sqrt2$ is chosen so that the unit 2-sphere $\bS^2$ in
$\bR^3$ with coordinates $\xx=(x_1,x_2,x_3)$ is mapped onto $\Sigma$.
In these coordinates we have
$X = 2\sqrt2\,x_1x_2x_3$ and $Y=\frac12(x_1^4+x_2^4+x_3^4)-\frac14\|\xx\|^4$, thus
the restriction of $Y$ to $\bS^2$ is $\frac12(x_1^4+x_2^4+x_3^4)-\frac14$.
Then the image of $\De_{\bS^2}$ through $(X,Y)$ can be easily computed
using~\eqref{lapl.sph.homog} and the result will be $\LL_{1/2}$.

For other values of $p$, the  discriminant form of the boundary suggests that
one looks at Hermitian $4\times 4$ matrices $M$  with vanishing trace,
restricted on the sphere
$\tr(M^*M)= 2$, embedded with the induced spherical structure, in the real and
complex cases,
and look at the induced operator on the characteristic polynomial%
\footnote{
   The coefficient of $\lambda^{n-2}$ in $\det(\lambda I-M)$ is equal to
   $\frac12(\tr\,M)^2 - \frac12\tr(M^tM)$
   for any Hermitian $n\times n$ matrix $M$,
}
$P(\lambda) = \lambda^4 - \lambda^2 + \lambda X + Y$.
Then the image of $\Delta_{S^{2+6\cd}}$ ($\cd=1$ for $\bR$ and $\cd=2$ for $\bC$)
through $(X,Y)$ is $\LL_p$ with $2p=\cd+1$.
  The quaternionic case is left to the reader as an exercise
(\cite{bakry-zani} can be used as a hint).

Observe also that the mapping $(X,Y)\mapsto(X^2,Y)$ transforms all these models
into some models for
the cuspidal cubic with tangent (indeed, the spherical triangle with angles
$(\frac\pi2,\frac\pi3,\frac\pi3)$ folded along its axis of symmetry yields
a triangle with the angles $(\frac\pi2,\frac\pi3,\frac\pi4)$).

%%%%%%%%%%%%%%%%%%%%%%%%%%%%%%%%%%%%%%%%%%%%%%%%%%%%%%%%%%%%
%%%%%%%%%%%%%%%%%%%%%%%%%%%%%%%%%%%%%%%%%%%%%%%%%%%%%%%%%%%%
%%%%%%%%%%%%%%%%%%%%%%%%%%%%%%%%%%%%%%%%%%%%%%%%%%%%%%%%%%%%
%%%%%%%%%%%%%%%%%%%%%%%%%%%%%%%%%%%%%%%%%%%%%%%%%%%%%%%%%%%%

%                D E L T O I D

\subsection{ The deltoid\label{sec.A2}}
\begin{figure}[ht!]
\begin{center}
\includegraphics[width=.3\textwidth]{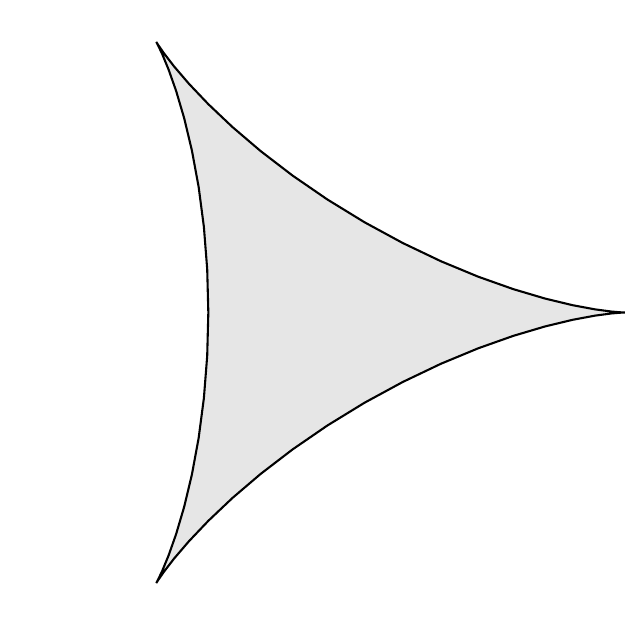}
\end{center}
\end{figure}
  
  In this case, up to affine transformation, we may choose the boundary equation to be
  $$
      (X^{2}+Y^{2})^{2}+ 18(X^{2}+Y^{2})- 8X^{3}+24XY^{2}-27=0.
  $$
  There is a unique metric $g$ up to scaling, which is
  $$
      G= \bpm 9+6X+Y^2-3X^2& -2Y(2X+3)\\ -2Y(2X+3)& 9-6X+X^2-3Y^2\epm.
  $$
  For the density  measure $\rho= \det(G)^{p-1}$, $p>1/6$, we have
  $$
     \LL_p(X)= -2(6p-1)X,\qquad
     \LL_p(Y)= -2(6p-1)Y.
  $$ 
  The operator looks simpler in complex variables: setting $Z= X+ iY$, one gets
  $$
     \Ga(Z,Z)= 12\bar Z-4Z^2,\qquad
     \Ga(Z, \bar Z)= 18-2Z\bar Z,\qquad
     \Ga(\bar Z, \bar Z)= 12 Z-4\bar Z^2,
  $$ 
  and 
  $$
     \LL_p(Z)= -2(6p-1)Z,\qquad
     \LL_p(\bar Z)= -2(6p-1)\bar Z.
  $$
  Under this form, it is easier to check the eigenvalues for $\LL_p$,
 since the action of the operator on the highest degree term of a polynomial
 is diagonal, and we see that, for any degree $p$, the highest degree part
 of any eigenvector is a monomial, say $Z^q\bar Z^r$, and the same use of complex variables gives that the eigenvalue
 corresponding to a polynomial whose highest degree term is $Z^p\bar Z^{q}$ is 
 $-3(q+r)(q+r+4p+2)-(q-r)^2$.
  
     Once again, as it is the case whenever the boundary is degree $4$,
 the density measure $\det(G)^{-1/2}$ is an admissible solution, and this
 corresponds to a Laplace operator for a metric which has zero scalar curvature.
 We may represent this operator from a Euclidean Laplacian in dimension 2.
 For this choice of the measure, one has $\LL_{1/2}(Z)= -4Z$, and if we 
 identify $\bR^2$ with the complex plane $\bC$, one may represent this using
 $$
     Z= e^{2i (1\cdot z)}+ e^{2i(\omega\cdot z)}+ e^{2i(\bar\omega\cdot z)},
 $$
 where $z= x+iy\in \bC$, $1$, $\omega$, $\bar\omega$ are the three third roots
 of unity
 (solution of $z^3=1$) and $z_1\cdot z_2$ is the Euclidean scalar product
 $\Re(z_1\bar z_2)$. 
 One may directly check the $\LL_{1/2}$ is the image of $\De_{\bR^2}$
 through $Z$. (The interior of the deltoid is indeed the image of $\bR^2$
 through $Z$.)   
 Moreover, the function $Z$ is invariant under the symmetries with respect
 to the lines 
 $$
    D_1= \{\Im(z)=0\}, \qquad
    D_2= \{te^{i\pi/3}, t\in \bR\},\qquad
    D_3= \{ae^{i\pi/6}+ te^{2i\pi/3}\},
 $$
 with $a= \pi/\sqrt{3}$. Those three lines determine a equilateral triangle
 $(ABC)$ in the plane, and any function which has those symmetries is also
 invariant under all the symmetries with respect to the lines of the triangular
 lattice generated by $A,B,C$ (that is all the lines parallel to
 $D_1$, $D_2$, $D_3$ which are distant from $ka$, $k\in \bN$).
 (This group of symmetries is the affine Weyl group associated with the root system
 $A_2$.

  The deltoid is then the image of the boundary of the triangle $(ABC)$ through
 $Z$, and it is not hard to see that the restriction of $Z$ to $(ABC)$ is
 injective. Then, the functions of the form  $F(Z)$ are nothing else than the
 functions which are invariant under the symmetries of the triangular lattice,
 and $\LL_{1/2}$ is just $\De$ acting on functions invariant under
 the affine Weyl group associated with $A_2$.

  As usual, this model extends to rank 2 symmetric spaces with restricted
 root system $A_2$. For a study of this case, see also
 Koornwinder~\cite{Koorn3,Koorn4}. Indeed, it is perhaps worth to notice
 that the boundary equation is the discriminant of the polynomial
 $T^3-ZT^2+ \bar Z T + 1$, putting forward the interest of representing
 $Z$ as $\lambda_1+\lambda_2+ \lambda_3$ where $|\lambda_i|=1$ and
 $\lambda_1\lambda_2\lambda_3=1$. In particular, one may consider
 the Casimir operator on $SU(3)$. Indeed, if $Z$ denotes the trace of
 a $SU(3)$ matrix, due to the fact that eigenvalues
 $(\lambda_1, \lambda_2, \lambda_3)$ satisfy $|\lambda_i|=1$,
 $\lambda_1\lambda_2\lambda_3=1$, one sees that the characteristic polynomial
 $P(\lambda)$ writes
 $P(\lambda)= -\lambda^3+ Z\lambda^2 -\lambda \bar Z +1$.
 Then, using formulae~\eqref{lapl.su(n)1} and \eqref{lapl.su(n)2}, one sees that
 $$
   \De_{SU(3)} Z = - \tfrac{16}{3}Z,\qquad
   \De_{SU(3)} \bar Z=- \tfrac{16}{3}\bar Z,
 $$
 and 
 \beqnas
    &\Ga_{SU(3)}(Z,Z)= \tfrac43(3\bar Z - Z^2),
    & ~\Ga_{SU(3)}(Z,Z)=  \tfrac43(3Z-\bar Z^2), \\
    &\Ga_{SU(3)}(Z,\bar Z))= \tfrac23({Z\bar Z }-9),&
 \eeqnas
 which shows that $\LL_{3/2}$ is the image of $3\De_{SU(3)}$ through $Z$.

 It is worth to observe that this operator preserves the functions which are
 symmetric under the conjugacy $Z\mapsto \bar Z$. Indeed, choosing as new
 variables $X= Z+\bar Z$ and $Y= Z\bar Z$, one gets for the cometric $\Ga$
 in those variables
 $$
   \bpm         1+X-X^2+Y            & \frac12 X+X^2-2Y-\frac32 XY  \\
        \frac12 X+X^2-2Y-\frac32 XY  &      Y-3XY-3Y^2+X^3          \epm
 $$
 with 
 $$
   \LL X= -\lambda X, \qquad
   \LL (Y)= 1-(2\lambda+1)Y.
 $$

 The determinant of this matrix is, up to some constant,
 $(4X-Y^2)(4X^3-3Y^2-12XY-6Y+1)$, which corresponds to a domain delimited by
 a parabola and a cuspidal cubic which have a double tangent at their
 intersection points (and is a degree 5 curve). This domain is the image
 of the deltoid domain through the map $(X,Y)$. We may now describe
 a two parameters family of orthogonal polynomials  associated with this domain,
 but with a weighted degree  $\deg(X^p Y^q)= 2p+q$. For the associated Laplace operator,
 this corresponds to the 2 dimensional Euclidean Laplacian  associated with
 the symmetries of the root system $G_2$.

%%%%%%%%%%%%%%%%%%%%%%%%%%%%%%%%%%%%%%%%%%%%%%%%%%%%%%%%%%%%%%%%%%%%
%%%%%%%%%%%%%%%%%%%%%%%%%%%%%%%%%%%%%%%%%%%%%%%%%%%%%%%%%%%%%%%%%%%%
%%%%%%%%%%%%%%%%%%%%%%%%%%%%%%%%%%%%%%%%%%%%%%%%%%%%%%%%%%%%%%%%%%%%
%%%%%%%%%%%%%%%%%%%%%%%%%%%%%%%%%%%%%%%%%%%%%%%%%%%%%%%%%%%%%%%%%%%%
%%%%%%%%%%%%%%%%%%%%%%%%%%%%%%%%%%%%%%%%%%%%%%%%%%%%%%%%%%%%%%%%%%%%

%            N O    B O U N D A R Y     ( F U L L   R^2 )

\section{The full $\bR^2$ case\label{no.boundary}}

In this Section, we consider the case where $\Omega= \bR^2$, and concentrate on the SDOP problem. From the one dimensional models and the tensorization procedure, we already know that Gaussian measures provide such orthogonal polynomials, with the two-dimensional Ornstein-Uhlenbeck operator as associated diffusion operator.
We shall prove in this section the following

\bthm If $(\bR^2, G, \rho dx)$ is a solution of the SDOP problem, then, up to an affine change of coordinates, one has $\rho= \frac{1}{2\pi}\exp(-\frac12(X^2+ Y^2))$, and 
$$G= \bpm a+ \mu X^2& -\mu XY \\ -\mu XY& c+ \mu Y^2\epm,$$ where  $a$ and $c$ are positive real numbers. The corresponding operator $\LL$ is a sum of the general Ornstein-Uhlenbeck operator
$$
    L_{a,c} = a(\partial^2_X-X\partial_X) + c(\partial^2_Y - Y\partial_Y)
$$
and a squared rotation $\mu (Y\partial_X-X\partial_Y)^2$.

\ethm

In this situation,  we do not have boundary equation to restrict the analysis of the (co-)metric $(G^{ij})$.  We therefore  look for $3$ polynomials $G^{11},G^{12},G^{22}$ of degree at most $2$ in the variables $(X,Y)$ with $\Delta=G^{11}G^{22}-(G^{12})^2>0$ on $\mathbb R^2$ and for a function $h$ such that  for the measure $d\rho=e^hdXdY$, the polynomials of 
 $(X,Y)$ are dense in $\cL^2(\rho,\mathbb R^2)$. From \eqref{mesure}, there exists $L_X$ and $L_Y$ affine forms in $\mathbb R^2$ such that
\begin{equation}
   \label{relfond}
   \partial_X h=\frac 
   {G^{22}L_X-G^{12}L_Y}{\Delta},
  ~\partial_Y h = \frac{-G^{12}L_X+ G^{11}L_Y}{\Delta}
\end{equation}
 Let us first show that $\Delta$ has degree at most $2$. If $\Delta$ is of degree $4$, and since $\Delta>0$ on $\mathbb R^2$,  there is at least a cone in which, for some constant $c$,  $\Delta\geq c (X^2+Y^2)^2$ at infinity and $d\rho$ cannot integrate any polynomial.
Hence $\Delta$ is of degree at most $3$ and since it is positive on $\mathbb R^2$, $\Delta$ is of an even degree, thus of degree $2$ or zero.

%%%%%%%%%%%%%%%%%%%%%%%%%%%%%%%%%%%%%%%%%%%%%%%%%%%%%%%%%%%%%%%%%%%%%%%%
%%%%%%%%%%%%%%%%%%%%%%%%%%%%%%%%%%%%%%%%%%%%%%%%%%%%%%%%%%%%%%%%%%%%%%%%
%%%%%%%%%%%%%%%%%%%%%%%%%%%%%%%%%%%%%%%%%%%%%%%%%%%%%%%%%%%%%%%%%%%%%%%%

%           NO BOUNDARY    Deg(Delta) = 2

  \subsection{Case where  $\Delta$ is degree  2}
  
  Let us first consider the case where $\Delta$ is complex irreducible. From the form of the measure (Proposition~\ref{prop.measure}), we know that 
\beq
     \label{eq.rho.R2}
     h=\log\rho= -P-\alpha  \log\Delta,
\eeq
where $P$ is a polynomial with degree at most $2$. The terms of highest degrees
in $P$ is a positive definite quadratic form. Indeed, otherwise the measure
$\rho\, dx$ would not integrate the polynomials. Thus $\deg P=2$.

Let us show that $\alpha=0$. Indeed, suppose that $\alpha\ne 0$.
Then by Proposition~\ref{rho=infty} $(G^{11},G^{12},G^{22},\Delta)$
is a solution of the $\bR$-AlgDOP problem (see Definition~\ref{gal.pb.compact}).
Since $\deg\Delta=2$ and $\Delta\ne0$ on $\bR^2$, up to affine linear transformation,
we have either $\Delta=X^2+Y^2+1$ or $\Delta=X^2+1$.
Then Proposition~\ref{AlgDOP.circle} and Lemma \ref{lem.ParLines} imply that
$$
      G = \bpm 1+X^2 & XY \\ XY & 1+Y^2 \epm
      \qquad\text{or}\qquad
      G = \bpm 1+X^2 & 0 \\ 0 & 1 \epm
$$
(we do the change of coordinates $X=ix$, $Y=iy$ in the corresponding solution
in Proposition~\ref{AlgDOP.circle} and Lemma \ref{lem.ParLines}).
Thus condition~\eqref{eq.rho.1} takes the form
$$
    (1+X^2)\,\partial_X P +    XY  \,\partial_Y P = L_X,      \qquad
       XY  \,\partial_X P + (1+Y^2)\,\partial_Y P = L_Y,
$$
or, respectively,
$$
    (1+X^2)\,\partial_X P = L_X,      \qquad
           \,\partial_Y P = L_Y,
$$
with $\deg L_X\le 1$ and $\deg L_Y\le 1$. One easily checks that in both
cases this is
possible only when $P$ is a constant which contradicts the condition $\deg P=2$.
Thus $\alpha=0$.

Recall that the highest order homogeneous component of $P$ is
a positive definite quadratic form. Hence by an affine change
of coordinates we may reduce to the case $P=(X^2+Y^2)/2$.
Then condition~\eqref{eq.rho.1} reads
$$
       G^{11}X + G^{12}Y = L_X,      \qquad
       G^{21}X + G^{22}Y = L_Y,      \qquad
$$
with $\deg L_X\le 1$ and $\deg L_Y\le 1$.
This is a linear system on the coefficients of $G^{ij}$.
By solving it, we obtain, up to a scalar factor,
$$
   G= \bpm a+ Y^2& b-XY\\b-XY& c+ X^2\epm,
$$
with some constant positive definite matrix $\bpm a& b\\b&c\epm$
and the associated operator is the sum of the general Ornstein-Uhlenbeck operator
$$
   L_{a,b,c} = a\partial^2_X + 2b\,\partial_X\partial_Y + c\partial^2_Y
             - (aX+bY)\partial_X - (bX+cY)\partial_Y
$$
and a plane squared rotation $(Y\partial_X-X\partial_Y )^2$.

With a further rotation, one may reduce to the case where $b=0$.

%%%%%%%%%%%%%%%%%%%%%%%%%%%%%%%%%%%%%%%%%%%%%%%%%%%%%%%%%%%%%%%
%%%%%%%%%%%%%%%%%%%%%%%%%%%%%%%%%%%%%%%%%%%%%%%%%%%%%%%%%%%%%%%
%%%%%%%%%%%%%%%%%%%%%%%%%%%%%%%%%%%%%%%%%%%%%%%%%%%%%%%%%%%%%%%
%%%%%%%%%%%%%%%%%%%%%%%%%%%%%%%%%%%%%%%%%%%%%%%%%%%%%%%%%%%%%%%
%%%%%%%%%%%%%%%%%%%%%%%%%%%%%%%%%%%%%%%%%%%%%%%%%%%%%%%%%%%%%%%

 \subsection{ $\Delta$ is constant}

We can boil down to  $\Delta=1$.  Our aim is to prove in this section
that $G^{ij}$ are constant -- in this case \eqref{relfond} implies that
$P$ is of degree 2 which corresponds to Ornstein-Uhlenbeck operators.
Suppose that $G$ is non-constant. Then one of $G^{11}$, $G^{22}$ is
non-constant. Let it be $G^{11}$.
  
  Since $\Delta=1$, we have $G^{11}G^{22}= (G^{12}-i)(G^{12}+i)$.
  If $G^{11}$ is irreducible,
  then $G^{11}= \lambda(G^{12} \pm i)$, which is impossible since $G^{11}$ is real.
  Therefore,  $G^{11}=l_1 l_2$ with $\deg l_1=\deg l_2=1$.
  Since $G^{11}>0$
  on $\bR^2$, the only solution is $G^{11}=(l_a+\alpha)(l_a+\bar\alpha)$
  where $l_a$ is a nonzero real linear form
  and $\alpha$ is a non-real complex number.
  Similarly, $G^{22}=(l_c+\gamma)(l_c+\bar\gamma)$,
  where $l_c$ is a real linear form,
  and if $l_c$ is nonzero, then $\gamma$ is a non-real complex number.
Hence, changing if necessary the sign of $l_c+\gamma$, we have
  $$
      G^{12}+i = (l_a+    \alpha)(l_c+    \gamma), \qquad
      G^{12}-i = (l_a+\bar\alpha)(l_c+\bar\gamma).
  $$
  
  We know that $(l_a+\alpha)(l_c+\gamma)\pm i$ is a real polynomial. Hence
  all its homogeneous forms are real, in particular, its linear form is real, i.e., 
  $\alpha l_c+ \gamma l_a = \bar \alpha  l_c+ \bar \gamma l_a$ whence
  $$
      (\alpha-\bar\alpha) l_c = (\bar\gamma-\gamma)l_a.
  $$
  Since $\alpha$ is non-real, we obtain $l_c= \nu l_a$
  for some real number $\nu$. From now on we denote $l_a$ just by $l$.
  So, we have  
$$
     \bpm G^{11} & G^{12}\\
          G^{12}& G^{22}\epm
   = \bpm l^2+ (\alpha+ \bar \alpha)l + \alpha\bar \alpha & \nu l^2+c\,l+b\\
    \nu l^2+c\,l+b & \nu^2l^2+\nu (\gamma+ \bar \gamma)l + \gamma\bar \gamma\epm,
$$
where $l$ is some real non zero linear form,
$\nu$, $c= (\alpha\nu + \gamma)$, $b= \alpha\gamma-i$ are real numbers whereas
$\alpha$ and $\gamma$ some non-real complex numbers.
  
  Assume first  that $\nu \neq 0$. Then
  one may reduce to the case $c=0$ through a translation,
  and eventually to the case that $\nu=\alpha\bar\alpha=\gamma\bar\gamma =1$, $b=0$,
  through a linear change of coordinates. In this case, the determinant
  is easy to compute and the only solution, up to a change of $Y$ into $-Y$
  and an exchange of $X$ and $Y$, is    
  $$
     \bpm G^{11} & G^{12}\\
          G^{12} & G^{22}
     \epm
    =\bpm l^2+\sqrt2\,l+1  &        l^2     \\
              l^2          &   l^2-\sqrt2\,l+1  \epm.
  $$  
  Let us show that then the measure $\rho\, dx$ cannot be integrable.

  We know from Proposition~\ref{prop.measure} that $\rho= \exp (P) dx$, where $P$
  is some polynomial of degree at most $4$. By \eqref{mesure} we have
  $$
    \bpm G^{11} & G^{12}\\ G^{12} & G^{22} \epm
    \bpm \partial_X P\\\partial_Y P\epm
    = \bpm l_X+c_X\\ l_Y+c_Y\epm,
  $$
  where $l_X$ and $l_Y$ are linear forms, and $c_X$ and $c_Y$ are constants.
  
  From this, we get 
  $$
    \bpm \partial_X P \\ \partial_Y P\epm= \bpm G^{22} & -G^{12} \\ 
                                               -G^{12} & G^{11} \epm
    \bpm l_X+c_X\\ l_Y+c_Y\epm.
  $$
  
  Writing $P= P_4+P_3+P_2+P_1+P_0$ where $P_k$ is homogeneous of degree $k$,
  one sees from this equation that $(\partial_X+\partial_Y)P_4=0$.
  Since the measure has to be integrable, this also requires that
  $(\partial_X+\partial_Y) P_3= 0$. But, looking precisely at
  $(\partial_X+\partial_Y) P_3$, we see that it is equal to
  $-\sqrt2\,l(l_X-l_Y)$, thus  $l_X-l_Y=0$.
  From this, one sees that $\partial_X P_4= \partial_Y P_4=0$,
  and hence $P_4=0$. This implies that $P_3=0$ too (since once again
  the measure has to be finite). We have
$$
    \partial_X P_3 = (c_X-c_Y)l^2 - \sqrt2\,l\,l_X, \qquad
   -\partial_Y P_3 = (c_X-c_Y)l^2 - \sqrt2\,l\,l_Y,
$$
thus $P_3=0$ implies $l_X=l_Y = \frac{\sqrt2}2(c_X-c_Y)l$.
Since $\partial_X P_2=l_X-\sqrt2\,c_Xl$ and
$\partial_Y P_2=l_Y+\sqrt2\,c_Yl$, we obtain
$\partial_X P_2 = -\partial_Y P = -\frac{\sqrt2}2(c_X+C_Y)l$ whence
$(\partial_X+\partial_Y)P_2=0$.
Therefore $P_2$ is a degenerate quadratic form,
and $\exp(P)$ is non integrable on $\bR^2$.

  In the case where $\nu=0$ with the same argument, we boil down to
  $$
     \bpm G^{11} & G^{12}\\ G^{12} & G^{22} \epm= \bpm l^2+1& l\;\\ l& 1\;\epm.
  $$
  The same reasoning shows that the measure may not be integrable on $\bR^2$.

 %%%%%%%%%%%%%%%%%%%%%%%%%%%%%%%%%%%%%%%%%%%%%%%%%%%%%%%%%%%%%%%%%%%%
 %%%%%%%%%%%%%%%%%%%%%%%%%%%%%%%%%%%%%%%%%%%%%%%%%%%%%%%%%%%%%%%%%%%%
 %%%%%%%%%%%%%%%%%%%%%%%%%%%%%%%%%%%%%%%%%%%%%%%%%%%%%%%%%%%%%%%%%%%%
 %%%%%%%%%%%%%%%%%%%%%%%%%%%%%%%%%%%%%%%%%%%%%%%%%%%%%%%%%%%%%%%%%%%%
 %%%%%%%%%%%%%%%%%%%%%%%%%%%%%%%%%%%%%%%%%%%%%%%%%%%%%%%%%%%%%%%%%%%%
 %%%%%%%%%%%%%%%%%%%%%%%%%%%%%%%%%%%%%%%%%%%%%%%%%%%%%%%%%%%%%%%%%%%%
 %%%%%%%%%%%%%%%%%%%%%%%%%%%%%%%%%%%%%%%%%%%%%%%%%%%%%%%%%%%%%%%%%%%%
 %%%%%%%%%%%%%%%%%%%%%%%%%%%%%%%%%%%%%%%%%%%%%%%%%%%%%%%%%%%%%%%%%%%%

 %       N O N    C O M P A C T    W I T H    B O U N D A R Y

\section{ Non compact cases with boundaries\label{non.compact.with.boundary}.}

\begin{figure}[ht!]
\begin{center}
\includegraphics[width=.9\textwidth]{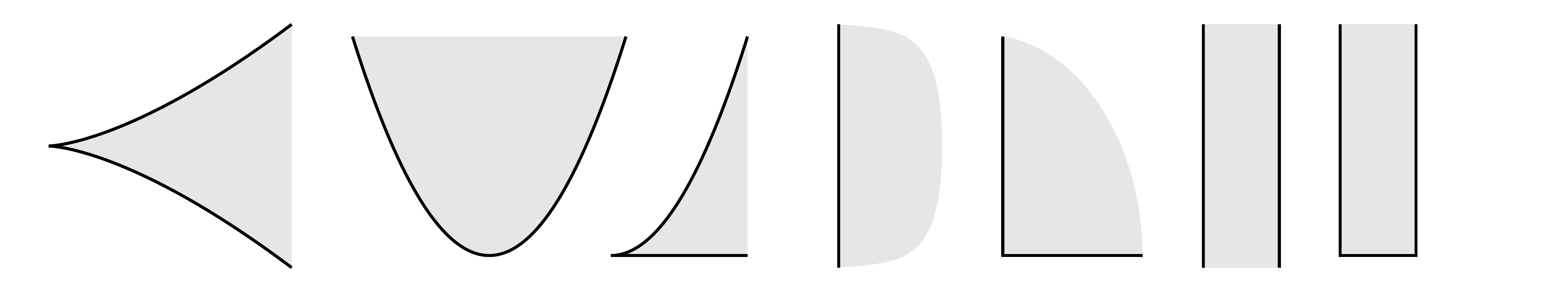}

\hskip4mm
(1)\hskip18mm
(2)\hskip13mm
(3)\hskip10mm
(4)\hskip11mm
(5)\hskip9mm
(6)\hskip7mm
(7)

\if01{ %===============  for unbounded.eps ================
\hskip10mm
(1)\hskip16mm
(2)\hskip13mm
(3)\hskip12mm
(4)\hskip12mm
(5)\hskip12mm
(6)\hskip8mm
(7)
}\fi %========================================
\end{center}
\end{figure}

In this section, we again consider the SDOP problem, which is perhaps not enough to describe all the possible solutions of the general DOP problem (although we have no example of solution of the latter beyond the cases described here).  
\par We describe all the possible models, but we do not give any geometric interpretation, and do not detail for which values of the parameters appearing in the measure the polynomials are  dense in $\cL^2(\mu)$. However, in all the cases described below, it is indeed the case for at least some values of these parameters. Moreover, one could  give a geometric construction for many   models as images of  Ornstein-Uhlenbeck operators in some Euclidean space, associated with Gaussian measures.
\par
Following the results of Section~\ref{section2},  we reduce to the  cases where every factor $\Delta_p$ appearing in the boundary satisfies the fundamental equations \eqref{eq.thm.gal}.    We also need, for  the measure $d\mu= e^hdx$, that $\cL^2(\mu)$ contains every polynomial. Hence in any case, we have to look for the existence of  such a measure, which will turn out to be the main restriction. 
We indeed require more, namely that polynomials are dense in $\cL^2(\mu)$. 
We know from Proposition \ref{prop.measure} the general form of the measure. In addition to the boundary terms, there appear in $h$ a polynomial term $P$ which will be crucial when integrating the polynomials on the domain (see previous section) . This constraint will help us to restrict the number of cases for the metric $(G)$.
Moreover, if the determinant $\Delta$ of $(G)$ has no multiple factors and the domain contains an open cone, the degree of $\Delta$ is at most $3$.  When there are multiple factors, the same kind of analysis can be undergone.
\par
The algebraic analysis undergone in Section~\ref{dim2} still holds, and produces the following list of possible boundaries.

%%%%%%%%%%%%%%%%%%%%%%%%%%%%%%%%%%%%%%%%%%%%%%%%%%%%%%%%%%%

\smallskip\par
\benum[(\ncCu)]
\item $\partial \Omega= \{ Y^2-X^3=0\}$. 
In this case, the general metric is given by
$$
  G=\alpha\bpm 4X^2 & 6XY  \\ 6XY  & 9Y^2 \epm
   + \beta\bpm 4X   & 6Y   \\ 6Y   & 9X^2 \epm
   +\gamma\bpm 4Y   & 6X^2 \\ 6X^2 & 9XY  \epm
$$
Here, the determinant $\Delta$ is
$36(X^3 - Y^2)\big((\alpha\beta-\gamma^2)X - \alpha\gamma Y + \beta^2\big)$.
By Corollary~\ref{cor.gal.2} we have $\deg \Delta \leq 3$, hence $\alpha=\gamma=0$, and by homogeneity we may  restrict to  $\beta=1$.  Since $\Delta>0$ in the interior of the domain $\Omega$, one sees that the domain must be $X^3> Y^2$. This leads to measures on the domain $X^3>Y^2$
of the form
$$
   d\mu_{a,b}=C_{a,b} (X^3-Y^2)^{a-1}\exp(-bX)\, dXdY,\quad a>1/6,\;\; b>0.  . 
$$

%%%%%%%%%%%%%%%%%%%%%%%%%%%%%%%%%%%%%%%%%%%%%%%%%%%%%%%%%%%

\item[(\ncPa)] $\partial \Omega= \{ Y-X^2=0\}$. 
Then $\Delta= (Y-X^2)\Delta_1$ where, by Corollary~\ref{cor.gal.2},
either $\Delta_1$ is proportional to $(Y-X^2)$, or $\deg \Delta_1 \leq 1$.
\benum 
\item [(\ncPa i)] $\Delta= c (Y-X^2)^2$. 
In this case, the general metric is given in Proposition~\ref{AlgDOP.parab.2}(1).
We may show that there is no finite measure solution for the problem.

\item[(\ncPa ii)]  $\deg \Delta =3$ (cf.~Proposition~\ref{AlgDOP.parab.2}(2)).
 The metric for which there exist a measure solution for the problem  may be written as 
 \beq\label{eq.ncPa}
     G= \bpm  1+\alpha(Y-X^2) & 2 X\\
                   2 X        & 4 Y\epm, \qquad \alpha > 0.
 \eeq
 We have $\det G=4(Y-X^2)(1+\alpha Y)$, this is why we impose the condition $\alpha>0$.
 By a change of coordinates $X\mapsto cX, Y\mapsto c^2Y$, we may  reduce to $\alpha=1$.
 The existence of a finite measure solution imposes $\Omega= \{Y>X^2\}$ and
 measures to be
$$
   (Y-X^2)^{a-1}\exp\big(-bY\big)dXdY,\quad a,b>0.
$$

\item[(\ncPa iii)] $\Delta$ has degree $2$.
The only metric for which there exists a measure solution is
\eqref{eq.ncPa} with $\alpha=0$.
This is a limit case of the previous one.
Notice that any measure of the form $(Y-X^2)^{a-1}\exp\big(-bY+cY\big)dXdY$
with $a,b>0$ is admissible in this case. However we may reduce to $c=0$
by $(X,Y)\mapsto(X+q,Y+2qX+q^2)$ which is an isometry of $(\Omega,g)$.
\eenum

%%%%%%%%%%%%%%%%%%%%%%%%%%%%%%%%%%%%%%%%%%%%%%%%%%%%%%%%%%%

\item[(\ncPaTa)] $\Omega$ is bounded by $\{ Y-X^2=0\}$ and a line.
Up to an affine linear transformation, it is enough to consider
the following three cases for the line.

\benum
\item[(\ncPaTa i)] $\partial \Omega\subset \{ Y(Y-X^2)=0\}$.
Then the metric is \eqref{eq.parab.2} with $\beta=\gamma=\lambda=0$.
If $\Delta$ has a multiple factor, then $r=0$.
One can check that a finite admissible measure does not exist.

If there are no multiple factors of $\Delta$, then the integrability condition implies
$\deg\Delta=3$ (see Corollary~\ref{cor.gal.2}).
This occurs only for $\alpha=\mu=0$, hence
$$
    G = r\bpm X & 2Y\\ 2Y & 4XY\epm.
$$
This matrix is positive definite in $\Omega=\{rX>0$ and $0<Y<X^2\}$.
The admissible measures for $r>0$ are
$$
   C\,Y^{a-1}(X^2-Y)^{b-1}\,\exp(-cX),\quad a,b,c>0,\quad a+b>1/2.
$$

%%%%%%%%%%%%%%%%%%%%%%%%%%%%%%%%%%%%%%%%%%%%%%%%%%%%%%%%%%%

\item[(\ncPaTa ii)] $\partial \Omega\subset \{ X(Y-X^2)=0\}$.
The metric is \eqref{eq.parab.2} with $\lambda=0$,
$\alpha=-\mu$, $\beta=-2r$.
We have always $\deg\Delta=4$,
and a multiple factor $X^2$ appears
only if $\alpha=\mu=0$.
One can check that a finite admissible measure does not exist.

%%%%%%%%%%%%%%%%%%%%%%%%%%%%%%%%%%%%%%%%%%%%%%%%%%%%%%%%%%%

\item[(\ncPaTa iii)] $\partial \Omega\subset \{ Y(Y-X^2+\varepsilon)=0\}$,
$\varepsilon=\pm1$.
In the case $\varepsilon=1$,
the metric is $G'_\alpha$ computed in Section~\ref{parab1}.
We have $\deg\Delta=3$ only for $\alpha=0$,
and a multiple factor $Y^2$ appears
only for $\alpha=-1$.
One can check that a finite admissible measure does not exist
when $\Omega$ is unbounded. The case $\varepsilon=-1$ is similar.
\eenum

%%%%%%%%%%%%%%%%%%%%%%%%%%%%%%%%%%%%%%%%%%%%%%%%%%%%%%%%%%%
%%%%%%%%%%%%%%%%%%%%%%%%%%%%%%%%%%%%%%%%%%%%%%%%%%%%%%%%%%%

\item[(\ncX)]
$\partial \Omega= \{ X=0\}$. In this case, $G^{11}$ and $G^{12}$ are multiples of $X$.
The integrability condition implies that $\deg Q\ge 2$ in \eqref{measure.form1}.
Then \eqref{mes.deg} implies the sum of the degrees of distinct irreducible
factors of $\Delta$ is at most $2$. Since $X$ divides $\Delta$,
the other factor (if exists) is linear,
and the non-vanishing of $\Delta$ on $\Omega$ implies that it is a function of $X$.
Thus $\deg_Y\Delta=0$. Then \eqref{mes.deg} implies
$\deg\Delta=\deg_X\Delta\le 3$.
Let us write $G^{11}=(a_1+a_0 Y)X$, $G_{12}=(b_1+b_0 Y)X$,
$G^{22}=c_2+c_1Y+c_0 Y^2$ where $a_k$, $b_k$, $c_k$ are
polynomials in $X$ of degree $\le k$.
Since $G$ is positive definite on $\Omega$, we have $a_1+a_0Y>0$ when $x>0$,
whence $a_0=0$ and $a_1\ne 0$.
Further, $\deg_Y\Delta=0$ implies $Xa_1c_0-X^2b_0^2=0$
(the coefficient of $Y^2$).

\benum
  \item[(\ncX i)] $b_0\ne 0$.
  Then $a_1c_0=Xb_0^2$, hence up to scaling,
  we have $a_1 = X$ and $c_0=b_0^2$.
  By equating the coefficient of $Y$ to zero, we deduce that
  $G^{22}=\tilde c_2 + (b_1+b_0Y)^2$ where $\tilde c_2$ is a polynomial in $X$
  of degree $\le 2$.
  The change $(X,Y)\mapsto(X,Y+pX+q)$
  transforms $b_1$ into $b_1-qb_0+p(1-b_0)X$, thus we may reduce to
  $G^{11}=X^2$, $G^{12}=Xl_b$, $G^{22}=\tilde c_2l_b^2$ where
  $\deg_Y\tilde c_2=0$, and $l_b=Y$ or $l_b=\beta X+Y$.
  Then $\Delta = X^2\tilde c_2(X)$. Since $\deg\Delta\le 3$, we have
  $\deg\tilde c_2\le 1$.
  The case $l_b=\beta Y$ contradicts
  Corollary~\ref{cor.measure}, thus $l_b = \beta X+Y$ with $\beta\ne 0$.
  One checks that there is no integrable measure solution.
  \item[(\ncX  ii)]
     $b_0=0$. Then, since $a_1\ne 0$ and $\deg_Y\Delta=0$,
     we have $c_0=c_1=0$, thus $G$ depends on $X$ only.
     Note that the change of coordinates $(X,Y)\mapsto(X,Y+pX)$
     transforms $b_1$ into $b_1+p a_1$.
     \benum
        \item[(a)] $\deg a_1=1$. Then, up to change of
        coordinates and rescaling, we may assume $a_1=X+\alpha_0$ and
        $b_1=\text{const}$, hence $\deg\Delta\le 3$ implies $\deg c_2\le 1$
        but then Corollary~\ref{cor.measure} implies $\deg\Delta\le 2$,
        i.e., $c_2=\text{const}$. One checks that there is no
        integrable measure solution.
       \item[(b)] $\deg a_1=0$.
        Then, by the aforementioned change of
        coordinates and rescaling, we may achieve that $a_1=1$ and
        $b_1=\beta_1 X$.
        Then the condition $\deg\Delta\le 3$ implies $b_1=0$.
        One easily checks that the only measure solution is the
        product of Laguerre and Hermite polynomials.
     \eenum
\eenum

%%%%%%%%%%%%%%%%%%%%%%%%%%%%%%%%%%%%%%%%%%%%%%%%%%%%%%%%%%%
%%%%%%%%%%%%%%%%%%%%%%%%%%%%%%%%%%%%%%%%%%%%%%%%%%%%%%%%%%%

\item[(\ncXY)] $\partial\Omega\subset \{XY=0\}$.
 The boundary equations  imply $G^{11}$ and $G^{12}$ are multiples
 of $X$ while $G^{12}$ and $G^{22}$ are multiples of $Y$. Hence
 the metric is
$$
   G= \bpm Xl_a& -\beta XY\\ -\beta XY& Y l_c\epm,\qquad l_a,l_c\in\cP_1^2,
$$
thus $\Delta=XY(l_al_c-\beta^2 XY)$.
By Corollary~\ref{cor.measure} we have
$\deg_X\Delta\le 2$ and $\deg_Y\Delta\le 2$,
whence $\deg_X l_al_c\le 1$ and $\deg_Y l_al_c\le 1$.
The integrability condition combined with
\eqref{mes.deg} implies that the sum of the degrees of the irreducible factors of
of $\Delta$ is at most $3$. Therefore, since $l_al_c-\beta XY$ cannot be a square of
a polynomial of degree $1$, it is either affine linear or
divisible by $X$ or by $Y$. The ellipticity also implies that $l_a$ and $l_c$ cannot
be identically zero.
Hence, up to scaling and exchange of $X$ and $Y$,
one of the following cases occurs:

\benum
  \item[(\ncXY i)] $l_a = 1$, $\beta=0$.
  A computation shows that a measure solution exists only when
  $l_c=\text{const}$. This corresponds to the product of Laguerre polynomials.

  \item[(\ncXY ii)]
    $l_a = \alpha+\alpha_1 X$,
    $l_c = \gamma+\gamma_1 Y$, and:
    \benum
      \item[(a)] $\alpha_1\gamma_1-\beta^2=0$.
        If $\beta=0$, then $\alpha_1\gamma_1=0$, thus we fall into Case (\ncXY i).
        So, we may assume that $\alpha_1=1$ and $\gamma_1=\beta^2$. Then
        a long but routine case-by-case consideration shows that
        there is no measure solution.
      \item[(b)] $l_a=X$: no measure solution.
    \eenum
  \item[(\ncXY iii)]
    $l_a = \alpha+\alpha_1 Y$,
    $l_c = \gamma+\gamma_1 X$, and:
    \benum
      \item[(a)] $\alpha_1\gamma_1-\beta^2=0$. As in (\ncXY ii.a), we may
        assume $(\alpha_1,\gamma_1)=(1,\beta^2)$, thus
        $$
            G = \bpm X(\alpha+Y) & -\beta XY                \\
                     -\beta XY      & (\gamma + \beta^2 X)Y \epm,
            \qquad \alpha,\gamma\ge 0,
           \;\; (\alpha,\gamma)\ne(0,0).
        $$
        A measure solution exists when $\beta>0$, $\Omega=(\bR_+)^2$. It is
        $$
           C X^{p-1} Y^{q-1} \exp( -\lambda\beta X - \lambda Y),
           \qquad C,p,q,\lambda > 0. 
        $$
        The curvature is non-constant.
      \item[(b)] $l_a=Y$: no measure solution.
    \eenum
\eenum

\item[(\ncXX)] $\partial\Omega= \{X^2=1\}$. 
Lemma~\ref{lem.ParLines} combined with the condition that $G^{11}>0$ on $\Omega$
(since $\LL$ is elliptic) implies that, up to scaling and change of
coordinates, we may assume $G^{11}=1-X^2$ and $G^{12}=0$. Then
\eqref{eq.rho.1} reads $G^{22}\partial_Y\log\rho\in\cP_1^2$. 
Hence the integrability condition implies that $\deg_Y G^{22}=0$.
We also see from \eqref{mes.deg} that the sum of the degrees of the irreducible
factors of $\Delta$ is $\le 2$, hence $G^{22}=(X-1)^r(X+1)^s$ with $r+s\le 2$.
One easily checks that a measure solution exists only when $G^{22}$ is constant,
and it is a product of Hermite and Jacobi polynomials.
Notice however that the coordinate change $(X,Y)\mapsto (X,Y+\beta X)$ transforms
the (co)metric into
$$
   (1-X^2)\bpm 1&\beta\\\beta&\beta^2\epm + \bpm 0&0\\0&\gamma\epm.
$$

\item[(\ncXXY)] $\partial\Omega\subset\{XY(1-X)=0\}$.
The metric solution is  (up to homothety and affine change)
$$G=\bpm X(1-X)&0\\0& Y(\alpha X+\beta Y+\gamma)\epm.$$  
Except in the case $\alpha=0$, the curvature is non constant, and the additional factor in $\Delta$: $\alpha X+\beta Y+\gamma$ does not satisfy the boundary equation. The only case when there is a measure solution on the domain is $\alpha=\beta=0$, which is a product of Jacobi and Laguerre polynomials.

\item[(8)] $\partial\Omega\subset\{XY-1\}$.
We see from Proposition~\ref{AlgDOP.hyperb} that $\deg\Delta<4$ only when
$(G^{11},G^{12},G^{22}) = (X^2,\,XY-2,\,Y^2)$. One easily checks that
there is no measure solution.
\eenum

\section{Two fold covers, surfaces of revolution, etc.\label{2fold.covering}}

\subsection{ Simple double covers \label{sect.sdc} }
For many examples in dimension $2$, with domain $\Omega$ described by the equation $P(X,Y)\geq 0$, one may look at models in dimension $3$ given by the equation $Z^2\leq P(X,Y)$. 
It turns out that, in every case where no cusp or double tangent appears in $\partial\Omega$, this provides a new domain in dimension $3$  which  is again a solution of the problem. This is therefore the case for the circle, the triangle, the double parabola and  the double point cubic.

Those new three dimensional models present the same pathology than the circle and triangle models in dimension 2: the metric is not in general  unique up to scaling, the curvature is not constant (except for specific values of the parameters). In fact, in those models, the boundary of the domain has degree at most $4$, whereas the maximal degree of the boundary in general  is $6$. The Laplace operator associated with the metric does not in general belong to the admissible operators.

For example, if one starts with the nodal cubic described
in Section  \ref{sec:DlePtCub}, one gets  for the metric, up to scaling,
$$
  G= \bpm 4X(1-X)  &  2Y(2-3X)               &  2Z(2-3X)\\
         2Y(2-3X)  &  4X-3X^2-9Y^2-(9+a)Z^2  &  aYZ\\
         2Z(2-3X)  &  aYZ                    &  4X-3X^2-(9+a)Y^2-9Z^2\epm.
$$
For the double cover of the triangle, however, one gets a unique metric up to scaling, which is
$$
  G = \bpm 4X(1-X) & -4XY      &  2Z(1-3X)          \\
          -4XY     &  4Y(1-Y)  &  2Z(1-3Y)           \\
         2Z(1-3X)  &  2Z(1-3Y) &  X+Y-X^2-XY-Y^2-9Z^2 \epm,
$$
which has no constant curvature.
For the double cover of the square $[-1,1]\times[-1,1]$ we get
$$
  G = \bpm a(1-X^2) & 0        & -aXZ \\
           0        & b(1-Y^2) & -bYZ \\
           -aXZ     & -bYZ     & b(X^2-Z^2-1)+a(Y^2-Z^2-1) \epm.
$$

 We did not try to push the analysis of these models any further, but this shows that one may construct in higher dimension some models which are not direct extensions of the 2 dimensional models, and that the higher dimension analysis of the problem seems much more complex.

\subsection{ Weighted double covers }

One can observe that in all bounded solutions in dimension 2 except four of them,
the domain $\Omega$ is of the form
$y^2-x^r(1-x)^s=0$, or $(ay^2 - x^r)(by^2 - (1-x)^s)=0$.
We may consider the domains in $\bR^{d+1}$ of the form
\beq
   \label{eq.wdc}
   \prod_k \Big(a_k z^2 - \prod_l F_{k,l}(x_1,\dots,x_d)^{p_{k,l}}\Big) = 0.
\eeq
An easy consequence of Theorem~\ref{thm.gal.bord.Omega} is that if such a domain
admits a solution to the SDOP problem, then its intersection with the
hyperplane $z=0$ admits at least a solution of the algebraic counterpart
of the DOP problem (the $\bR$-AlgDOP problem according to the terminology of
Section~\ref{dim2}).

The following is a complete list of all bounded domains in $\bR^3$ of the form
\eqref{eq.wdc}
which admit a solution of the DOP problem.
We do not include
the direct products of plane domains by a segment (in these cases all admissible
metrics are also direct products).
In the angular brackets we indicate the dimension of the set $\cG$ of
admissible metrics (by Theorem~\ref{thm.gal.bord.Omega} it is always an
open cone in a linear subspace of the space of metrics);
``S.R." means ``surface of revolution".
\def\<{\langle}
\def\>{\rangle}
\benum
\item[(1a)] $\<2\>$: $xy(1-x)(1-y)-z^2$;       %$\dim\cG=2$;
\item[(1b)] $\<3\>$: $(xy - z^2)(1-x)(1-y)$ (same as (3j));   %$\dim\cG=3$;
\item[(1c)] $\<2\>$: $(xy^2 - z^2)(1-x)(1-y)$; %$\dim\cG=2$;
\item[(1d)] $\<2\>$: $(xy(1-x) - z^2)(1-y)$;   %$\dim\cG=2$;
\item[(1e)] $\<4\>$: $(xy - z^2)\big((1-x)(1-y) - z^2\big)$ (S.R.);

\medskip

\item[(2a)] $\<12\>$: $1-x^2-y^2-z^2$  (S.R.),     %$\dim\cG=12$;
\item[(2b)] $\< 2\>$: $(1-x^2-y^2)^2-z^2$ (S.R., same as (4j) with $a=c$ and $b=d$);
                  %$\dim\cG=2$;

\medskip

\item[(3a)] $\<1\>$: $xy(1-x-y)-z^2$;     %$\dim\cG=1$;
\item[(3b)] $\<4\>$: $xy( 1-x-y   -z^2)$ (same as (4f)); %$\dim\cG=4$;
\item[(3c)] $\<6\>$: $xy((1-x-y)^2-z^2)$ (tetrahedron); %$\dim\cG=6$;
\item[(3d)] $\<2\>$: $xy((1-x-y)^3-z^2)$ (same as (8c)); %$\dim\cG=2$;
\item[(3e)] $\<7\>$: $(xy-z^2)(1-x-y)$ (S.R.);   %$\dim\cG=7$;
\item[(3f)] $\<1\>$: $(xy^2-z^2)(1-x-y)$; %$\dim\cG=1$;
\item[(3g)] $\<1\>$: $(x-z^2)(y-az^2)(1-x-y)$, $a\ne0$, $a>-1$
            ($a=1$ $\Rightarrow$ same as (4d)); %$\dim\cG=1$;
\item[(3h)] $\<1\>$: $x(y^2-4z^2)(1-x-y+z^2)$;         %$\dim\cG=1$;
\item[(3i)] $\<3\>$: $(xy-z^2)( 1-x-y    +z^2)$ (same as (6d));       %$\dim\cG=3$;
\item[(3j)] $\<3\>$: $(xy-z^2)((1-x-y)^2-4z^2)$ (same as (1b));       %$\dim\cG=3$;
\item[(3k)] $\<1\>$: $(x-az^2)(y-bz^2)(1-x-y-z^2)$, $ab\ne0$, $a+b>-1$; %$\dim\cG=1$;

\medskip

\item[(4a)] $\<1\>$: $(x-y^2)(1-x-ay^2)-z^2$, $a\ne0$, $a>-1$
            ($a=1$ $\Rightarrow$ $\dim\cG=2$);
\item[(4b)] $\<2\>$: $x( 1-x-y^2)   - z^2$;     %$\dim\cG=2$;
\item[(4c)] $\<5\>$: $x( 1-x-y^2    - z^2)$ (S.R.);    %$\dim\cG=5$;
\item[(4d)] $\<1\>$: $x((1-x-y^2)^2 - z^2)$ (same as (3g) with $a=1$); %$\dim\cG=1$;
\item[(4e)] $\<3\>$: $(x  -z^2)(1-x-y^2)$;    %$\dim\cG=3$;
\item[(4f)] $\<4\>$: $(x^2-z^2)(1-x-y^2)$ (same as (3b));    %$\dim\cG=4$;
\item[(4g)] $\<2\>$: $(x^3-z^2)(1-x-y^2)$ (same as (8b));    %$\dim\cG=2$;
\item[(4h)] $\<2\>$: $(x-az^2)(1-x-y^2-z^2)$, $a\ne0$, $a>-1$;  %$\dim\cG=2$;
\item[(4i)] $\<3\>$: $(x^2-4z^2)(1-x-y^2+z^2)$ (same as (6c));    %$\dim\cG=3$;
\item[(4j)] $\<1\>$: $(x-ay^2-cz^2)(1-x-by^2-dz^2)$, $abcd\ne0$, $a+b>0$, $c+d>0$
            \\(if $a=b$ and $c=d$, then $\dim\cG=2$, same as (2b), and S.R.);

\medskip

\item[(5a)] $\<2\>$: $y(x^2-y)(1-x - z^2)$;
\item[(5b)] $\<2\>$: $y(x^2-y)((1-x)^2-z^2)$;
\item[(5c)] $\<1\>$: $y(x^2-y)((1-x)^3-z^2)$;
\item[(5d)] $\<2\>$: $(x^2-y)(1-x)(y-z^2)$;
\item[(5e)] $\<2\>$: $(x^2-y)(1-x)(y^2-z^2)$;
\item[(5f)] $\<2\>$: $y(1-x)(x^2-y-z^2)$;
\item[(5g)] $\<2\>$: $(x^2-y)(y(1-x) - z^2)$;
\item[(5h)] $\<1\>$: $y\big((x^2-y)(1-x) - z^2\big)$;
\item[(5i)] $\<1\>$: $y\big((x^2-y)(1-x)^2 - z^2\big)$;
\item[(5j)] $\<2\>$: $y((1-x)^2 - z^2)(x^2-y -z^2)$;
\item[(5k)] $\<1\>$: $(x^2-y - z^2)((x-1)^2 - z^2)(y + z^2)$;

\medskip

\item[(6a)] $\<2\>$: $(1+y-2x)(1+y+2x - z^2)(x^2-y)$;
\item[(6b)] $\<2\>$: $(1+y-2x)((1+y+2x)^2 - z^2)(x^2-y)$;
\item[(6c)] $\<3\>$: $(1+y-2x)(1+y+2x)(x^2-y-z^2)$ (same as (4i));
\item[(6d)] $\<3\>$: $\big((1+y-2x)(1+y+2x)-z^2\big)(x^2-y)$ (same as (3i));
\item[(6e)] $\<1\>$: $(1+y-2x)(1+y+2x + z^2)(x^2-y-z^2)$;
\item[(6f)] $\<2\>$: $\big((1+y-2x)(1+y+2x) - 4z^2\big)(x^2-y + z^2)$ (S.R.);

\medskip

\item[(7a)] $\<2\>$: $x^2-x^3-y^2-z^2$ (S.R.);

\medskip

\item[(8a)] $\<2\>$: $(x^3-y^2-z^2)(1-x)$ (S.R.);
\item[(8b)] $\<2\>$: $(x^3-y^2)(1-x-z^2)$ (same as (4g));
\item[(8c)] $\<2\>$: $(x^3-y^2)((1-x)^2-z^2)$ (same as (3d));
\item[(8d)] $\<1\>$: $(x^3-y^2)((1-x)^3-z^2)$;

\medskip

\item[(9a)] $\<1\>$: $(x^3-y^2 - z^2)(2y-3x+1)$;
\item[(9b)] $\<2\>$: $(x^3-y^2)(2y-3x+1 - z^2)$;
\item[(9c)] $\<2\>$: $(x^3-y^2)((2y-3x+1)^2 - z^2)$;

\eenum

\subsection{ Surfaces of revolution \label{sect.revol}}

Each bounded two-dimensional solution admitting an axial symmetry
(thus all of them except the cubic with a tangent and the parabola with 
the axis and a tangent)
provides a three-dimensional solution obtained by rotation around the axis of symmetry.
This observation has the following higher-dimensional generalization.

\bprop
Let $\Omega$ be a domain in $\bR^n$ which is symmetric
with respect to the coordinate hyperplanes $x_i=0$, $i=1,\dots,m$,
and let $\widetilde\Omega\in\bR^{d_1}\times\dots\times\bR^{d_m}\times\bR^{n-m}$ be
given by
$$
  \widetilde\Omega=\big\{(\xx_1,\dots,\xx_m,\yy)\;\big|\;
              \big(\|\xx_1\|,\dots,\|\xx_m\|,\yy\big)\in\Omega\big\}.
$$
Then $\widetilde\Omega$ admits a solution of the SDOP problem
if and only if so does $\Omega$. 
\eprop

\bpf
The ``only if" statement easily follows from Theorem~\ref{thm.gal.bord.Omega}.
Let us prove the ``if" statement. By induction, it is enough to do it for $m=1$.
Let $(\Omega,g,\rho)$ be a solution to the SDOP problem.
We may assume that $g$ is invariant under the symmetry because otherwise we
replace $g$ by its sum with its image under the symmetry (here the positive
definiteness of $g$ is crusial).
Then $g^{ij}$ is even (resp. odd) with respect to $x_1$
if $1$ occurs even (resp. odd) number of times in the pair $(i,j)$, i.~e.,
$g^{ij}=h^{ij}(x_1^2,x_2,\dots,x_d)$ when $i=j=1$ or $2\le i\le j$,
and $g^{1j}=x_1 h^{1j}(x_2,\dots,x_d)$ for $j\ge 2$.
To simplify the notation, we assume that $n=3$ and $d_1=2$, and we denote
the coordinates in $\bR^3$ and in $\bR^2\times\bR^2$ by
$(x,y_1,y_2)$ and $(x_1,x_2,y_1,y_2)$ respectively.
The general case is similar.
So, we have
$$
   g(x,y_1,y_2)=\bpm
     h^{11}(x^2,y_1,y_2) & xh^{12}(y_1,y_2)    & xh^{13}(y_1,y_2)    \\
                *        & h^{22}(x^2,y_1,y_2) & h^{23}(x^2,y_1,y_2) \\
                *        &         *           & h^{33}(x^2,y_1,y_2) \epm.
$$
Let us set
$$
   \tilde g(x_1,x_2,y_1,y_2)=\bpm
     h^{11} & 0      & x_1h^{12} & x_1h^{13} \\
        0   & h^{11} & x_2h^{12} & x_2h^{13} \\
        *   &   *    & h^{22}    & h^{23}    \\
        *   &   *    &    *      & h^{33}    \epm.
$$
where the $x^2$ in the arguments of $h^{ij}$ is replaced by $x_1^2+x_2^2$.
Since $\partial\Omega$ is symmetric, its reduced equation is of the form
$F(x^2,y_1,y_2)=0$. By \eqref{eq.thm.gal} we have
$$
\begin{aligned}
   x\Big(2h^{11}\partial_1F
   + h^{12}\partial_2F
   + h^{13}\partial_3F\Big)&= S^1 F\\
   2x^2h^{i1}\partial_1F
   + h^{i2}\partial_2F
   + h^{i3}\partial_3F &= S^i F, \qquad i=2,3,
\end{aligned}
$$
where the arguments of $F$ and $\partial_jF$ are $(x^2,y_1,y_2)$.
The domain $\Omega$ is connected and symmetric with respect to the plane
$x=0$. Hence $\Omega\cap\{x=0\}$ is non-empty. Therefore $x$ cannot divide $F$
because $F$ is a factor of $\det(g)$ and $g$ is non-degenerate on $\Omega$.
It follows that $S^1=cx$ for some constant $c$, hence we can cancel the
both sides of the first equation by $x$. Then the result follows from
Theorem~\ref{thm.gal.bord.Omega} because the left hand sides of
\eqref{eq.thm.gal} for $\tilde g$ are the same as for $g$ except that
the first equation is replaced by two equations with $x_1$ or $x_2$ standing
instead of the factor $x$, and the arguments of $F$ and $\partial_jF$ in
all the equations are $(x_1^2+x_2^2,y_1,y_2)$.
\epf

\bigskip
{\bf Aknowledgment.}  The authors thank an anonymous referee for valuable remarks on a previous version of this paper.
 
\bibliographystyle{amsplain}   
\bibliography{bib}

  \end{document}